\pdfoutput=1
\documentclass[11pt]{article}
\usepackage[left=1in,right=1in,top=1in,bottom=1in]{geometry}
\usepackage{times}
\usepackage{expl3}
\usepackage{cite}
\usepackage[table]{xcolor}
\usepackage{multirow}
\usepackage{stackengine} 
\usepackage{hhline}
\usepackage{mathtools}
\usepackage{tikz}
\usetikzlibrary{matrix, calc}
\usepackage{lipsum}
\usepackage{titlesec}
\usepackage{wrapfig}
\usepackage{enumerate}
\usepackage{epsfig}
\usepackage{amsmath}
\usepackage{tabularx}
\usepackage{array}
\usepackage{booktabs}
\usepackage{enumitem}
\usepackage{bbm}
\usepackage{calc}
\usepackage{graphicx}
\usepackage{amsmath}
\usepackage[title]{appendix}
\usepackage{amssymb}
\usepackage{epstopdf}
\usepackage{boldline}
\usepackage{arydshln}
\usepackage{calligra}
\usepackage{bm}
\usepackage{url}
\usepackage{blindtext}
\usepackage{accents}

\newcommand{\define}{\stackrel{\mbox{\tiny def}}{=}}

\newtheorem{definition}{Definition}
\newtheorem{theorem}{Theorem}

\newtheorem{corollary}{Corollary}
\newtheorem{lemma}{Lemma}

\newtheorem{example}{Example}
\newtheorem{remark}{Remark}
\newcommand{\myop}{\sqcup}
\usepackage{mathtools}
\usepackage{epstopdf}
\usepackage{balance}
\usepackage{thmtools}
\usepackage{thm-restate}
\usepackage{hyperref}
\usepackage{cleveref}
\usepackage[mathscr]{euscript}

\usepackage[ruled,vlined]{algorithm2e}
\include{pythonlisting}
\newcommand{\ostar}{\mathbin{\mathpalette\make@circled\star}}

\makeatletter
\newcommand{\removelatexerror}{\let\@latex@error\@gobble}
\makeatother
\makeatletter
\newcommand*{\rom}[1]{\expandafter\@slowromancap\romannumeral #1@}
\makeatother

\ExplSyntaxOn
\newcommand\latinabbrev[1]{
  \peek_meaning:NTF . {
    #1\@}%
  { \peek_catcode:NTF a {
      #1.\@ }%
    {#1.\@}}}
\ExplSyntaxOff

\titleclass{\subsubsubsection}{straight}[\subsubsection]

\begin{document}
\vspace{1cm}
\title{Matrix Ordering through Spectral and Nilpotent Structures in Totally Ordered Complex Number Fields}
\vspace{1.8cm}
\author{Shih-Yu~Chang
\thanks{Shih Yu Chang is with the Department of Applied Data Science,
San Jose State University, San Jose, CA, U. S. A. (e-mail: {\tt
shihyu.chang@sjsu.edu}). 
           }}

\maketitle

\begin{abstract}
Matrix inequalities play a pivotal role in mathematics, generalizing scalar inequalities and providing insights into linear operator structures. However, the widely used Lowner ordering, which relies on real-valued eigenvalues, is limited to Hermitian matrices, restricting its applicability to non-Hermitian systems increasingly relevant in fields like non-Hermitian physics. To overcome this, we develop a total ordering relation for complex numbers, enabling comparisons of the spectral components of general matrices with complex eigenvalues. Building on this, we introduce the Spectral and Nilpotent Ordering (SNO), a partial order for arbitrary matrices of the same dimensions. We further establish a theoretical framework for majorization ordering with complex-valued functions, which aids in refining SNO and analyzing spectral components. An additional result is the extension of the Schur–Ostrowski criterion to the complex domain. Moreover, we characterize Jordan blocks of matrix functions using a generalized dominance order for nilpotent components, facilitating systematic analysis of non-diagonalizable matrices. Finally, we derive monotonicity and convexity conditions for functions under the SNO framework, laying a new mathematical foundation for advancing matrix analysis.
\end{abstract}


\section{Introduction}\label{sec: Introduction}

Matrix inequalities are essential tools in both pure and applied mathematics due to their ability to generalize scalar inequalities and provide insights into the structure of linear operators. In pure mathematics, matrix inequalities are studied to explore spectral theory, convex optimization, and functional analysis. Inequalities like the L\"owner-Heinz inequality help establish fundamental properties of operator monotone functions \cite{bhatia2013matrix}. Furthermore, they play a crucial role in proving results in matrix analysis, such as majorization theory, which underpins inequalities like Ky Fan's inequality \cite{marshall1979inequalities}. In applied mathematics, matrix inequalities are indispensable in control theory, signal processing, and machine learning. For instance, the Lyapunov inequality provides the stability conditions of dynamic systems \cite{boyd1994linear}. In quantum information theory, matrix inequalities such as those arising from the von Neumann entropy provide bounds on quantum state transformations \cite{nielsen2010quantum}. Moreover, numerical algorithms to solve semidefinite programing problems utilize matrix inequality constraints \cite{vandenberghe1996semidefinite} heavily. The interplay between theory and application demonstrates the profound impact of matrix inequalities across disciplines.

L\"owner ordering, a fundamental concept in matrix analysis, provides a partial order on Hermitian matrices based on their spectra and is widely utilized in quantum information theory, optimization, and operator theory \cite{bhatia2013matrix}. Many existing studies on inequalities involving functions of operators or matrices rely on the framework of L\"owner ordering~\cite{chang2023BiRanTenPartI,chang2022randomMOI,chang2022randomPDT,chang2022randomDTI,chang2024multivariateMP,chang2024generalizedJensen,chang2024generalizedCDJ}. However, its applicability is limited to Hermitian matrices, as it relies on the real-valued nature of eigenvalues to establish the order. This limitation has become increasingly evident with the rise of non-Hermitian physics, a field gaining significant attention for its ability to model phenomena such as open quantum systems, exceptional points, and non-reciprocal wave propagation \cite{bergholtz2021exceptional}. Non-Hermitian matrices often possess complex eigenvalues, which lack the spectral properties required for L\"owner ordering. Consequently, alternative mathematical frameworks, such as pseudospectra theory \cite{trefethen2005spectra}, have emerged to analyze non-Hermitian systems. Recent studies propose generalized orderings and other spectral measures tailored for non-Hermitian operators \cite{gong2018topological,ashida2020non}. These approaches are crucial for capturing the novel physics encoded in the non-Hermitian eigenstructure, which is inaccessible through traditional Hermitian methods.

Spectral analysis reveals that a general matrix can be decomposed into two distinct parts: the spectral part, associated with its eigenvalues, and the nilpotent part, represented by its Jordan blocks. For Hermitian matrices, which are diagonalizable, the nilpotent part vanishes, leaving only the spectral part. This distinction motivates the extension of the L\"owner ordering, traditionally defined for Hermitian matrices based solely on their spectra, to arbitrary matrices by jointly considering both the spectral and nilpotent components. Such an approach would enable a broader and more nuanced framework for matrix comparison, applicable even to non-Hermitian and non-normal matrices \cite{horn1990matrix, bhatia2013matrix}.

A promising direction is to leverage \emph{majorization ordering} for comparing the spectral components and \emph{dominance ordering} for the nilpotent parts. Majorization provides a natural way to compare the eigenvalue distributions of matrices, capturing information about their ``spread" or ``concentration" \cite{marshall1979inequalities}. This is particularly relevant for Hermitian or normal matrices, whose eigenvalues are real or have well-defined complex magnitudes. For the nilpotent part, dominance ordering of Jordan block sizes offers a systematic way to measure the structure and rank effects of the non-diagonalizable component, which is crucial in understanding matrix dynamics in stability analysis and perturbation theory \cite{horn1990matrix, trefethen2005spectra}. By combining these two perspectives, a unified framework emerges, providing insights into the interplay between spectral and nilpotent contributions, with potential applications in operator theory, functional analysis, matrix analysis, quantum information, numerical analysis, and control theory \cite{gong2018topological}.

Our first contribution is the development of a total ordering relation over complex numbers, which enables the comparison of the spectral part of a general matrix as eigenvalues become complex numbers. Building on this, we introduce the spectral and nilpotent ordering (SNO), a partial order designed to compare arbitrary matrices of the same dimensions. The second major contribution focuses on establishing a theoretical framework and deriving properties for majorization ordering with complex-valued functions. These properties play a crucial role in analyzing the spectral part of matrices and in refining SNO. A key advancement lies in characterizing Jordan blocks of matrix functions by leveraging a generalized dominance order for the nilpotent part. This approach facilitates a systematic analysis of the structure of non-diagonalizable matrices. Finally, we identify monotonicity and convexity conditions for functions under the SNO framework, providing a new mathematical foundation for its applications for matrix analysis.

The paper is organized as follows. Section~\ref{sec: Total Order for Complex Numbers} introduces a total ordering relation over complex numbers of the form \( x + y\iota \), where \( \iota = \sqrt{-1} \). This foundational result enables the comparison of complex numbers and serves as a basis for subsequent discussions. In Section~\ref{sec: Matrix Ordering}, we propose the spectral and nilpotent ordering (SNO) framework, which provides a systematic approach to compare finite-dimensional matrices by jointly considering their spectral and nilpotent structures. Section~\ref{sec: Majorization Ordering for Complex-Valued Functions} focuses on majorization ordering for complex-valued functions. Two classes of functions are explored: functions with multiple input variables in Section~\ref{sec: Function with Multiple Variables}, and functions with a single input variable in Section~\ref{sec: Function with a Single Variable}. Section~\ref{sec: Joedan Blocks for Function of a Matrix} investigates the relationship between the nilpotent parts of a matrix \(\bm{X}\) and its transformation under a function \(f(\bm{X})\), characterized via Jordan blocks. Finally, Section~\ref{sec: Monotone and Convex Conditions} establishes conditions for monotonicity and convexity of functions under the SNO framework. Monotonicity is discussed in Section~\ref{sec: Monotonicity}, while convexity is addressed in Section~\ref{sec: Convexity}.

\begin{remark}
All inequalities used to compare two complex numbers in this paper are derived from the newly proposed total ordering for complex numbers introduced in Section~\ref{sec: Total Order for Complex Numbers}.
\end{remark}

\section{Total Ordering for Complex Numbers}\label{sec: Total Order for Complex Numbers}

In this section, we will propose a total ordering over complex numbers $x+y\iota$, where $\iota = \sqrt{-1}$.

\subsection{Cartesian Product of Totally Ordered Sets}

Given two complex numbers  $z_1=x_1 + y_1 \iota$ and   $z_2=x_2 + y_2 \iota$, we define $z_1 \leq z_2$ if we have the following relations:
\begin{eqnarray}\label{eq1: complex num total order}
x_1&<&x_2,
\end{eqnarray}
or
\begin{eqnarray}\label{eq2: complex num total order}
x_1&=&x_2, \mbox{~~but $y_1 \leq y_2$}.
\end{eqnarray}
Note that the zero in the complex numbers will be expressed as $0+0\iota$.

From the order relations given by Eq.~\eqref{eq1: complex num total order} and Eq.~\eqref{eq2: complex num total order}, we have Lemma~\ref{lma: complex num total order} below to show that these relations can provide a total ordering among all complex numbers. 

\begin{lemma}\label{lma: complex num total order}
The order relations defined in Eq.~\eqref{eq1: complex num total order} and Eq.~\eqref{eq2: complex num total order} establish the complex number set as a totally ordered set.
\end{lemma} 
\textbf{Proof:}
From Eq.~\eqref{eq1: complex num total order} and Eq.~\eqref{eq2: complex num total order}, we have the lexicographical ordering $\leq$ on pairs $(x_1, y_1)$ and $(x_2, y_2)$ used to represent the complex numbers $z_1=x_1 + y_1 \iota$ and   $z_2=x_2 + y_2 \iota$, respectively, as:
\[
(x_1, y_1) \leq (x_2, y_2) \iff x_1 < x_2 \text{ or } (x_1 = x_2 \text{ and } y_1 \leq y_2).
\]
To prove this is a total ordering, we verify the following properties:

\subsection*{1. Reflexivity}
For any pair $(x_1, y_1)$, we must show $(x_1, y_1) \leq (x_1, y_1)$.  
Since $x_1 = x_1$, and $y_1 \leq y_1$ (because $\leq$ is reflexive for real numbers), the condition $(x_1 = x_2 \text{ and } y_1 \leq y_2)$ holds. Thus, $(x_1, y_1) \leq (x_1, y_1)$, proving reflexivity.

\subsection*{2. Antisymmetry}
Suppose $(x_1, y_1) \leq (x_2, y_2)$ and $(x_2, y_2) \leq (x_1, y_1)$.  
From the definition of $\leq$, we have two cases:
\begin{itemize}
    \item Case 1: $x_1 < x_2$. This contradicts $(x_2, y_2) \leq (x_1, y_1)$, which requires $x_2 < x_1$ or $x_2 = x_1$ and $y_2 \leq y_1$.
    \item Case 2: $x_1 = x_2$. In this case, $y_1 \leq y_2$ (from $(x_1, y_1) \leq (x_2, y_2)$) and $y_2 \leq y_1$ (from $(x_2, y_2) \leq (x_1, y_1)$). Thus, $y_1 = y_2$.
\end{itemize}
In either case, $(x_1, y_1) = (x_2, y_2)$, proving antisymmetry.

\subsection*{3. Transitivity}
Suppose $(x_1, y_1) \leq (x_2, y_2)$ and $(x_2, y_2) \leq (x_3, y_3)$. We need to show $(x_1, y_1) \leq (x_3, y_3)$.  
From the definition of $\leq$, we have:
\begin{itemize}
    \item If $x_1 < x_2$ and $x_2 < x_3$, then $x_1 < x_3$, so $(x_1, y_1) \leq (x_3, y_3)$.
    \item If $x_1 < x_2$ and $x_2 = x_3$, then $x_1 < x_3$, so $(x_1, y_1) \leq (x_3, y_3)$.
    \item If $x_1 = x_2$ and $x_2 < x_3$, then $x_1 < x_3$, so $(x_1, y_1) \leq (x_3, y_3)$.
    \item If $x_1 = x_2$ and $x_2 = x_3$, then $x_1 = x_3$. Furthermore, $y_1 \leq y_2$ (from $(x_1, y_1) \leq (x_2, y_2)$) and $y_2 \leq y_3$ (from $(x_2, y_2) \leq (x_3, y_3)$), so $y_1 \leq y_3$. Hence, $(x_1, y_1) \leq (x_3, y_3)$.
\end{itemize}
In all cases, transitivity holds.

\subsection*{4. Comparability}
For any two pairs $(x_1, y_1)$ and $(x_2, y_2)$, we must show either $(x_1, y_1) \leq (x_2, y_2)$ or $(x_2, y_2) \leq (x_1, y_1)$.  
From the definition of $\leq$, if $x_1 < x_2$, then $(x_1, y_1) \leq (x_2, y_2)$.  
If $x_1 > x_2$, then $(x_2, y_2) \leq (x_1, y_1)$.  
If $x_1 = x_2$, then $y_1 \leq y_2$ or $y_2 \leq y_1$ (because $\leq$ is a total order for real numbers). Thus, comparability is satisfied.

Since reflexivity, antisymmetry, transitivity, and comparability are satisfied, the order relations defined in Eq.~\eqref{eq1: complex num total order} and Eq.~\eqref{eq2: complex num total order} is a total ordering.
$\hfill\Box$

\subsection{New Basic Inequalities}

\subsubsection{Properties on Total Ordering of Complex Numbers}\label{sec: Properties on Total Ordering of Complex Numbers}

Inequalities are governed by the following properties. These properties also hold if all non-strict inequalities (\(\leq\) and \(\geq\)) are replaced by their corresponding strict inequalities (\(<\) and \(>\)).

\textbf{Converse:} 

The relations \(\leq\) and \(\geq\) are each other's converse. This means that for any complex numbers \(z_1\) and \(z_2\):  
\begin{eqnarray}\label{eq: Converse:} 
z_1 \leq z_2 \iff z_2 \geq z_1.
\end{eqnarray}
 
\textbf{Addition and Subtraction:}

A common constant \(z_3\) may be added to or subtracted from both sides of an inequality. For any real numbers \(z_1\), \(z_2\), and \(z_3\), we have 
\begin{eqnarray}\label{eq: Addition and Subtraction:}
z_1 + z_3 \leq z_2 + z_3 \quad \text{and} \quad z_1 - z_3 \leq z_2 - z_3.
\end{eqnarray}
In other words, the inequality relation is preserved under addition (or subtraction), and the complex numbers form an ordered group under addition.

\textbf{Multiplication and Division:}

If $z_1(=(x_1+y_1\iota))\leq z_2(=(x_2+y_2\iota))$, it is clear that $z_1 z_3 = z_2 z_3$ for any $z_3$ with $x_1=x_2$ and $y_1=y_2$. Below, we want to obtain the conditions for $z_3=(x_3+y_3\iota)$ such that $z_1 z_3 \leq z_2 z_3$. Since we have $z_1(=(x_1+y_1\iota))\leq z_2(=(x_2+y_2\iota))$, there are two cases:
\begin{eqnarray}\label{eq1: Multiplication}
z_1 z_3 \leq z_2 z_3 \Longleftrightarrow
\begin{cases}
x_1 < x_2 & \text{given $\frac{(y_2 - y_1)y_3}{x_2 - x_1}< x_3$,} \\
& ~ \text{or $\frac{(y_2 - y_1)y_3}{x_2 - x_1}=x_3$ with $\frac{(y_1 - y_2)x_3}{x_2 - x_1}\leq y_3$},\\
\mbox{or}\\
x_1 = x_2, y_1 < y_2 & \text{given $y_3 \leq 0$,} \\
& ~ \text{or $y_3 = 0$ with $x_3 \geq 0$}.\\
\end{cases}
\end{eqnarray}
Similarly,  we want to obtain the conditions for $z_3=(x_3+y_3\iota)$ such that $z_1 z_3 > z_2 z_3$. Since we assume $z_1(=(x_1+y_1\iota))\leq z_2(=(x_2+y_2\iota))$, there are two cases:
\begin{eqnarray}\label{eq2: Multiplication}
z_1 z_3 > z_2 z_3 \Longleftrightarrow
\begin{cases}
x_1 < x_2 & \text{given $\frac{(y_2 - y_1)y_3}{x_2 - x_1}> x_3$,} \\
& ~ \text{or $\frac{(y_2 - y_1)y_3}{x_2 - x_1}=x_3$ with $\frac{(y_1 - y_2)x_3}{x_2 - x_1} > y_3$},\\
\mbox{or}\\
x_1 = x_2, y_1 < y_2 & \text{given $y_3 > 0$,} \\
& ~ \text{or $y_3 = 0$ with $x_3 < 0$}.\\
\end{cases}
\end{eqnarray}

For division, we will assume that $z_3 \neq 0 + 0 \iota$. If $z_1(=(x_1+y_1\iota))=z_2(=(x_2+y_2\iota))$, it is clear that $z_1/z_3 = z_2/z_3$ for any $z_3 \neq 0 + 0 \iota$ with $x_1=x_2$ and $y_1=y_2$. Below, we want to obtain the conditions for $z_3=(x_3+y_3\iota)$ such that $z_1/z_3 \leq z_2/z_3$. Since we have $z_1(=(x_1+y_1\iota))\leq z_2(=(x_2+y_2\iota))$, there are two cases:
\begin{eqnarray}\label{eq1: Division}
z_1/z_3 \leq z_2/z_3 \Longleftrightarrow
\begin{cases}
x_1 < x_2 & \text{given $\frac{(y_1 - y_2)y_3}{x_2 - x_1}< x_3$,} \\
& ~ \text{or $\frac{(y_1 - y_2)y_3}{x_2 - x_1}=x_3$ with $\frac{(y_2 - y_1)x_3}{x_2 - x_1}\geq y_3$},\\
\mbox{or}\\
x_1 = x_2, y_1 < y_2 & \text{given $y_3 \geq 0$,} \\
& ~ \text{or $y_3 = 0$ with $x_3 \geq 0$}.\\
\end{cases}
\end{eqnarray}
Similarly,  we want to obtain the conditions for $z_3=(x_3+y_3\iota)$ such that $z_1/ z_3 > z_2/ z_3$. Since we assume $z_1=(x_1+y_1\iota)\leq z_2=(x_2+y_2\iota)$, there are two cases:
\begin{eqnarray}\label{eq2: Division}
z_1/z_3 > z_2/z_3 \Longleftrightarrow
\begin{cases}
x_1 < x_2 & \text{given $\frac{(y_1 - y_2)y_3}{x_2 - x_1}> x_3$,} \\
& ~ \text{or $\frac{(y_1 - y_2)y_3}{x_2 - x_1}=x_3$ with $\frac{(y_2 - y_1)x_3}{x_2 - x_1}< y_3$},\\
\mbox{or}\\
x_1 = x_2, y_1 < y_2 & \text{given $y_3 < 0$,} \\
& ~ \text{or $y_3 = 0$ with $x_3 < 0$}.\\
\end{cases}
\end{eqnarray}

\textbf{Additive Inverse:}

Given two complex numbers $z_1$ and $z_2$ with $z_1=(x_1+y_1\iota)\leq z_2=(x_2+y_2\iota)$, we have
\begin{eqnarray}\label{eq: Additive Inverse:}
-z_2 \leq - z_1. 
\end{eqnarray}

\textbf{Multiplicative Inverse:} 

Given two complex numbers $z_1=(x_1+y_1\iota) \neq 0+0\iota$ and $z_2=(x_2+y_2\iota) \neq 0+0\iota$, we have
\begin{eqnarray}\label{eq1: Multiplicative Inverse:} 
1/z_1 \leq  1/z_2 \Longleftrightarrow
\begin{cases}
\frac{x_1}{x_1^2+ y_1^2} < \frac{x_2}{x_2^2+ y_2^2}, \\
\mbox{or}\\
\frac{x_1}{x_1^2+ y_1^2} = \frac{x_2}{x_2^2+ y_2^2}~\mbox{and $\frac{-y_1}{x_1^2+ y_1^2} \leq \frac{-y_2}{x_2^2+ y_2^2}$.}\\
\end{cases}
\end{eqnarray}
Similarly, we also have
\begin{eqnarray}\label{eq2: Multiplicative Inverse:} 
1/z_1 >  1/z_2 \Longleftrightarrow
\begin{cases}
\frac{x_1}{x_1^2+ y_1^2} > \frac{x_2}{x_2^2+ y_2^2}, \\
\mbox{or}\\
\frac{x_1}{x_1^2+ y_1^2} = \frac{x_2}{x_2^2+ y_2^2}~\mbox{and $\frac{-y_1}{x_1^2+ y_1^2} > \frac{-y_2}{x_2^2+ y_2^2}$.}\\
\end{cases}
\end{eqnarray}

\textbf{Monotonicity:}

A complex-valued function \(f\) is termed \textbf{monotonically non-decreasing} if, for all \(z_1=(x_1+y_1\iota)\) and \(z_2=(x_2+y_2\iota\) such that \(z_1 \leq z_2\), one has:
\[
f(z_1) \leq f(z_2),
\]
which means that \(f\) preserves the order.  Likewise, a function is called \textbf{monotonically non-increasing} if, whenever \(z_1 \leq z_2\), it holds that:
\[
f(z_1) \geq f(z_2),
\]
which means that \(f\) reverses the order.

If the ordering \(\leq\) in the definition of monotonicity is replaced by the strict ordering \(<\), one obtains a stronger condition. A function satisfying this property is called \textbf{strictly increasing} (also simply \textit{increasing}). That is, \(f\) is strictly increasing if, for all \(z_1 < z_2\):
\[
f(z_1) < f(z_2).
\]
Similarly, by reversing the inequality symbol, one defines the concept of a \textbf{strictly decreasing} function (also called \textit{decreasing}), where for all \(z_1 < z_2\):
\[
f(z_1) > f(z_2).
\]

A function that is either strictly increasing or strictly decreasing is called \textbf{strictly monotone}.  

Functions that are strictly monotone are \textbf{one-to-one} (injective). This is because for any \(z_1 \neq z_2\), either \(z_1 < z_2\) or \(z_1 > z_2\). By monotonicity, it follows that either:
\[
f(z_1) < f(z_2) \quad \text{or} \quad f(z_1) > f(z_2),
\]
which implies:
\[
f(z_1) \neq f(z_2).
\]

\textbf{Convexity:}

Before defining convexity for complex-valued functions of complex variables, we have to define the notion of \emph{Convexity of the Domain}. A set \(D \subset \mathbb{C}\) is convex if for any two points \(z_1, z_2 \in D\) and \(t \in [0, 1]\),
\[
(1 - t)z_1 + t z_2 \in D.
\]
This ensures that the straight line connecting \(z_1\) and \(z_2\) lies entirely within \(D\).

The convexity of a function \(f: D \to \mathbb{C}\), where \(D\) is a convex subset of \(\mathbb{C}\), can be defined as:
\[
f((1 - t)z_1 + t z_2) \leq (1 - t)f(z_1) + t f(z_2),
\]
for all \(z_1, z_2 \in D\) and \(t \in [0, 1]\).
This is similar to the convexity definition for real-valued functions but applied to complex-valued function with the proposed total ordering over complex numbers. 

\subsubsection{Ordered Field of Complex Numbers}\label{sec: Ordered Field Relation}

Recall that an \textbf{ordered field} is a field \(F\) equipped with a total ordering \(\leq\) that satisfies the following properties:

\begin{enumerate}
    \item \textbf{Compatibility with Addition}:  
    For all \(a, b, c \in F\), if \(a \leq b\), then:
    \[
    a + c \leq b + c.
    \]

    \item \textbf{Compatibility with Multiplication}:  
    For all \(a, b \in F\), if \(0 \leq a\) and \(0 \leq b\), then:
    \[
    0 \leq a \cdot b.
    \]

    \item \textbf{Totality of the Order}:  
    The ordering \(\leq\) is a total ordering, meaning for all \(a, b \in F\), exactly one of the following holds:
    \[
    a < b, \quad a = b, \quad \text{or} \quad a > b.
    \]
\end{enumerate}

Now, the field $F$ is the complex number $\mathbb{C}$, the property of  \textbf{Compatibility with Multiplication} will be modified as the following:

Given two complex numbers $0+0\iota \leq z_1=(x_1+y_1\iota)$ and $0+0\iota \leq z_2=(x_2+y_2\iota)$, we have
\begin{eqnarray}\label{eq1: ordered field}
0+0\iota \leq z_1 z_2 \Longleftrightarrow
\begin{cases}
x_1 x_2 - y_1 y_2 \geq 0, \\
\mbox{or}\\
x_1 x_2 = y_1 y_2, \mbox{~with $x_1y_2+x_2y_1 \geq 0$;}\\
\end{cases}
\end{eqnarray}
and
\begin{eqnarray}\label{eq2: ordered field}
z_1 z_2 <  0+0\iota \Longleftrightarrow
\begin{cases}
x_1 x_2 - y_1 y_2 < 0, \\
\mbox{or}\\
x_1 x_2 = y_1 y_2, \mbox{~with $x_1y_2+x_2y_1 < 0$;}\\
\end{cases}
\end{eqnarray}

\section{Matrix Ordering}\label{sec: Matrix Ordering}

\subsection{Majorization Ordering of Complex-Valued Vectors}\label{sec: Majorization Ordering of Complex-Valued Vectors}

The majorization relation is a partial order that compares two vectors based on their sorted components. It is a key concept in mathematics, particularly in optimization, inequalities, and matrix analysis. Here, we will extend the majorization from real numbers to complex numbers with total order relation. 

Let \(\mathbf{x} = (x_1, x_2, \ldots, x_n)\) and \(\mathbf{y} = (y_1, y_2, \ldots, y_n)\) be two vectors in \(\mathbb{C}^n\). Assume that their components are sorted in non-increasing order:
\[
x_1 \geq x_2 \geq \cdots \geq x_n \quad \text{and} \quad y_1 \geq y_2 \geq \cdots \geq y_n.
\]

We say that \(\mathbf{x}\) is majorized by \(\mathbf{y}\), written as \(\mathbf{x} \prec \mathbf{y}\), if the following two conditions hold:\\

1. Partial sums condition:
   \[
   \sum_{i=1}^k x_i \leq \sum_{i=1}^k y_i \quad \text{for all } k = 1, 2, \ldots, n-1,
   \]
   meaning the cumulative sums of the largest \(k\) components of \(\mathbf{x}\) do not exceed those of \(\mathbf{y}\).
   
2. Total sum condition:
   \[
   \sum_{i=1}^n x_i = \sum_{i=1}^n y_i,
   \]
   ensuring the total sum of the components of both vectors is the same.

However, we say that \(\mathbf{x}\) is majorized by \(\mathbf{y}\) weakly, written as \(\mathbf{x} \prec_{w} \mathbf{y}\), if only the partial sum condition is satisifed. 

For example, let \(\mathbf{x} = (4,1+\iota, 3)\) and \(\mathbf{y} = (2+\iota, 5, 1)\), we first sort components: \(\mathbf{x} = (4, 3, 1+\iota)\) and \(\mathbf{y} = (5, 2+\iota, 1)\). Then, we check partial sums:\\
\begin{eqnarray}
4 &\leq&5, \nonumber \\
4 + 3 = 7& \leq & 5 + 2+\iota = 7+\iota \nonumber .
\end{eqnarray}
Finally, we find that total sums: \(4 + 3 + 1+\iota = 8  + \iota= 5 + 2+\iota + 1\). Thus, \(\mathbf{x} \prec \mathbf{y}\).

\subsection{Generalized Dominance Ordering of Partitions}\label{sec: Generalized Dominance Ordering of Partitions}

The dominance ordering (or majorization order) among partitions of a natural number compares partitions based on how ``spread out" their parts are, using partial sums of their components. It is a partial order that is especially important in combinatorics, representation theory, and the study of symmetric functions. Note that we extend conventional dominance order from majorization ordering, where the original two natural numbers under partition are identical, to weak majorization order between any two natural numbers under partition.

Let \( m,n \) be two natural numbers, and let \( \bm{p}= (p_1, p_2, \ldots, p_k) \) and \( \bm{q} = (q_1, q_2, \ldots, q_l) \) be two partitions of $m$ and $n$, respectively, written in non-increasing order:
\[
p_1 \geq p_2 \geq \cdots \geq p_k > 0, \quad q_1 \geq q_2 \geq \cdots \geq q_l > 0.
\]

We say that $\bm{p}$ precedes $\bm{q}$ in the dominance order, written $\bm{p}\trianglelefteq\bm{q}$, if the following conditions are satisfied:\\

1. Partial sums condition:
\begin{eqnarray}\label{eq1: dom ordering}
   \sum_{i=1}^j p_i \leq \sum_{i=1}^j q_i \quad \text{for all } j \geq 1,
\end{eqnarray}
where parts beyond the length of a partition are treated as \(0\) (e.g., if \(\bm{p}\) has 3 parts but \(\bm{q}\) has 4 parts, then \(p_4 = 0\)). \\

2. Total sum condition:
   \[
   \sum_{i=1}^k p_i = m, \quad \sum_{i=1}^l q_i = n,
   \]
   ensuring both are partitions of numbers $m$ and $n$. 

Consider the following example, given $m = 5, n=6$ and the partitions:
\[
\bm{p}= (3, 2), \quad \bm{q} = (4, 2).
\]
we first check partial sums: 
\begin{eqnarray}
   \mbox{For \(\bm{p}\): \(3, 3+2=5\)},\nonumber \\
   \mbox{For \(\bm{q}\): \(4, 4+2=6\)}.\nonumber 
\end{eqnarray}
Then, we also check dominance:
\begin{eqnarray}
   \mbox{First partial sum: \(3 \leq 4\)},\nonumber \\
   \mbox{Second partial sum: \(5 \leq 6\)}.\nonumber 
\end{eqnarray}
Thus, we have $\bm{p}\trianglelefteq\bm{q}$.

Besides, if two partitions with $\bm{p}\trianglelefteq\bm{q}$, their difference distance upto the $j$-th position, denoted by $\mathfrak{D}_{\bm{p},\bm{q}}(j)$, which is defined by
\begin{eqnarray}\label{eq: dom order distance}
\mathfrak{D}_{\bm{p},\bm{q}}(j)\define \sum_{i=1}^j q_i - \sum_{i=1}^j p_i,
\end{eqnarray}
where $\sum_{i=1}^j q_i=n$ and $\sum_{i=1}^j p_i=m$.  We call the quantity $\mathfrak{D}_{\bm{p},\bm{q}}(j)$ as the \emph{generalized dominance ordering distance} (GDOD) between the vector $\bm{p}$ and $\bm{q}$ upto the $j$-th position. 

\subsection{Matrix Ordering via Spectral and Nilpotent Structures}\label{sec: Matrix Ordering via Spectral and Nilpotent Structures}

From Jordan decomposition theorem, a square matrix $\bm{X} \in \mathbb{C}^{m \times m}$ can be decomposed as~\cite{gohberg1996simple,ashida2020non,chang2024operatorCH}:  
\begin{eqnarray}\label{eq: Jordan decompostion}
\bm{X}&=& \bm{U}\left(\bigoplus\limits_{k=1}^{K}\bigoplus\limits_{i=1}^{\alpha_k^{(\mathrm{G})}}\bm{J}_{m_{k,i}}(\lambda_k)\right)\bm{U}^{-1}.
\end{eqnarray}
where $\bm{U} \in \mathbb{C}^{m \times m}$ is an invertible matrix, $\alpha_k^{(\mathrm{G})}$ is the geometry multiplicity with respect to the $k$-th eigenvalue $\lambda_k$, and $\bm{J}_{m_{k,i}}(\lambda_k)$ is a Jordan Block with dimension $m_{k,i} \times m_{k,i}$ and the eigenvalue $\lambda_k$ expressed by  
\begin{eqnarray}\label{eq: Jordan block def}
\bm{J}_{m_{k,i}}(\lambda)= \begin{pmatrix}
   \lambda & 1 & 0 & \cdots & 0 \\
   0 & \lambda & 1 & \cdots & 0 \\
   \vdots & \vdots & \ddots & \ddots & \vdots \\
   0 & 0 & \cdots & \lambda & 1 \\
   0 & 0 & \cdots & 0 & \lambda
   \end{pmatrix}_{m_{k,i} \times m_{k,i}}
\end{eqnarray}
where $\lambda_k$ is a scalar (an eigenvalue of the original matrix), and the off-diagonal elements are either 0 or 1. The representation given by Eq.~\eqref{eq: Jordan decompostion} is unique if we further assume that $\lambda_1 \geq \lambda_2 \geq \ldots \geq \lambda_K$ based on the complex number total ordering, and $m_{k,1} \geq m_{k,2} \geq \ldots \geq m_{k,\alpha_k^{(\mathrm{G})}}$ for all $k=1,2,\ldots,K$ (partition ordering).

We also have the following relationships about $\alpha_k^{(\mathrm{A})}$, the algebraic multiplicity for the $k$-th eigenvalue,  and $\alpha_k^{(\mathrm{G})}$:
\begin{eqnarray}
\sum\limits_{k=1}^K \alpha_k^{(\mathrm{A})}&=&m,
\end{eqnarray}
and 
\begin{eqnarray}
\sum\limits_{i=1}^{\alpha_k^{(\mathrm{G})}}m_{k,i}&=&\alpha_k^{(\mathrm{A})}.
\end{eqnarray}

From Jordan decomposition given by Eq.~\eqref{eq: Jordan decompostion}, we can represent the matrix $\bm{X}$ via its spectral and nilpotent structures, denoted as $\mathfrak{R}(\bm{X})$, by
\begin{eqnarray}\label{eq: matrix SN representation1}
\mathfrak{R}(\bm{X})&\define&[\underbrace{\lambda_1,\ldots,\lambda_1}_{\alpha_1^{(\mathrm{A})}},\underbrace{\lambda_2,\ldots,\lambda_2}_{\alpha_2^{(\mathrm{A})}},\ldots,\underbrace{\lambda_K,\ldots,\lambda_K}_{\alpha_K^{(\mathrm{A})}},\nonumber \\
&&\underbrace{m_{1,1},\ldots,m_{1,\alpha_{1}^{(\mathrm{G})}}}_{\sum\limits_{i=1}^{\alpha_1^{(\mathrm{G})}}m_{1,i}=\alpha_1^{(\mathrm{A})}},\underbrace{m_{2,1},\ldots,m_{2,\alpha_{2}^{(\mathrm{G})}}}_{\sum\limits_{i=1}^{\alpha_2^{(\mathrm{G})}}m_{2,i}=\alpha_2^{(\mathrm{A})}},\ldots,
\underbrace{m_{K,1},\ldots,m_{K,\alpha_{K}^{(\mathrm{G})}}}_{\sum\limits_{i=1}^{\alpha_K^{(\mathrm{G})}}m_{K,i}=\alpha_K^{(\mathrm{A})}}]
\end{eqnarray}
We call such representation as spectral and nilpotent structures representation of the matrix $\bm{X}$ and $\mathfrak{R}(\bm{X})$ is unique as we require that $\lambda_1 \geq \lambda_2 \geq \ldots \geq \lambda_K$ based on the complex number total ordering, and $m_{k,1} \geq m_{k,2} \geq \ldots \geq m_{k,\alpha_k^{(\mathrm{G})}}$ for all $k=1,2,\ldots,K$. Note that the upper row in Eq.~\eqref{eq: matrix SN representation1} represents the spectral structure, i.e., eigenvalues, and the lower row in Eq.~\eqref{eq: matrix SN representation1} represents the nilpotent structure, i.e., nilpotent dimensions information. 

However, in order to distinguish representation $\mathfrak{R}(\bm{X})$ by different input arguments $\bm{X}$, we can rewrite Eq.~\eqref{eq: matrix SN representation1} by adding argument $\bm{X}$ with respect to those eigenvalues $\lambda_k$, multiplicities $m_{k,i}$, $\alpha_k^{(\mathrm{A})}$ and $\alpha_k^{(\mathrm{G})}$, as 
\begin{eqnarray}\label{eq: matrix SN representation2}
\mathfrak{R}(\bm{X})&=&[\underbrace{\lambda_1(\bm{X}),\ldots,\lambda_1(\bm{X})}_{\alpha_1^{(\mathrm{A})}(\bm{X})},\underbrace{\lambda_2(\bm{X}),\ldots,\lambda_2(\bm{X})}_{\alpha_2^{(\mathrm{A})}(\bm{X})},\ldots,\underbrace{\lambda_K(\bm{X}),\ldots,\lambda_K(\bm{X})}_{\alpha_K^{(\mathrm{A})}(\bm{X})},\nonumber \\
&&\underbrace{m_{1,1}(\bm{X}),\ldots,m_{1,\alpha_{1}^{(\mathrm{G})}(\bm{X})}(\bm{X})}_{\sum\limits_{i=1}^{\alpha_1^{(\mathrm{G})}(\bm{X})}m_{1,i}(\bm{X})=\alpha_1^{(\mathrm{A})}(\bm{X})},\underbrace{m_{2,1}(\bm{X}),\ldots,m_{2,\alpha_{2}^{(\mathrm{G})}(\bm{X})}(\bm{X})}_{\sum\limits_{i=1}^{\alpha_2^{(\mathrm{G})}(\bm{X})}m_{2,i}(\bm{X})=\alpha_2^{(\mathrm{A})}(\bm{X})},\ldots,\nonumber \\
&&\underbrace{m_{K,1}(\bm{X}),\ldots,m_{K,\alpha_{K}^{(\mathrm{G})}(\bm{X})}(\bm{X})}_{\sum\limits_{i=1}^{\alpha_K^{(\mathrm{G})}(\bm{X})}m_{K,i}(\bm{X})=\alpha_K^{(\mathrm{A})}(\bm{X})}] \nonumber \\
&\define_1&[\bm{\lambda_1}(\bm{X}),\ldots,\bm{\lambda_K}(\bm{X}),\bm{m}_1(\bm{X}),\ldots,\bm{m}_K(\bm{X})],
\end{eqnarray}
where we set $[\underbrace{\lambda_k(\bm{X}),\ldots,\lambda_k(\bm{X})}_{\alpha_k^{(\mathrm{A})}(\bm{X})}]\define\bm{\lambda_k}(\bm{X})$ and $[\underbrace{m_{k,1}(\bm{X}),\ldots,m_{k,\alpha_{k}^{(\mathrm{G})}(\bm{X})}(\bm{X})}_{\sum\limits_{i=1}^{\alpha_k^{(\mathrm{G})}(\bm{X})}m_{k,i}(\bm{X})=\alpha_k^{(\mathrm{A})}(\bm{X})}]\define\bm{m}_k(\bm{X})$ for $k=1,2,\ldots,K$ in $\define_1$.

Given two matrices $\bm{X}_1$ and $\bm{X}_2$ with the same dimensions $m \times m$ and Jordan decomposition, these two matrices $\bm{X}_1$ and $\bm{X}_2$ can be expressed by
\begin{eqnarray}\label{eq: matrices setup}
\bm{X}_1&=& \bm{U}_1\left(\bigoplus\limits_{k_1=1}^{K_1}\bigoplus\limits_{i_1=1}^{\alpha_{k_1}^{(\mathrm{G})}}\bm{J}_{m_{k_1,i_1}}(\lambda_{k_1})\right)\bm{U}_1^{-1}, \nonumber \\
\bm{X}_2&=& \bm{U}_2\left(\bigoplus\limits_{k_2=1}^{K_2}\bigoplus\limits_{i_2=1}^{\alpha_{k_2}^{(\mathrm{G})}}\bm{J}_{m_{k_2,i_2}}(\lambda_{k_2})\right)\bm{U}_2^{-1}.
\end{eqnarray}   

We have the following definition to compare two matrices $\bm{X}_1$ and $\bm{X}_2$ by their representations $\mathfrak{R}(\bm{X}_1)$ and $\mathfrak{R}(\bm{X}_2)$.
\begin{definition}\label{def: SNO}
From Eq.~\eqref{eq: matrix SN representation2}, we can represent matrices $\bm{X}_1$ and $\bm{X}_2$ as 
\begin{eqnarray}\label{eq1: def: SNO}
\mathfrak{R}(\bm{X}_1)&=&[\bm{\lambda_1}(\bm{X}_1),\ldots,\bm{\lambda}_{K_1}(\bm{X}_1),\bm{m}_1(\bm{X}_1),\ldots,\bm{m}_{K_1}(\bm{X}_1)], \nonumber \\
\mathfrak{R}(\bm{X}_2)&=&[\bm{\lambda_1}(\bm{X}_2),\ldots,\bm{\lambda}_{K_2}(\bm{X}_2),\bm{m}_1(\bm{X}_2),\ldots,\bm{m}_{K_2}(\bm{X}_2)].
\end{eqnarray}   
We say that $[\bm{m}_1(\bm{X}_1),\ldots,\bm{m}_{K_1}(\bm{X}_1)] \preceq_{\mbox{\tiny N}} [\bm{m}_1(\bm{X}_2),\ldots,\bm{m}_{K_2}(\bm{X}_2)]$, where $\preceq_{\mbox{\tiny N}}$ represents the order between $\max(K_1, K_2)$ distinct eigenvalues by considering their nilpotent structures~\footnote{If $K_1 \neq K_2$, the shorter list will be added $|K_1 - K_2|$ zero entries before applying $\trianglelefteq$ comparison.}, if and only if 
there is some $k \in \{1,2,\ldots, \max(K_1, K_2)\}$ such that $\bm{m}_j(\bm{X}_1)=\bm{m}_j(\bm{X}_2)$ for all $j < k$ and $\bm{m}_k(\bm{X}_1)\trianglelefteq\bm{m}_k(\bm{X}_2)$.

We say that $\bm{X}_1 \preceq_{\mbox{\tiny SN}} \bm{X}_2$, where $\preceq_{\mbox{\tiny SN}}$ represents the order between matrices by considering their spectral and nilpotent structures, if and only if 
\begin{eqnarray}\label{eq3: def: SNO}
[\bm{\lambda_1}(\bm{X}_1),\ldots,\bm{\lambda}_{K_1}(\bm{X}_1)] \preceq_w  
[\bm{\lambda_1}(\bm{X}_2),\ldots,\bm{\lambda}_{K_2}(\bm{X}_2)] 
\end{eqnarray}   
or
\begin{eqnarray}\label{eq4: def: SNO}
[\bm{\lambda_1}(\bm{X}_1),\ldots,\bm{\lambda}_{K_1}(\bm{X}_1)] &=&  
[\bm{\lambda_1}(\bm{X}_2),\ldots,\bm{\lambda}_{K_2}(\bm{X}_2)] \nonumber \\
\mbox{with} & & [\bm{m}_1(\bm{X}_1),\ldots,\bm{m}_{K_1}(\bm{X}_1)] \preceq_{\mbox{\tiny N}} [\bm{m}_1(\bm{X}_2),\ldots,\bm{m}_{K_2}(\bm{X}_2)].
\end{eqnarray}   

Moreover, 
We say that $\bm{X}_1 \prec_{\mbox{\tiny SN}} \bm{X}_2$, where $\prec_{\mbox{\tiny SN}}$ represents the order between matrices by considering their spectral and nilpotent structures, if and only if 
\begin{eqnarray}\label{eq5: def: SNO}
[\bm{\lambda_1}(\bm{X}_1),\ldots,\bm{\lambda}_{K_1}(\bm{X}_1)] \prec_w  
[\bm{\lambda_1}(\bm{X}_2),\ldots,\bm{\lambda}_{K_2}(\bm{X}_2)] 
\end{eqnarray}   
or
\begin{eqnarray}\label{eq6: def: SNO}
[\bm{\lambda_1}(\bm{X}_1),\ldots,\bm{\lambda}_{K_1}(\bm{X}_1)] &=&  
[\bm{\lambda_1}(\bm{X}_2),\ldots,\bm{\lambda}_{K_2}(\bm{X}_2)] \nonumber \\
\mbox{with} & & [\bm{m}_1(\bm{X}_1),\ldots,\bm{m}_{K_1}(\bm{X}_1)] \prec_{\mbox{\tiny N}} [\bm{m}_1(\bm{X}_2),\ldots,\bm{m}_{K_2}(\bm{X}_2)].
\end{eqnarray}   
\end{definition}
The order between matrices provided by Definition~\ref{def: SNO} about $\preceq_{\mbox{\tiny SN}}$ is named as \emph{spectral and nilpotent structures ordering}, abbreviated by \textbf{SNO}. 

We have to prove the following Lemma~\ref{lma: weak maj is a po} and Lemma~\ref{lma: preceq N is po} about the partial order relation for $\preceq_{\mbox{\tiny N}}$ before proving that SNO is a partial ordering. 

\begin{lemma}\label{lma: weak maj is a po}
Both $\preceq_w$ and $\trianglelefteq$ are partial orderings.
\end{lemma}
\textbf{Proof:}
For the part of $\preceq_w$, let \(\mathbf{x} = (x_1, x_2, \ldots, x_n)\) and \(\mathbf{y} = (y_1, y_2, \ldots, y_n)\) be two complex vectors of length \(n\), sorted in non-increasing order:
\[
x_1 \geq x_2 \geq \cdots \geq x_n, \quad y_1 \geq y_2 \geq \cdots \geq y_n.
\]
We say that \(\mathbf{x}\) weakly majorizes \(\mathbf{y}\), written as \(\mathbf{x} \prec_w \mathbf{y}\), if:
\[
\sum_{i=1}^k x_i \leq \sum_{i=1}^k y_i \quad \text{for all } k = 1, 2, \ldots, n.
\]

To prove such weak majorization is a partial order on the set of complex vectors of fixed length \(n\), we must demonstrate that it satisfies the three defining properties of a partial order: reflexivity, antisymmetry, and transitivity.

1. Reflexivity: \\
This is to show \(\mathbf{x} \prec_w \mathbf{x}\) for any vector \(\mathbf{x}\). Since for any \(k \in \{1, 2, \ldots, n\}\), we have \(\sum_{i=1}^k x_i = \sum_{i=1}^k x_i\). Thus, \(\mathbf{x} \prec_w \mathbf{x}\) holds.

2. Antisymmetry:\\
To show that if \(\mathbf{x} \prec_w \mathbf{y}\) and \(\mathbf{y} \prec_w \mathbf{x}\), then \(\mathbf{x} = \mathbf{y}\). From \(\mathbf{x} \prec_w \mathbf{y}\), we have:
  \[
  \sum_{i=1}^k x_i \leq \sum_{i=1}^k y_i \quad \text{for all } k = 1, 2, \ldots, n.
  \]
From \(\mathbf{y} \prec_w \mathbf{x}\), we have:
  \[
  \sum_{i=1}^k y_i \leq \sum_{i=1}^k x_i \quad \text{for all } k = 1, 2, \ldots, n.
  \]
Combining the two conditions:
  \[
  \sum_{i=1}^k x_i = \sum_{i=1}^k y_i \quad \text{for all } k = 1, 2, \ldots, n.
  \]
Since the cumulative sums are equal for all \(k\), and the components of \(\mathbf{x}\) and \(\mathbf{y}\) are sorted in non-increasing ordering, this implies \(\mathbf{x} = \mathbf{y}\).

3. Transitivity:\\
To show that if \(\mathbf{x} \prec_w \mathbf{y}\) and \(\mathbf{y} \prec_w \mathbf{z}\), then \(\mathbf{x} \prec_w \mathbf{z}\). From \(\mathbf{x} \prec_w \mathbf{y}\), we know:
  \[
  \sum_{i=1}^k x_i \leq \sum_{i=1}^k y_i \quad \text{for all } k = 1, 2, \ldots, n.
  \]
From \(\mathbf{y} \prec_w \mathbf{z}\), we know:
  \[
  \sum_{i=1}^k y_i \leq \sum_{i=1}^k z_i \quad \text{for all } k = 1, 2, \ldots, n.
  \]
Combining these inequalities:
  \[
  \sum_{i=1}^k x_i \leq \sum_{i=1}^k y_i \leq \sum_{i=1}^k z_i \quad \text{for all } k = 1, 2, \ldots, n.
  \]
Thus, \(\mathbf{x} \prec_w \mathbf{z}\).

To prove that the relation \(\trianglelefteq\) defined above is a partial order relation, we also need to verify that it satisfies the following properties of a partial order: reflexivity, antisymmetry, and transitivity. We now verify each property.

1. Reflexivity: 

For any partition \(\bm{p} = (p_1, p_2, \ldots, p_k)\) of \(m\):  
- The partial sums condition in \(\eqref{eq1: dom ordering}\) holds trivially because:
  \[
  \sum_{i=1}^j p_i \leq \sum_{i=1}^j p_i \quad \text{for all } j \geq 1.
  \]
The total sum condition also holds because \(\bm{p}\) is a partition of \(m\). Thus, \(\bm{p} \trianglelefteq \bm{p}\), so the relation is reflexive.

2. Antisymmetry: 

Assume \(\bm{p} \trianglelefteq \bm{q}\) and \(\bm{q} \trianglelefteq \bm{p}\), where \(\bm{p} = (p_1, p_2, \ldots, p_k)\) and \(\bm{q} = (q_1, q_2, \ldots, q_l)\).  
- From the partial sums condition \(\sum_{i=1}^j p_i \leq \sum_{i=1}^j q_i\) (for all \(j\)) and \(\sum_{i=1}^j q_i \leq \sum_{i=1}^j p_i\), it follows that:
  \[
  \sum_{i=1}^j p_i = \sum_{i=1}^j q_i \quad \text{for all } j \geq 1.
  \]
Since both partitions are in non-increasing order, equality in all partial sums implies that the sequences \((p_1, p_2, \ldots, p_k)\) and \((q_1, q_2, \ldots, q_l)\) must be identical. Thus, \(\bm{p} = \bm{q}\), and the relation is antisymmetric.

3. Transitivity:

Assume \(\bm{p} \trianglelefteq \bm{q}\) and \(\bm{q} \trianglelefteq \bm{r}\), where \(\bm{p} = (p_1, p_2, \ldots, p_k)\), \(\bm{q} = (q_1, q_2, \ldots, q_l)\), and \(\bm{r} = (r_1, r_2, \ldots, r_m)\).  
- From \(\bm{p} \trianglelefteq \bm{q}\), we know:
  \[
  \sum_{i=1}^j p_i \leq \sum_{i=1}^j q_i \quad \text{for all } j \geq 1.
  \]
- From \(\bm{q} \trianglelefteq \bm{r}\), we know:
  \[
  \sum_{i=1}^j q_i \leq \sum_{i=1}^j r_i \quad \text{for all } j \geq 1.
  \]
- Combining these, we get:
  \[
  \sum_{i=1}^j p_i \leq \sum_{i=1}^j r_i \quad \text{for all } j \geq 1.
  \]
Thus, \(\bm{p} \trianglelefteq \bm{r}\), and the relation is transitive.

This lemma is proved as both $\prec_w$ and $\trianglelefteq$ are partial ordering relations. 
$\hfill\Box$

The next lemma will show that $\preceq_{\mbox{\tiny N}}$ is a partial orderinig.

\begin{lemma}\label{lma: preceq N is po}
The order relation given below:
\begin{eqnarray}\label{eq1: lma: preceq N is po}
[\bm{m}_1(\bm{X}_1),\ldots,\bm{m}_{K_1}(\bm{X}_1)] \preceq_{\mbox{\tiny N}} [\bm{m}_1(\bm{X}_2),\ldots,\bm{m}_{K_2}(\bm{X}_2)],
\end{eqnarray}
is a partial ordering.
\end{lemma}
\textbf{Proof:}
Without lost of generality, we may assume that $K_1= K_2$ since we can add $|K_1 - K_2|$ zero entries before applying $\trianglelefteq$ comparison to the shorter list of $[\bm{m}_1(\bm{X}_1),\ldots,\bm{m}_{K_1}(\bm{X}_1)]$or $ [\bm{m}_1(\bm{X}_2),\ldots,\bm{m}_{K_2}(\bm{X}_2)]$.

To prove such $\preceq_{\mbox{\tiny N}}$ order is a partial order on the set of Cartesian product of $K_1$ partition numbers, we must demonstrate that it satisfies the three defining properties of a partial order: reflexivity, antisymmetry, and transitivity.

1. Reflexivity: \\
Because $\bm{m}_k(\bm{X}_1) \trianglelefteq \bm{m}_k(\bm{X}_1)$ for any $k= 1,2,\ldots,K_1$, we have \\
$[\bm{m}_1(\bm{X}_1),\ldots,\bm{m}_{K_1}(\bm{X}_1)] \preceq_{\mbox{\tiny N}} [\bm{m}_1(\bm{X}_1),\ldots,\bm{m}_{K_1}(\bm{X}_1)]$.

2. Antisymmetry:\\
If we have $[\bm{m}_1(\bm{X}_1),\ldots,\bm{m}_{K_1}(\bm{X}_1)] \preceq_{\mbox{\tiny N}} [\bm{m}_1(\bm{X}_2),\ldots,\bm{m}_{K_1}(\bm{X}_2)]$ and $[\bm{m}_1(\bm{X}_2),\ldots,\bm{m}_{K_1}(\bm{X}_2)] \preceq_{\mbox{\tiny N}} [\bm{m}_1(\bm{X}_1),\ldots,\bm{m}_{K_1}(\bm{X}_1)]$, and we further assume that \\
$[\bm{m}_1(\bm{X}_1),\ldots,\bm{m}_{K_1}(\bm{X}_1)] \neq [\bm{m}_1(\bm{X}_2),\ldots,\bm{m}_{K_1}(\bm{X}_2)]$. We have the following two cases between \\
$[\bm{m}_1(\bm{X}_1),\ldots,\bm{m}_{K_1}(\bm{X}_1)]$ and $[\bm{m}_1(\bm{X}_2),\ldots,\bm{m}_{K_1}(\bm{X}_2)]$.

First case, we can assume that $[\bm{m}_1(\bm{X}_1),\ldots,\bm{m}_{K_1}(\bm{X}_1)] \prec_{\mbox{\tiny N}} [\bm{m}_1(\bm{X}_2),\ldots,\bm{m}_{K_1}(\bm{X}_2)]$. Then, we must have some $k_1= 1,2,\ldots,K_1$ such that $\bm{m}_i(\bm{X}_1) = \bm{m}_i(\bm{X}_2)$ for all $i < k_1$, and $\bm{m}_{k_1}(\bm{X}_1) \triangleleft \bm{m}_{k_1}(\bm{X}_2)$.  

Second case, we can assume that $[\bm{m}_1(\bm{X}_2),\ldots,\bm{m}_{K_1}(\bm{X}_2)] \prec_{\mbox{\tiny N}} [\bm{m}_1(\bm{X}_1),\ldots,\bm{m}_{K_1}(\bm{X}_1)]$. Then, we must have some $k_2= 1,2,\ldots,K_1$ such that $\bm{m}_i(\bm{X}_1) = \bm{m}_i(\bm{X}_2)$ for all $i < k_2$, and $\bm{m}_{k_2}(\bm{X}_2) \triangleleft \bm{m}_{k_2}(\bm{X}_1)$.  

Without loss of generality, we may assume that $k_1 \leq k_2$. If $k_1 < k_2$, we will have $\bm{m}_k(\bm{X}_1)= \bm{m}_{k_1}(\bm{X}_2)$, which is a contradition. On the other hand, if $k_1 = k_2$, we will have $\bm{m}_{k_1}(\bm{X}_1) \triangleleft \bm{m}_{k_1}(\bm{X}_2) = \bm{m}_{k_2}(\bm{X}_2) \triangleleft \bm{m}_{k_2}(\bm{X}_1)=\bm{m}_{k_1}(\bm{X}_1)$, which is an another contradition. Our assumption about \\
$[\bm{m}_1(\bm{X}_1),\ldots,\bm{m}_{K_1}(\bm{X}_1)] \neq [\bm{m}_1(\bm{X}_2),\ldots,\bm{m}_{K_1}(\bm{X}_2)]$ is incorrect. This implies that \\
$[\bm{m}_1(\bm{X}_1),\ldots,\bm{m}_{K_1}(\bm{X}_1)] = [\bm{m}_1(\bm{X}_2),\ldots,\bm{m}_{K_1}(\bm{X}_2)]$, which is antisymmetry. 

3. Transitivity:\\
If we have $[\bm{m}_1(\bm{X}_1),\ldots,\bm{m}_{K_1}(\bm{X}_1)] \preceq_{\mbox{\tiny N}} [\bm{m}_1(\bm{X}_2),\ldots,\bm{m}_{K_1}(\bm{X}_2)]$ and $[\bm{m}_1(\bm{X}_2),\ldots,\bm{m}_{K_1}(\bm{X}_2)] \preceq_{\mbox{\tiny N}} [\bm{m}_1(\bm{X}_3),\ldots,\bm{m}_{K_1}(\bm{X}_3)]$, then for some $k_1$ and $k_2$, we have
\begin{eqnarray}
\bm{m}_i(\bm{X}_1)&=&\bm{m}_i(\bm{X}_2), \mbox{~~for all $i < k_1$ and $\bm{m}_{k_1}(\bm{X}_1) \triangleleft \bm{m}_{k_1}(\bm{X}_2)$}, \nonumber \\
\bm{m}_i(\bm{X}_2)&=&\bm{m}_i(\bm{X}_3), \mbox{~~for all $i < k_2$ and $\bm{m}_{k_2}(\bm{X}_2) \triangleleft \bm{m}_{k_2}(\bm{X}_3)$}.\nonumber 
\end{eqnarray}
Without loss of generality, we can assume that $k_1 \leq k _2$. Then, we have $\bm{m}_i(\bm{X}_1)=\bm{m}_i(\bm{X}_2)=\bm{m}_i(\bm{X}_3)$ for all $i < k_1$, and $\bm{m}_{i}(\bm{X}_1) \trianglelefteq \bm{m}_{i}(\bm{X}_2)$ and $\bm{m}_{i}(\bm{X}_2) \trianglelefteq \bm{m}_{i}(\bm{X}_3)$ for $i \geq k_1$. Because we can perform partial ordering $\bm{m}_{i}(\bm{X}) \trianglelefteq \bm{m}_{i}(\bm{Y})$ for any $i$ and any pair of matrices $\bm{X}$ and $\bm{Y}$. Then, we also have $\bm{m}_{i}(\bm{X}_1) \trianglelefteq \bm{m}_{i}(\bm{X}_3)$. This implies that $[\bm{m}_1(\bm{X}_1),\ldots,\bm{m}_{K_1}(\bm{X}_1)] \preceq_{\mbox{\tiny N}} [\bm{m}_1(\bm{X}_3),\ldots,\bm{m}_{K_1}(\bm{X}_3)]$.

Since the order relation $ \prec_{\mbox{\tiny N}}$ satisfies reflexivity, antisymmetry, and transitivity, it is a partial ordering. 
$\hfill\Box$

\begin{theorem}[SNO is a partial ordering]\label{thm: SNO is a partial ordering}
Given two matrices $\bm{X}_1$ and $\bm{X}_2$ with same dimensions, the ordering $\preceq_{\mbox{\tiny SN}}$ provided by Definition~\ref{def: SNO} is a partial order relation.
\end{theorem}
\textbf{Proof:}
From spectral and nilpotent structures representation of matrices $\bm{X}_1$ and $\bm{X}_2$, we have
\begin{eqnarray}\label{eq1: thm: SNO is a partial ordering}
\mathfrak{R}(\bm{X}_1)&=&[\bm{\lambda_1}(\bm{X}_1),\ldots,\bm{\lambda}_{K_1}(\bm{X}_1),\bm{m}_1(\bm{X}_1),\ldots,\bm{m}_{K_1}(\bm{X}_1)], \nonumber \\
\mathfrak{R}(\bm{X}_2)&=&[\bm{\lambda_1}(\bm{X}_2),\ldots,\bm{\lambda}_{K_2}(\bm{X}_2),\bm{m}_1(\bm{X}_2),\ldots,\bm{m}_{K_2}(\bm{X}_2)].
\end{eqnarray}   
To prove that the lexicographic order defined on the Cartesian product of two partially ordered relations $[\bm{\lambda_1}(\bm{X}_1),\ldots,\bm{\lambda}_{K_1}(\bm{X}_1)] \preceq_w  
[\bm{\lambda_1}(\bm{X}_2),\ldots,\bm{\lambda}_{K_2}(\bm{X}_2)] $ and \\
 $[\bm{m}_1(\bm{X}_1),\ldots,\bm{m}_{K_1}(\bm{X}_1)] \preceq_{\mbox{\tiny N}} [\bm{m}_1(\bm{X}_2),\ldots,\bm{m}_{K_2}(\bm{X}_2)]$ is a partial ordering, we must verify the three properties of a partial order: reflexivity ,  antisymmetry , and  transitivity .

Let $P_1$ be the set of $[\bm{\lambda_1}(\bm{X}),\ldots,\bm{\lambda}_{K_1}(\bm{X})]$, and $P_2$ be the set of $[\bm{m}_1(\bm{X}),\ldots,\bm{m}_{K_1}(\bm{X})]$, given $(a, b)$, $(c, d) \in P = P_1 \times P_2$ (Cartesian product), the  lexicographic ordering  is defined as:
\[
[a, b] \preceq_{\mbox{\tiny SN}} [c, d] \quad \text{if and only if} \quad \begin{cases}
a \preceq_{w} c, & \text{or} \\
a = c \text{ and } b \preceq_{\mbox{\tiny N}}d.
\end{cases}
\]

1. Reflexivity: \\
A relation is reflexive if for all \((a, b) \in P\), \((a, b) \preceq_{\mbox{\tiny SN}} (a, b)\). From the definition of the lexicographic order: \(a = a\) and \(b \preceq_{\mbox{\tiny N}} b\). Because \(\preceq_{\mbox{\tiny N}}\) is reflexive on \(P_2\), the $\preceq_{\mbox{\tiny SN}}$ order is reflexive.

2. Antisymmetry:\\
A relation is antisymmetric if \((a, b) \preceq_{\mbox{\tiny SN}} (c, d)\) and \((c, d) \preceq_{\mbox{\tiny SN}} (a, b)\) imply \((a, b) = (c, d)\).

Considering the following two cases:\\
A: \((a, b) \preceq_{\mbox{\tiny SN}}  (c, d)\), so either: (1) \(a \preceq_w c\), or (2) \(a = c\) and \(b \preceq_{\mbox{\tiny N}} d\).\\
B: \((c, d) \preceq_{\mbox{\tiny SN}} (a, b)\), so either: (1) \(c \preceq_w a\), or (2) \(c = a\) and \(d \preceq_{\mbox{\tiny N}} b\).

From A and B: \\
- If \(a \prec_w c\) and \(c \prec_w a\), we have a contradiction.\\
- If \(a = c\), then \(b \preceq_{\mbox{\tiny N}}  d\) and \(d \preceq_{\mbox{\tiny N}}  b\). By the antisymmetry of \(\preceq_{\mbox{\tiny N}}\) in \(P_2\), \(b = d\). Thus, \((a, b) = (c, d)\). Hence, the $\preceq_{\mbox{\tiny SN}}$ order is antisymmetric.

3. Transitivity:\\
A relation is transitive if \((a, b) \preceq_{\mbox{\tiny SN}} (c, d)\) and \((c, d) \preceq_{\mbox{\tiny SN}} (e, f)\) imply \((a, b) \preceq_{\mbox{\tiny SN}} (e, f)\).

Considering the following two cases:\\
A: \((a, b) \preceq_{\mbox{\tiny SN}} (c, d)\), so either: (1) \(a \preceq_w c\), or (2) \(a = c\) and \(b \preceq_{\mbox{\tiny N}} d\).\\
B: \((c, d) \preceq_{\mbox{\tiny SN}} (e, f)\), so either: (1) \(c \preceq_w e\), or (2) \(c = e\) and \(d \preceq_{\mbox{\tiny N}} f\).

From A and B:
\begin{itemize}
\item If \(a \preceq_w c\) and \(c \preceq_w e\), then \(a \preceq_w e\), so \((a, b) \preceq_{\mbox{\tiny SN}} (e, f)\).
\item If \(a \preceq_w c\) and \(c = e\), then \(a \preceq_w e\), so \((a, b) \preceq_{\mbox{\tiny SN}} (e, f)\).
\item If \(a = c\) and \(b \preceq_{\mbox{\tiny N}} d\), and \(c \preceq_w e\), then \(a \preceq_w e\), so \((a, b) \preceq_{\mbox{\tiny SN}} (e, f)\).
\item If \(a = c\) and \(b \preceq_{\mbox{\tiny N}} d\), and \(c = e\) and \(d \preceq_{\mbox{\tiny N}} f\), then \(a = e\) and \(b \preceq_{\mbox{\tiny N}} f\), so \((a, b) \preceq_{\mbox{\tiny SN}} (e, f)\).
\end{itemize}
Then, the lexicographic order is transitive.

Since the $\preceq_{\mbox{\tiny SN}}$ order satisfies  reflexivity,  antisymmetry, and  transitivity, this theorem is proved.  
$\hfill\Box$

Next Corollary~\ref{cor: Similarity of two matrices with same SN representation} describes that any two matrices with identical SN representation are similar.

\begin{corollary}[Similarity of two matrices with same SN representation]\label{cor: Similarity of two matrices with same SN representation}
Given two matrices $\bm{X}_1$ and $\bm{X}_2$ with same size, if~~$\mathfrak{R}(\bm{X}_1)=\mathfrak{R}(\bm{X}_2)$, we have that these two matrices $\bm{X}_1$ and $\bm{X}_2$ are similar.
\end{corollary}
\textbf{Proof:}
From SN representation provided by Eq.~\eqref{eq: matrix SN representation1} (or Eq.~\eqref{eq: matrix SN representation2}), we have 
\begin{eqnarray}\label{eq: cor: Similarity of two matrices with same SN representation}
\bm{X}_1&=& \bm{U}_1\left(\bigoplus\limits_{k=1}^{K}\bigoplus\limits_{i=1}^{\alpha_{k}^{(\mathrm{G})}}\bm{J}_{m_{k,i}}(\lambda_{k})\right)\bm{U}_1^{-1}, \nonumber \\
\bm{X}_2&=& \bm{U}_2\left(\bigoplus\limits_{k=1}^{K}\bigoplus\limits_{i=1}^{\alpha_{k}^{(\mathrm{G})}}\bm{J}_{m_{k,i}}(\lambda_{k})\right)\bm{U}_2^{-1}.
\end{eqnarray}   
Then, we have $\bm{X}_1 = \bm{U}_1\bm{U}^{-1}_2\bm{X}_2 \bm{U}_2 \bm{U}^{-1}_1$, which shows that these two matrices $\bm{X}_1$ and $\bm{X}_2$ are similar since both matrices $\bm{U}_1$ and $\bm{U}_2$ are invertible matrices.
$\hfill\Box$

Loewner order is a special case of SNO by considering only $\preceq_w$ in $\preceq_{\mbox{\tiny SN}}$ and requiring two matrices $\bm{X}_1$ and $\bm{X}_2$ with $\mathfrak{R}(\bm{X}_1)=\mathfrak{R}(\bm{X}_2)$ having to be similar by identity matrix only. But, we have to note that $[\bm{\lambda_1}(\bm{X}_1),\ldots,\bm{\lambda}_{K_1}(\bm{X}_1)]\preceq_w [\bm{\lambda_1}(\bm{X}_2),\ldots,\bm{\lambda}_{K_2}(\bm{X}_2)$ does not have $\bm{X}_1 \leq_{Loewner} \bm{X}_2$ always. 

\section{Majorization Ordering for Complex-Valued Functions}\label{sec: Majorization Ordering for Complex-Valued Functions}

In this section, we will consider two kinds of functions. The first one are functions with multiple input variables discussed in Section~\ref{sec: Function with Multiple Variables} and the second one are functions with a single input variable discussed in Section~\ref{sec: Function with a Single Variable}.

\subsection{Function with Multiple Variables}\label{sec: Function with Multiple Variables}

In this section, we explore properties of Schur-convex functions defined over the complex domain, incorporating total order relations among complex numbers. We begin by presenting the definition of Schur-convex functions for the complex domain with total order relations. Building on the foundations laid in Chapters 2 and 3 of~\cite{marshall1979inequalities}, we extend key lemmas and properties of majorization from the real domain to the complex domain.

\begin{definition}\label{def: Schur convex}
Let \( f : \mathbb{C}^n \to \mathbb{C} \) be a symmetric function. \( f \) is said to be \emph{Schur-convex over complex domain} if, for any two vectors \(\mathbf{x}, \mathbf{y} \in \mathbb{C}^n\),  
\[
\mathbf{x} \prec \mathbf{y} \implies f(\mathbf{x}) \leq f(\mathbf{y}),
\]  
where \(\mathbf{x} \prec \mathbf{y}\) means that \(\mathbf{x}\) is majorized by \(\mathbf{y}\). In this context, we sort entries in $\bm{x}, \bm{y}$ by $x_1 \geq x_2 \geq \ldots \geq x_n$ and $y_1 \geq y_2 \geq \ldots \geq y_n$, and require:\\
1. \(\sum_{i=1}^k x_i \leq \sum_{i=1}^k y_i\) for all \(k = 1, 2, \ldots, n-1\), and \\
2. \(\sum_{i=1}^n x_i = \sum_{i=1}^n y_i\).

If the inequality is reversed (\(f(\mathbf{x}) \geq f(\mathbf{y})\)), \(f\) is said to be \emph{Schur-concave over complex domain.}
\end{definition}

The next lemma is about a special transformation acting on a complex-valued vector, named \emph{affine T-transformaion}. The matrix of a T-transform is expressed as:

\[
\bm{T} = \beta\bm{I} + (1 - \beta)\bm{Q},
\]
where $\beta \in \mathbb{C}$, $\bm{I}$ is an identity matrix and \( \bm{Q} \) is a permutation matrix that interchanges two coordinates only. Note that we relax the conventional requirement of $\beta$ with $ 0 \leq \beta \leq 1$, i.e., from convex combination to affine combination for matrices $\bm{I}$ and $\bm{Q}$.

Given a complex-valued vector $\bm{x}=[x_1,\ldots,x_{i-1},x_i,x_{i+1},\ldots,x_{j-1},x_j,x_{j+1},\ldots,x_n]$, we have $\bm{x}\bm{T}$, i.e., T-transformation for the vector $\bm{x}$, given by:

\[
\bm{x}\bm{T} = \left[x_1, \ldots, x_{i-1}, \beta x_i + (1 - \beta)x_j, x_{i+1}, \ldots, x_{j-1}, \beta x_j + (1 - \beta)x_i, x_{j+1}, \ldots, x_n\right].
\]

\begin{lemma}[Affine T-transformation over complex numbers]\label{lma: Affine T-transformation over Complex Numbers}
For two complex-valued vectors $\bm{x} \in \mathbb{C}^n$ and $\bm{y}  \in \mathbb{C}^n$ with $\bm{x}\prec\bm{y}$, the vector $\bm{x}$ can be obtained from the vector $\bm{y}$ by finite number of affine T-transformation operations. 
\end{lemma}
\textbf{Proof:}
There are two cases. The first case is that all entries in the vector $\bm{x}$ are permutation from entries in the vector $\bm{y}$. 

\textbf{Case I}

In such case, $\beta=0$ and we wish to show that any permutation of \( n \) elements can be achieved by a finite number of swaps. Recall a swap is represented by the  matrix \( \bm{Q} \) in a T-transformation. We will prove this by induction. 

Base Case: \( n = 2 \)\\
The set \( \{1, 2\} \) has two permutations: \( \{1, 2\} \) and \( \{2, 1\} \). To obtain \( \{2, 1\} \) from \( \{1, 2\} \), a single swap suffices (swap \( 1 \) and \( 2 \)). Thus, the statement holds for \( n = 2 \).

Inductive Hypothesis:\\
Assume that for \( n = k \), any permutation of \( k \) elements can be achieved using a finite number of swaps.

For \( n = k + 1 \), we consider a set \( \{1, 2, \dots, k+1\} \) with an arbitrary permutation of its elements, say \( \pi \). Let \( \pi = (\pi_1, \pi_2, \dots, \pi_{k+1}) \), where \( \pi_i \) indicates the element at the \( i \)-th position in the permutation. If \( \pi_{k+1} = k+1 \) (the last element is in its correct position), then the problem reduces to permuting the first \( k \) elements, which can be done by the inductive hypothesis using a finite number of swaps. On the other hand, if \( \pi_{k+1} \neq k+1 \), then we swap \( \pi_{k+1} \) with the element \( k+1 \) to place \( k+1 \) in its correct position. After this swap, the remaining \( k \) elements need to be permuted to match \( \pi \). By the inductive hypothesis, this can be achieved using a finite number of swaps. Therefore, any permutation of \( k+1 \) elements can be constructed using a finite number of swaps. 

By induction, the statement holds for all \( n \geq 2 \): any permutation of \( n \) elements can be achieved using a finite number of swaps. Thus, if entries in $\bm{x}$ are permutation from entries in the vector $\bm{y}$, we can convert the vector $\bm{y}$ to the vector $\bm{x}$ by finite number of T-transformation operations. 

\textbf{Case II}

The second case is that the vector $\bm{x}$ cannot be obtained by permuting entries in the vector $\bm{y}$. Without loss of generality, we assume that $x_1\geq x_2 \geq \ldots \geq x_n$ and $y_1\geq y_2 \geq \ldots \geq y_n$. Because $\bm{x} \prec \bm{y}$, we can find the index $i$ such that the value of $i$ is the largest index with $x_i < y_i$, and we also can find the index $j$ such that $i < j$ and $j$ is the smallest index with $x_j > y_j$. Then, we can set $\beta$ in the T-transformation as
\begin{eqnarray}\label{eq1: lma: Affine T-transformation over Complex Numbers}
\beta&=&1 - \frac{\min(y_i-x_i, x_j-y_j)}{y_i - y_j}.
\end{eqnarray}
Note that $y_i-x_i$, $x_j-y_j$ and $y_i-y_j$ are greater than $0+0\iota$. For sitations that $x_i,x_j,y_i$ and $y_j$ are real numbers, we have $0 \leq \beta \leq 1$. Let the vector $\bm{y}' = \beta\bm{y}\bm{I} + (1-\beta)\bm{y}\bm{Q}$ ( a T-transformation with the parameter $\beta$ given by Eq.~\eqref{eq1: lma: Affine T-transformation over Complex Numbers}), where $\bm{Q}$ is the permutation matrix to swap the index $i$ and the index $j$. Then, we have 
\begin{eqnarray}\label{eq2 lma: Affine T-transformation over Complex Numbers}
\bm{y}' &=&[y_1,\ldots,y_{i-1},y_i -\epsilon ,y_{i+1},\ldots,y_{j-1},y_j+\epsilon,y_{j+1},\ldots,x_n],
\end{eqnarray}
where $\epsilon = \min(y_i-x_i, x_j-y_j)$.

Now, we claim that $\bm{x} \prec \bm{y}' \prec \bm{y}$.  Because $y_i-x_i > 0+0 \iota$, $x_j-y_j > 0 + 0 \iota$, we have $\bm{y}' \prec \bm{y}$ immediatedly. On the other hand, we have the following relationships:
\begin{eqnarray}\label{eq3 lma: Affine T-transformation over Complex Numbers}
\sum\limits_{\ell=1}^p x_\ell \leq \sum\limits_{\ell=1}^p y'_\ell, 
\end{eqnarray}
where $p=1,2,\ldots,i-1$. For $p=i$, since $x_i \leq y_i-\epsilon (=y'_i)$, we have 
\begin{eqnarray}\label{eq4 lma: Affine T-transformation over Complex Numbers}
\sum\limits_{\ell=1}^i x_\ell \leq \sum\limits_{\ell=1}^i y'_\ell. 
\end{eqnarray}
Because $y'_\ell=y_\ell$ for $\ell=i+1,\ldots,j-1$ and $\sum\limits_{\ell=1}^n x_\ell = \sum\limits_{\ell=1}^n y'_\ell$, we have 
\begin{eqnarray}\label{eq5 lma: Affine T-transformation over Complex Numbers}
\sum\limits_{\ell=1}^p x_\ell \leq \sum\limits_{\ell=1}^p y'_\ell 
\end{eqnarray}
where $p=j+1,\ldots,n$. Hence, Eq.~\eqref{eq3 lma: Affine T-transformation over Complex Numbers} to Eq.~\eqref{eq5 lma: Affine T-transformation over Complex Numbers} imply that $\bm{x} \prec \bm{y}'$.

If we define $\gamma(\bm{x}, \bm{y})$ be the number of non-equal entries between $x_\ell$ and $y_\ell$ for $\ell=1,2,\ldots,n$. From Eq.~\eqref{eq2 lma: Affine T-transformation over Complex Numbers} and $\epsilon = \min(y_i-x_i, x_j-y_j)$, we have $y'_i = x_i$ or $y'_j = x_j$. This shows that $\gamma(\bm{x}, \bm{y}') \leq \gamma(\bm{x}, \bm{y})-1$ after one T-transformation operation. Since $\gamma(\bm{x}, \bm{y}) \leq n$, therefore, we also can convert the vector $\bm{y}$ to the vector $\bm{x}$ by finite number of T-transformation operations for Case II.
$\hfill\Box$

Below, we will define generalized doubly stochastic matrices for later Schur-convexity properties discussion. 
\begin{definition}\label{def: GDSM}
A matrix \( A = [a_{i,j}] \) is called \emph{generalized doubly stochastic} matrix if:
\begin{enumerate}[label=(\roman*)]
        \item All entries are complex numbers : \( a_{i, j} \in \mathbb{C} \) for all \( i, j \),
        \item The rows sum to 1: \( \sum\limits_{j=1}^n a_{i, j} = 1 \)~for all \( i \),
        \item The columns sum to 1: \( \sum\limits_{i=1}^n a_{i, j} = 1 \)~for all \( j \).
\end{enumerate}
Different from conventional doubly stochastic matrices, we allow entries from a generalized doubly stochastic matrix to be any complex numbers. 
\end{definition}

Next Lemma~\ref{lma: product of GDSM is still GDSM} will show that the product of two generalized doubly stochastic matrices will be a generalized doubly stochastic matrix again. 
\begin{lemma}\label{lma: product of GDSM is still GDSM}
The product of two \( n \times n \) generalized doubly stochastic matrices is also a generalized doubly stochastic matrix.
\end{lemma}
\textbf{Proof:}
Let \( \bm{A} \) and \( \bm{B} \) be two \( n \times n \) generalized doubly stochastic matrices. We aim to prove that the product \( \bm{C} = \bm{A}\bm{B} \), where \( \bm{C} = [c_{i,j}] \) and \( c_{i,j} = \sum_{k=1}^n a_{i,k}b_{k,j} \), is also a generalized doubly stochastic.

1. Entries are complex numbers:\\
Each entry of \( C \) is given by:
\[
c_{i,j} = \sum_{k=1}^n a_{i,k} b_{k,j}.
\]
Since \( a_{i,k} \in \mathbb{C} \) and \( b_{k,j} \in \mathbb{C} \) for all \( i, j, k \) (as \( \bm{A} \) and \( \bm{B} \) are doubly stochastic), it follows that:
\[
c_{i,j} \in \mathbb{C}.
\]
Thus, all entries of \( C \) are complex numbers.

2. Row Sums of \( \bm{C} \)\\
Consider the sum of the entries in the \( i \)-th row of \( \bm{C} \):
\[
\sum_{j=1}^n c_{i, j} = \sum_{j=1}^n \sum_{k=1}^n a_{i,k} b_{k,j}.
\]
Reordering the sums (by switching the order of summation):
\[
\sum_{j=1}^n c_{i,j} = \sum_{k=1}^n a_{i,k} \left( \sum_{j=1}^n b_{k,j} \right).
\]
Since \( \bm{B} \) is a generalized doubly stochastic matrix, we know \( \sum_{j=1}^n b_{k,j} = 1 \) for all \( k \). Substituting this:
\[
\sum_{j=1}^n c_{i,j} = \sum_{k=1}^n a_{i,k} \cdot 1 = \sum_{k=1}^n a_{i,k}.
\]
Now, since \( \bm{A} \) is doubly stochastic, \( \sum_{k=1}^n a_{i,k} = 1 \) for all \( i \). Therefore:
\[
\sum_{j=1}^n c_{i,j} = 1.
\]
This shows that the rows of \( \bm{C} \) sum to 1.

3. Column Sums of \( \bm{C} \)\\
Consider the sum of the entries in the \( j \)-th column of \( \bm{C} \):
\[
\sum_{i=1}^n c_{i,j} = \sum_{i=1}^n \sum_{k=1}^n a_{i,k} b_{k,j}.
\]
Reordering the sums (by switching the order of summation):
\[
\sum_{i=1}^n c_{i,j} = \sum_{k=1}^n \left( \sum_{i=1}^n a_{i,k} \right) b_{k,j}.
\]
Since \( \bm{A} \) is doubly stochastic, we know \( \sum_{i=1}^n a_{i,k} = 1 \) for all \( k \). Substituting this:
\[
\sum_{i=1}^n c_{i,j} = \sum_{k=1}^n 1 \cdot b_{k,j} = \sum_{k=1}^n b_{k,j}.
\]
Now, since \( \bm{B} \) is doubly stochastic, \( \sum_{k=1}^n b_{k,j} = 1 \) for all \( j \). Therefore:
\[
\sum_{i=1}^n c_{i,j} = 1.
\]
This shows that the columns of \( \bm{C} \) sum to 1.
$\hfill\Box$

\begin{lemma}[ Marjoization characterization by generalized double stochastic matrix]\label{lma:  marjoization characterization by generalized double stochastic matrix}
Given two complex-valued vectors with $\bm{x} \prec \bm{y}$, there exist a generalized doubly stochastic matrix $\bm{P}$ such that $\bm{x}=\bm{y}\bm{P}$. 
\end{lemma}
\textbf{Proof:}
Because T-transforms are generealized doubly stochastic matrices, the product of T-transformations is again a generalized doubly stochastic matrix from Lemma~\ref{lma: product of GDSM is still GDSM}. Therefore, there exists a  generealized doubly stochastic matrix $\bm{P}$ to have $\bm{x}=\bm{y}\bm{P}$. 
$\hfill\Box$

Below, we will discuss functions compositions associated with Schur-convex and Schur-concave functions. We consider the following composition relation:
\begin{eqnarray}\label{eq1: func comp Schur-conv}
f(\bm{x})=g(h_1(\bm{x}),\ldots,h_m(\bm{x})),
\end{eqnarray}
where functions $f, h_1,\ldots,h_m$ are symmetric complex valued function defined on domain $\mathcal{D} \in \mathbb{C}^n$, and the function $g$ is a complex valued function defined on $\mathbb{C}^n$.  When we say ``increasing" for functions $f,g, h_1,\ldots,h_m$ if the function value increases with respect to each variable increasing. Similarly, we say ``decreasing" if these functions decrease with respect to each variable decreasing. 

We will have the following property related to Eq.~\eqref{eq1: func comp Schur-conv}. If the function $g$ is an increasing function, and each function $h_i$ is a Schur-convex function  over the domain $\mathcal{D}$, the function $f(\bm{x})$ will be a Schur-convex on the domain $\mathcal{D}$. This is the property (1) in Table~\ref{tab: comp properties I}. We can prove this property easily. If $\bm{x} \prec \bm{y}$ on $\mathcal{D}$, and each function $h_i$ is a Schur-convex function, i.e., $h_i(\bm{x}) \leq h_i(\bm{y})$ for $i=1,2,\ldots,m$, we have
\begin{eqnarray}\label{eq2: func comp Schur-conv}
g(h_1(\bm{x}),\ldots,h_m(\bm{x})) \leq g(h_1(\bm{y}),\ldots,h_m(\bm{y})),
\end{eqnarray}
because the function $g$ is an increasing function. Given basic conditions of $f$ and $h_i$ for $i=1,2,\ldots,m$, we have the following Table~\ref{tab: comp properties I} to summarize other properties about functions compositions associated with Schur-convex and Schur-concave functions. Other properties can also be proved similar to the Property (1) in Table~\ref{tab: comp properties I}.  

\begin{table}[h!]
\centering
\begin{tabular}{|c|c|c|c|}
\hline
Properties & $g$ & $h_i$ & $f$ \\ \hline
(1) & Increasing           & Schur-convex    & Schur-convex  \\ \hline
(2) & Increasing           & Schur-concave  & Schur-concave  \\ \hline
(3) & Decreasing          & Schur-convex   & Schur-concave  \\ \hline
(4) & Decreasing          & Schur-concave & Schur-convex  \\ \hline
(5) & Increasing           & Increasing Schur-convex    & Increasing Schur-convex  \\ \hline
(6) & Increasing           & Increasing Schur-concave  & Increasing Schur-concave  \\ \hline
(7) & Decreasing          & Decreasing Schur-convex   & Decreasing Schur-concave  \\ \hline
(8) & Decreasing          & Decreasing Schur-concave & Decreasing Schur-convex  \\ \hline
(9) & Increasing           & Decreasing Schur-convex    & Decreasing Schur-convex  \\ \hline
(10) & Increasing           & Increasing Schur-concave  & Increasing Schur-concave  \\ \hline
(11) & Decreasing          & Decreasing Schur-convex   & Increasing Schur-concave  \\ \hline
(12) & Decreasing          & Increasing Schur-concave & Decreasing Schur-convex  \\ \hline
\end{tabular}
\caption{Properties about Functions Compositions, I}
\label{tab: comp properties I}
\end{table}

Now, we will consider another format of functions compostion. We consider the following composition:
\begin{eqnarray}\label{eq3: func comp Schur-conv}
f(\bm{x})=g(h(x_1),\ldots,h(x_n)),
\end{eqnarray}
where functions $f,g$ are symmetric complx valued function defined on domain $\mathbb{C}^n$, and the function $h$ is a complx valued function defined on $\mathbb{C}$. 

For any given majorizatoin relation $[x_1, x_2] \prec [y_1, y_2]$, from Lemma~\ref{lma: Affine T-transformation over Complex Numbers}, we have
 \begin{eqnarray}\label{eq3-1: func comp Schur-conv}
x_1&=&\alpha y_1 + (1-\alpha) y_2, \nonumber \\
x_2&=&(1-\alpha)y_1  + \alpha y_2.
\end{eqnarray}
We also require the function $h$ to satisfy
 \begin{eqnarray}\label{eq3-2: func comp Schur-conv}
\alpha[h(y_1),h(y_2)]+ (1-\alpha)[h(y_2),h(y_1)] \prec [h(y_1), h(y_2)].
\end{eqnarray}
We call the condition for the function $h$ provided by Eq.~\eqref{eq3-1: func comp Schur-conv} and Eq.~\eqref{eq3-2: func comp Schur-conv} as \emph{complex domain majorization condition}. Since both $\alpha$ and the function $h$ are complex valuled, we have to apply Eq.~\eqref{eq1: Multiplication} and Eq.~\eqref{eq2: Multiplication} based on the relationship between $h(y_1)$ and $h(y_2)$.

We consider a complex valued convex function $h(x)$ with affine combination as the following condition is satisfied:
\begin{eqnarray}\label{eq4: func comp Schur-conv}
h(\alpha x_1 + (1-\alpha) x_2) \leq \alpha h(x_1) + (1-\alpha) h(x_2),
\end{eqnarray}
where $\alpha$ is any complex number.  

Given the function $g$ is an increasing and Schur-convex function, and $h$ is a convex function with affine combination and satisfying complex domain majorization condition, the function $f$ provided by Eq.~\eqref{eq3: func comp Schur-conv} is a Schur-convex function. This is the property (1) in Table~\ref{tab: comp properties II}. From Lemma~\ref{lma: Affine T-transformation over Complex Numbers}, it is sufficient to prove Schur-convex related properties by considering functions with only two argeuments. Consider $n=2$ in  Eq.~\eqref{eq3: func comp Schur-conv} with $[x_1, x_2] \prec [y_1, y_2]$, from Lemma~\ref{lma: Affine T-transformation over Complex Numbers}, we can express $x_1$ and $x_2$ by  
\begin{eqnarray}\label{eq5: func comp Schur-conv}
x_1&=&\alpha y_1 + (1-\alpha) y_2, \nonumber \\
x_2&=&(1-\alpha)y_1  + \alpha y_2.
\end{eqnarray}
Then, we have 
\begin{eqnarray}\label{eq6 func comp Schur-conv}
g(h(x_1),h(x_2))&=_1& g(h(\alpha y_1 + (1-\alpha) y_2), h((1-\alpha)y_1  + \alpha y_2)) \nonumber \\
&\leq_2& g(\alpha h(y_1) + (1-\alpha) h(y_2), (1-\alpha)h(y_1)  + \alpha h(y_2)) \nonumber \\
&=& g(\alpha[h(y_1),h(y_2)]+ (1-\alpha)[h(y_2),h(y_1)])
\end{eqnarray}
where we apply Eq.~\eqref{eq5: func comp Schur-conv} in $=_1$, and we apply $h$ as a convex function with affine combination to $\leq_2$.

Since $h$ satisfies \emph{complex domain majorization condition}, we have
\begin{eqnarray}\label{eq7 func comp Schur-conv}
\alpha[h(y_1),h(y_2)]+ (1-\alpha)[h(y_2),h(y_1)] \prec [h(y_1), h(y_2)].
\end{eqnarray}
From Schur-convex of the function $g$, we have
\begin{eqnarray}\label{eq7 func comp Schur-conv}
g(\alpha[h(y_1),h(y_2)]+ (1-\alpha)[h(y_2),h(y_1)]) < g(h(y_1), h(y_2)).
\end{eqnarray}
From Eq.~\eqref{eq6 func comp Schur-conv} and Eq.~\eqref{eq7 func comp Schur-conv}, we have
\begin{eqnarray}\label{eq8 func comp Schur-conv}
g(h(x_1),h(x_2))  <  g(h(y_1), h(y_2)).
\end{eqnarray}
Therefore, the property (1) in Table~\ref{tab: comp properties II} is proved.

We have the following Table~\ref{tab: comp properties II} to summarize other properties about functions compositions associated with Schur-convex and Schur-concave functions with the format provided by Eq.~\eqref{eq3: func comp Schur-conv}~\footnote{Convex and concave in Table~\ref{tab: comp properties II} is under affince combination.}. The function $h$ in Table~\ref{tab: comp properties II} satisfies complex domain majorization condition. Other properties can also be proved similar to the property (1). 

\begin{table}[h!]
\centering
\begin{tabular}{|c|c|c|c|}
\hline
Properties & $g$ & $h$ & $f$ \\ \hline
(1) & Increasing and Schur-convex  & Convex    & Schur-convex  \\ \hline
(2) & Decreasing and Schur-convex  & Concave    & Schur-convex  \\ \hline
(3) & Increasing and Schur-convex  & Increasing and Convex    & Increasing and Schur-convex  \\ \hline
(4) & Decreasing and Schur-convex  & Decreasing and  Concave    & Increasing and Schur-convex  \\ \hline
(5) & Increasing and Schur-convex  & Decreasing and Convex    & Decreasing and Schur-convex  \\ \hline
(6) & Decreasing and Schur-convex  & Increasing and  Concave    & Decreasing and Schur-convex   \\ \hline
(7) & Decreasing and Schur-concave & Convex    & Schur-concave  \\ \hline
(8) & Increasing and Schur-concave  & Concave    & Schur-concave  \\ \hline
(9) & Decreasing and Schur-concave  & Increasing and Convex    & Decreasing and Schur-concave \\ \hline
(10) & Increasing and Schur-concave  & Decreasing and  Concave    & Decreasing and Schur-concave  \\ \hline
(11) & Decreasing and Schur-concave  & Decreasing and Convex    &  Increasing and Schur-concave  \\ \hline
(12) & Increasing and Schur-concave  & Increasing and  Concave    &  Increasing and Schur-concave   \\ \hline
\end{tabular}
\caption{Properties about Functions Compositions, II}
\label{tab: comp properties II}
\end{table}

Next theorem will present Schur–Ostrowski criterion for function over a complex domain.
\begin{theorem}[Schur–Ostrowski criterion over complex domain]\label{thm: Schur–Ostrowski criterion}
Let $f(x_1, x_2, \ldots,x_n): \mathbb{C}^n \rightarrow \mathbb{C}$ be a symmetric function, and the partial derivatives $\frac{\partial f}{\partial x_i}$ exist. We assume that the domain of the function $f$ is $\mathcal{D}$, which is a bounded region. We define 
\begin{eqnarray}\label{eq1: thm: Schur–Ostrowski criterion}
DR_{i.j}(f)&=& \Re(\frac{\partial f}{\partial x_i}-\frac{\partial f}{\partial x_j}),\nonumber \\
DI_{i.j}(f)&=& \Im(\frac{\partial f}{\partial x_i}-\frac{\partial f}{\partial x_j}).
\end{eqnarray}
The entry update $\epsilon \in \mathbb{C}$ with $\epsilon \geq 0+0\iota$ used in T-transformation to transform any vector $\bm{x} \in \mathcal{D}$ to any vector $\bm{y} \in \mathcal{D}$ satisfies the following:
\begin{eqnarray}\label{eq2: thm: Schur–Ostrowski criterion}
0 \leq \Re(\epsilon) \leq C_1, \nonumber  \\
C_2 \leq \Im(\epsilon) \leq C_3, 
\end{eqnarray}
where $C_1, C_2$ and $C_3$ are real numbers. Constants $C_1, C_2$ and $C_3$ exists as the domain $\mathcal{D}$ is bounded.  

Then $f$ is a Schur-convex over a bounded domain in $ \mathbb{C}^n$ if and only if
\begin{eqnarray}\label{eq3: thm: Schur–Ostrowski criterion}
x_i \geq x_j,~\mbox{and~}     
  \begin{cases}
       C_3 DI_{i.j}(f) \leq 0 & \mbox{~if $DR_{i.j}(f) \geq 0$ and $DI_{i.j}(f) \geq 0$,}\\
       C_2 DI_{i.j}(f) \leq 0 & \mbox{~if $DR_{i.j}(f) \geq 0$ and $DI_{i.j}(f) < 0$,}\\
       C_1 DR_{i.j}(f) \geq C_3 DI_{i.j}(f) & \mbox{~if $DR_{i.j}(f) < 0$ and $DI_{i.j}(f) \geq 0$,}\\
       C_1 DR_{i.j}(f) \geq C_2 DI_{i.j}(f)  & \mbox{~if $DR_{i.j}(f) < 0$ and $DI_{i.j}(f) < 0$,}\\
    \end{cases}
\end{eqnarray}
where $1 \leq i, j \leq n$.
\end{theorem}
\textbf{Proof:}
We begin by demonstrating necessity. We need to show that if \( f \) is Schur-convex, then the condition  
\[
x_i \geq x_j,~\mbox{and~}     
  \begin{cases}
       C_3 DI_{i.j}(f) \leq 0 & \mbox{~if $DR_{i.j}(f) \geq 0$ and $DI_{i.j}(f) \geq 0$,}\\
       C_2 DI_{i.j}(f) \leq 0 & \mbox{~if $DR_{i.j}(f) \geq 0$ and $DI_{i.j}(f) < 0$,}\\
       C_1 DR_{i.j}(f) \geq C_3 DI_{i.j}(f) & \mbox{~if $DR_{i.j}(f) < 0$ and $DI_{i.j}(f) \geq 0$,}\\
       C_1 DR_{i.j}(f) \geq C_2 DI_{i.j}(f)  & \mbox{~if $DR_{i.j}(f) < 0$ and $DI_{i.j}(f) < 0$,}\\
    \end{cases}
\]
is satisfied.

If \( \mathbf{x} \prec \mathbf{y} \) (i.e., \( \mathbf{x} \) is majorized by \( \mathbf{y} \)), where $\bm{x}, \bm{y} \in \mathcal{D}$, Schur-convexity implies \( f(\mathbf{x}) \leq f(\mathbf{y}) \).  A majorization relation \( \mathbf{x} \prec \mathbf{y} \) can be decomposed into a sequence of T-transformation from Lemma~\ref{lma: Affine T-transformation over Complex Numbers}. It suffices to verify the condition for an elementary transposition of the form:
   \[
   \mathbf{x}' = (x_1, \dots, x_i + \epsilon, \dots, x_j - \epsilon, \dots, x_n),
   \]
   where \( \epsilon > 0+0\iota \) and \( x_i \geq x_j \). For \( f \) to be non-decreasing along this transformation:
   \[
   f(\mathbf{x}') - f(\mathbf{x}) \geq 0+0\iota.
   \]

Using a Taylor expansion, we approximate the change in \( f \):
   \[
   f(\mathbf{x}') - f(\mathbf{x}) \approx \frac{\partial f}{\partial x_i} (\epsilon) + \frac{\partial f}{\partial x_j} (-\epsilon) = \epsilon \left( \frac{\partial f}{\partial x_i} - \frac{\partial f}{\partial x_j} \right).
   \]
Schur-convexity requires:
\begin{eqnarray}\label{eq4: thm: Schur–Ostrowski criterion}
   \epsilon \left( \frac{\partial f}{\partial x_i} - \frac{\partial f}{\partial x_j} \right)&=&(\Re(\epsilon)+\Im(\epsilon)\iota)(DR_{i.j}(f)+DI_{i.j}(f)\iota)\nonumber \\
&=&(\Re(\epsilon)DR_{i.j}(f) - \Im(\epsilon)DI_{i.j}(f))+(\Im(\epsilon)DR_{i.j}(f) + \Re(\epsilon)DR_{i.j}(f))\iota \nonumber \\
&\geq&0+0\iota .
\end{eqnarray}
From the bounds of $\epsilon$ provided by Eq.~\eqref{eq2: thm: Schur–Ostrowski criterion}, we have to consider the following four cases with respect to the sign of $DR_{i.j}(f)$ and $DI_{i.j}(f)$:

Case I: $DR_{i.j}(f) \geq 0$ and $DI_{i.j}(f) \geq 0$\\
We have 
\begin{eqnarray}\label{eq5-1: thm: Schur–Ostrowski criterion}
0~~\leq~\Re(\epsilon) DR_{i.j}(f) &\leq& C_1 DR_{i.j}(f), \nonumber \\
C_2 DI_{i.j}(f)~\leq~~\Im(\epsilon) DI_{i.j}(f) &\leq& C_3 DI_{i.j}(f). 
\end{eqnarray}
The condition $C_3 DI_{i.j}(f) \leq 0$ ensures that $(\Re(\epsilon)DR_{i.j}(f) - \Im(\epsilon)DI_{i.j}(f))\geq 0$.  

Case II: $DR_{i.j}(f) \geq 0$ and $DI_{i.j}(f) < 0$\\
We have 
\begin{eqnarray}\label{eq5-2: thm: Schur–Ostrowski criterion}
0~~\leq~\Re(\epsilon) DR_{i.j}(f) &\leq& C_1 DR_{i.j}(f), \nonumber \\
C_2 DI_{i.j}(f)~\geq~~\Im(\epsilon) DI_{i.j}(f) &\geq& C_3 DI_{i.j}(f). 
\end{eqnarray}
The condition $C_2 DI_{i.j}(f) \leq 0$ ensures that $(\Re(\epsilon)DR_{i.j}(f) - \Im(\epsilon)DI_{i.j}(f))\geq 0$.  

Case III: $DR_{i.j}(f) < 0$ and $DI_{i.j}(f) \geq 0$\\
We have 
\begin{eqnarray}\label{eq5-3: thm: Schur–Ostrowski criterion}
0~~\geq~\Re(\epsilon) DR_{i.j}(f) &\geq& C_1 DR_{i.j}(f), \nonumber \\
C_2 DI_{i.j}(f)~\leq~~\Im(\epsilon) DI_{i.j}(f) &\leq& C_3 DI_{i.j}(f). 
\end{eqnarray}
The condition $C_1 DR_{i.j}(f) \geq C_3 DI_{i.j}(f) $ ensures that $(\Re(\epsilon)DR_{i.j}(f) - \Im(\epsilon)DI_{i.j}(f))\geq 0$.  

Case IV: $DR_{i.j}(f) < 0$ and $DI_{i.j}(f) < 0$\\
We have 
\begin{eqnarray}\label{eq5-4: thm: Schur–Ostrowski criterion}
0~~\geq~\Re(\epsilon) DR_{i.j}(f) &\geq& C_1 DR_{i.j}(f), \nonumber \\
C_2 DI_{i.j}(f)~\geq~~\Im(\epsilon) DI_{i.j}(f) &\geq& C_3 DI_{i.j}(f). 
\end{eqnarray}
The condition $C_1 DR_{i.j}(f) \geq C_2 DI_{i.j}(f) $ ensures that $(\Re(\epsilon)DR_{i.j}(f) - \Im(\epsilon)DI_{i.j}(f))\geq 0$.  
All these four cases make Eq.~\eqref{eq4: thm: Schur–Ostrowski criterion} valid. Thus, the condition is necessary for Schur-convexity.

Now, we wish to show sufficiency. We need to show that if the condition  
\[
x_i \geq  x_j  \mbox{~and~}  
  \begin{cases}
       C_3 DI_{i.j}(f) \leq 0 & \mbox{~if $DR_{i.j}(f) \geq 0$ and $DI_{i.j}(f) \geq 0$,}\\
       C_2 DI_{i.j}(f) \leq 0 & \mbox{~if $DR_{i.j}(f) \geq 0$ and $DI_{i.j}(f) < 0$,}\\
       C_1 DR_{i.j}(f) \geq C_3 DI_{i.j}(f) & \mbox{~if $DR_{i.j}(f) < 0$ and $DI_{i.j}(f) \geq 0$,}\\
       C_1 DR_{i.j}(f) \geq C_2 DI_{i.j}(f)  & \mbox{~if $DR_{i.j}(f) < 0$ and $DI_{i.j}(f) < 0$,}\\
    \end{cases}
\]
holds for all \( i, j \), then \( f \) is Schur-convex.

Since any majorization $\bm{x} \prec \bm{y}$, where $\bm{x}, \bm{y} \in \mathcal{D}$, can be decomposed into a finite sequence of T-transformation from Lemma~\ref{lma: Affine T-transformation over Complex Numbers}, it suffices to prove the result for a single T-transformation operation.

Suppose $\bm{x}'$ is obtained from $\bm{x}$ by transferring $\epsilon > 0+0\iota$ from $x_i$ to $x_j$, where $x_i \geq x_j$, i.e., $x_i' = x_i + \epsilon$ and $x_j' = x_j - \epsilon$ by keeping identical with remaining coordinates values. By using the symmetry and differentiability of $f$, the change in $f$ becomes:
\begin{eqnarray}
f(\bm{x}') - f(\bm{x})&\approx&\frac{\partial f}{\partial x_i} (\epsilon) + \frac{\partial f}{\partial x_j} (-\epsilon)\nonumber \\
&=&\epsilon \left( \frac{\partial f}{\partial x_i} - \frac{\partial f}{\partial x_j} \right)\nonumber \\
&=&(\Re(\epsilon)+\Im(\epsilon)\iota)(DR_{i.j}(f)+DI_{i.j}(f)\iota)\nonumber \\
&=&(\Re(\epsilon)DR_{i.j}(f) - \Im(\epsilon)DI_{i.j}(f))+(\Im(\epsilon)DR_{i.j}(f) + \Re(\epsilon)DR_{i.j}(f))\iota 
\end{eqnarray}
Then, the assumption given by Eq.~\eqref{eq3: thm: Schur–Ostrowski criterion} ensures that $f(\bm{x}') - f(\bm{x}) \geq 0+0\iota$, i.e., $f(\bm{x}') \geq f(\bm{x})$. This shows that $f(\bm{x})$ is non-decreasing under a T-transformation operation.

By repeating this argument for each step in the sequence of T-transformation operations required to transform $\bm{x}$ into $\bm{y}$, we have $\bm{x} \prec \bm{y} \implies f(\bm{x}) \leq f(\bm{y})$. Hence, the function $f$ is Schur-convex.
$\hfill \Box$

\begin{remark}
(1) For the function $f$ over the complex domain, we do not have conventional condition:
\begin{eqnarray}\label{eq: Schur–Ostrowski criterion conv}
(x_i - x_j) \left( \frac{\partial f}{\partial x_i} - \frac{\partial f}{\partial x_j} \right) \geq 0+0\iota.
\end{eqnarray}
Since $x_i \geq x_j \mbox{~and~}  \frac{\partial f}{\partial x_i} \geq \frac{\partial f}{\partial x_j}$ does not gurantee to have Eq.~\eqref{eq: Schur–Ostrowski criterion conv} valid because complex numbers with total order is not an ordered field. \\
(2) This conditions provided by Theorem~\ref{thm: Schur–Ostrowski criterion} means that \( f \) increases more rapidly in larger components of \( \mathbf{x} \) than in smaller ones, consistent with the majorization relation. 
\end{remark}

\subsection{Function with a Single Variable}\label{sec: Function with a Single Variable}

In this section, we will study the $\preceq_w$ relation under those functions with single variable. The results of this study is very crucial for us to characterize monotonicity and convexity for the function of general matrices under $\preceq_{\mbox{\tiny SN}}$ relation.   

Given a complex-valued function $f(x)$, we say that the function $f$ is \emph{$\prec_w$-preserving} if, given $\mathbf{x} \prec_w \mathbf{y}$ we have
\begin{eqnarray}\label{eq1: thm: prec w monotone f(x) cond}
\mathbf{f(x)}\define[f(x_1),\ldots,f(x_n)] \prec_w \mathbf{f(y)}\define[f(y_1),\ldots,f(y_n)],
\end{eqnarray} 
where $x_1 \geq x_2 \geq \cdots \geq x_n$, and $y_1 \geq y_2 \geq \cdots \geq y_n$. 

We have to prepare the following Lemma~\ref{lma: f(x,y) monotone wrt x at diff y} about the monotonicity conditions of a bi-variable real-valued function. 
\begin{lemma}\label{lma: f(x,y) monotone wrt x at diff y}
Given any two pairs of real variables $(x_1,y_1)$ and $(x_2,y_2)$ such that $x_1 \neq x_2$. If the real-valued function $f(x,y)$ satisfies:
\begin{subequations}\label{eq1-1: lma: f(x,y) monotone wrt x at diff y}
\begin{align}
i)~& \frac{\partial f}{\partial x} > 0; \label{eq1-1:subeq1} \\
ii)~& (y_2 - y_1)\frac{\partial f(x_1,y)}{\partial y} > 0~\mbox{for $y$ between $y_1$ and $y_2$}. \label{eq1-1:subeq2}
\end{align}
\end{subequations}
We have the following property:
\begin{eqnarray}\label{eq1-2: lma: f(x,y) monotone wrt x at diff y}
f(x_1,y_1) &< &f(x_2,y_2), \mbox{~whenever $x_1 < x_2$}.
\end{eqnarray}

Similarly,  If he function $f(x,y)$ satisfies:
\begin{subequations}\label{eq2-1: lma: f(x,y) monotone wrt x at diff y}
\begin{align}
i)~& \frac{\partial f}{\partial x} < 0; \label{eq2-1:subeq1} \\
ii)~& (y_2 - y_1)\frac{\partial f(x_1,y)}{\partial y} < 0~\mbox{for $y$ between $y_1$ and $y_2$}. \label{eq2-1:subeq2}
\end{align}
\end{subequations}
We have the following property:
\begin{eqnarray}\label{eq2-2: lma: f(x,y) monotone wrt x at diff y}
f(x_1,y_1) &> &f(x_2,y_2), \mbox{~whenever $x_1 < x_2$}.
\end{eqnarray}
\end{lemma}
\textbf{Proof:}
We prove conditions for Eq.~\eqref{eq1-2: lma: f(x,y) monotone wrt x at diff y} only. As proof of conditions for Eq.~\eqref{eq2-2: lma: f(x,y) monotone wrt x at diff y} can be obtained similarly.

Without loss of generality, if we assume $x_1 < x_2$, we have 
\begin{eqnarray}\label{eq3: lma: f(x,y) monotone wrt x at diff y}
f(x_1,y_2)&=&f(x_1,y_1) + \int_{y_1}^{y_2}\frac{\partial f(x_1, y)}{\partial y} dy.
\end{eqnarray}
From condition given by Eq.~\eqref{eq1-1:subeq2} and Eq.~\eqref{eq3: lma: f(x,y) monotone wrt x at diff y}, we have $f(x_1,y_2) > f(x_1,y_1)$.  Then, from the condition given by Eq.~\eqref{eq1-1:subeq1}, we have $f(x_1,y_2) < f(x_2,y_2)$. Therefore, we have the desired property given by Eq.~\eqref{eq1-2: lma: f(x,y) monotone wrt x at diff y}.
$\hfill\Box$

We are ready to present the following Theorem~\ref{thm: prec w monotone f(x) cond inc f} to characterize the monotonoicity under $\prec_w$ by monotone increasing funciton $f(z)$. 
\begin{theorem}[Majorization preserving by monotone increasing funciton]\label{thm: prec w monotone f(x) cond inc f}
Given \(\mathbf{x} = [x_1, x_2, \ldots, x_n] \prec_w \mathbf{y} = [y_1, y_2, \ldots, y_n]\), we consider a complex-valued function \(f(z = u + v \iota) = f_{\Re}(u, v) + f_{\Im}(u, v)\), where the partial derivatives of \(f_{\Re}(u, v)\) and \(f_{\Im}(u, v)\) exist. The function \(f(z)\) satisfies the following conditions~\footnote{The pair of real variables are selected from the real part and the imginary part from complex numbers $x_1, x_2, \ldots, x_n, y_1, y_2, \ldots, y_n$ when we apply Lemma~\ref{lma: f(x,y) monotone wrt x at diff y} conditions.}:  
\begin{align}\label{eq1: thm: prec w monotone f(x) cond inc f}
\begin{cases}
f_{\Re} \text{ satisfies conditions from Eq.~\eqref{eq1-1: lma: f(x,y) monotone wrt x at diff y}}, 
& \text{Case I: $z$ increment along } \Re(z); \\
f_{\Re} \text{ is constant, and } f_{\Im} \text{ satisfies conditions from } 
& \text{Eq.~\eqref{eq1-1: lma: f(x,y) monotone wrt x at diff y}}, \\
& \text{Case II: $z$ increment along } \Re(z); \\
\frac{\partial f_{\Re}}{\partial v} > 0, 
& \text{Case III: $z$ stable along } \Re(z) \text{ but increment along } \Im(z); \\
f_{\Re} \text{ is independent of } v, \text{ and } \frac{\partial f_{\Im}}{\partial v} > 0, 
& \text{Case IV: $z$ stable along } \Re(z) \text{ but increment along } \Im(z).
\end{cases}
\end{align}
and the following \emph{difference sum conidtion} for the increasing function $f$
\begin{eqnarray}\label{eq1-1: thm: prec w monotone f(x) cond inc f}
f(x_i)-f(y_i)&\leq&\sum\limits_{j=1}^{i-1}(f(y_j)-f(x_j)),
\end{eqnarray}
where $i=2,\ldots,n$. Then, the function $f$ is $\prec_w$ preserving.
\end{theorem}
\textbf{Proof:}
We first show that any complex-valued function $f(z)$ is monotone increasing with respect to the variable $z$ will have $\prec_w$ preserving property. 

Given $\mathbf{x}=[x_1,x_2,\ldots,x_n] \prec_w \mathbf{y}=[y_1,y_2,\ldots,y_n]$, we have
\begin{eqnarray}\label{eq2: thm: prec w monotone f(x) cond inc f}
\sum_{i=1}^k x_i \leq \sum_{i=1}^k y_i \quad \text{for all } k = 1, 2, \ldots, n.
\end{eqnarray}
If we consider vectors $\bm{f(x)}=[f(x_1),f(x_2),\ldots,f(x_n)]$ and $\bm{f(y)}=[f(y_1),f(y_2),\ldots,f(y_n)]$, from the monotone increasing property of the function $f$ and the difference sum condition provided by Eq.~\eqref{eq1-1: thm: prec w monotone f(x) cond inc f}, we have
\begin{eqnarray}\label{eq3: thm: prec w monotone f(x) cond inc f}
\sum_{i=1}^k f(x_i) \leq \sum_{i=1}^k f(y_i) \quad \text{for all } k = 1, 2, \ldots, n,
\end{eqnarray}
where $f(x_1) \geq f(x_2) \geq \ldots \geq f(x_n)$ and $f(y_1) \geq f(y_2) \geq \ldots \geq f(y_n)$ since $f(z)$ is monotone increasing with respect to the variable $z$, we have $f(z_1)<f(z_2)$ if $z_1 < z_2$. Therefore, we have $\bm{f(x)}\prec_w \bm{f(y)}$ given $\mathbf{x}\prec_w \mathbf{y}$.

Next, we will demonstrate that any complex-valued function $f(z)$ that satisfies condition provided by Eq.~\eqref{eq1: thm: prec w monotone f(x) cond inc f} is a monotone increasing function. There are four cases for the function $f(z)$ to be a monotone increasing function with respect to the increasing of the variable $z$: 
\begin{eqnarray}\label{eq4: thm: prec w monotone f(x) cond inc f}
\begin{cases}
f_\Re(u_1, v_1) < f_\Re(u_2,v_2), &~\mbox{if $u_1 < u_2$, case I;}\\
\mbox{or} \\
f_\Re(u_1, v_1) = f_\Re(u_2,v_2), f_\Im(u_1, v_1) < f_\Im(u_2,v_2), &~\mbox{if $u_1 < u_2$, case II;}\\
\mbox{or} \\
f_\Re(u_1, v_1) < f_\Re(u_2,v_2), &~\mbox{if $u_1 = u_2$, $v_1 < v_2$, case III;}\\
\mbox{or} \\
f_\Re(u_1, v_1) = f_\Re(u_2,v_2), f_\Im(u_1, v_1) < f_\Im(u_2,v_2), &~\mbox{if $u_1 = u_2$, $v_1 < v_2$, case IV;}\\
\end{cases}
\end{eqnarray}

For case I in Eq.~\eqref{eq4: thm: prec w monotone f(x) cond inc f}, the conditions provided by Eq.~\eqref{eq1-1: lma: f(x,y) monotone wrt x at diff y} will make $f_\Re(u_1, v_1) < f_\Re(u_2,v_2)$ if $u_1 < u_2$. 

For case II in Eq.~\eqref{eq4: thm: prec w monotone f(x) cond inc f}, we have to require $f_\Re(u,v)$ as a constant since $f_\Re$ is independent of the variable $u$ and $v$. Again, the conditions provided by  Eq.~\eqref{eq1-1: lma: f(x,y) monotone wrt x at diff y} will make $f_\Im(u_1, v_1) < f_\Im(u_2,v_2)$ if $u_1 < u_2$. 

For case III in Eq.~\eqref{eq4: thm: prec w monotone f(x) cond inc f}, this requires $f_\Re(u,v)$ to be  monotonic increasing in $v$ for all $u$. Specifically, for fixed $u$, $f_\Re(u,v)$ must satisfy:
  \[
  \frac{\partial f_\Re}{\partial v} > 0 \quad \text{for all } u.
  \]

Finally, for case IV in Eq.~\eqref{eq4: thm: prec w monotone f(x) cond inc f}, we need $f_\Re(u,v)$ independent of the variable $v$ first. We also need $f_\Im(u,v)$ to be monotonic increasing in $v$ for all $u$. Specifically, for fixed $u$, $f_\Im(u,v)$ must satisfy:
  \[
  \frac{\partial f_\Im}{\partial v} > 0 \quad \text{for all } u.
  \]
Therefore, any complex-valued function $f(z)$ that satisfies conditions provided by Eq.~\eqref{eq1: thm: prec w monotone f(x) cond inc f} must be a monotone increasing function.
$\hfill\Box$

Below, we will provide another characterization for the monotonicity under $\prec_w$ by monotone decreasing function $f(z)$ by Theorem~\ref{thm: prec w monotone f(x) cond dec f}.

\begin{theorem}[Majorization preserving by monotone decreasing funciton]\label{thm: prec w monotone f(x) cond dec f}
Given \(\mathbf{x} = [x_1, x_2, \ldots, x_n] \prec_w \mathbf{y} = [y_1, y_2, \ldots, y_n]\), we consider a complex-valued function \(f(z = u + v \iota) = f_{\Re}(u, v) + f_{\Im}(u, v)\), where the partial derivatives of \(f_{\Re}(u, v)\) and \(f_{\Im}(u, v)\) exist. The function \(f(z)\) satisfies the following conditions~\footnote{The pair of real variables are selected from the real part and the imginary part from complex numbers $x_1, x_2, \ldots, x_n, y_1, y_2, \ldots, y_n$ when we apply Lemma~\ref{lma: f(x,y) monotone wrt x at diff y} conditions.}:  
\begin{align}\label{eq1: thm: prec w monotone f(x) cond dec f}
\begin{cases}
f_{\Re} \text{ satisfies conditions from Eq.~\eqref{eq2-1: lma: f(x,y) monotone wrt x at diff y}}, 
& \text{Case I: $z$ increment along } \Re(z); \\
f_{\Re} \text{ is constant, and } f_{\Im} \text{ satisfies conditions from } 
& \text{Eq.~\eqref{eq2-1: lma: f(x,y) monotone wrt x at diff y}}, \\
& \text{Case II: $z$ increment along } \Re(z); \\
\frac{\partial f_{\Re}}{\partial v} < 0, 
& \text{Case III: $z$ stable along } \Re(z) \text{ but increment along } \Im(z); \\
f_{\Re} \text{ is independent of } v, \text{ and } \frac{\partial f_{\Im}}{\partial v} < 0, 
& \text{Case IV: $z$ stable along } \Re(z) \text{ but increment along } \Im(z).
\end{cases}
\end{align}
and
\begin{eqnarray}\label{eq1-1-1: thm: prec w monotone f(x) cond dec f}
f(x_n)&\leq& f(y_n).
\end{eqnarray}
Moreover, the following \emph{difference sum conidtion} for the decreasing function $f$
\begin{eqnarray}\label{eq1-1: thm: prec w monotone f(x) cond dec f}
f(x_i)-f(y_i)&\leq&\sum\limits_{j=i+1}^{n}(f(y_j)-f(x_j)),
\end{eqnarray}
where $i=n-1,\ldots,1$. Then, the function $f$ is $\prec_w$ preserving.
\end{theorem}
\textbf{Proof:}
We first show that any complex-valued function $f(z)$ is monotone decreasing with respect to the variable $z$ will have $\prec_w$ preserving property. 

Given $\mathbf{x}=[x_1,x_2,\ldots,x_n] \prec_w \mathbf{y}=[y_1,y_2,\ldots,y_n]$, we have
\begin{eqnarray}\label{eq2: thm: prec w monotone f(x) cond inc f}
\sum_{i=1}^k x_i \leq \sum_{i=1}^k y_i \quad \text{for all } k = 1, 2, \ldots, n.
\end{eqnarray}
If we consider vectors $\bm{f(x)}=[f(x_1),f(x_2),\ldots,f(x_n)]$ and $\bm{f(y)}=[f(y_1),f(y_2),\ldots,f(y_n)]$, from the monotone decreasing property of the function $f$ with Eq.~\eqref{eq1-1-1: thm: prec w monotone f(x) cond dec f} and the difference sum condition provided by Eq.~\eqref{eq1-1: thm: prec w monotone f(x) cond dec f}, we have
\begin{eqnarray}\label{eq3: thm: prec w monotone f(x) cond dec f}
\sum_{i=1}^k f(x_{n-i+1}) \leq \sum_{i=1}^k f(y_{n-i+1}) \quad \text{for all } k = 1, 2, \ldots, n,
\end{eqnarray}
where $f(x_n) \geq f(x_{n-1}) \geq \ldots \geq f(x_1)$ and $f(y_n) \geq f(y_{n-1}) \geq \ldots \geq f(y_1)$ since $f(z)$ is monotone decreasing with respect to the variable $z$, we have $f(z_1) > f(z_2)$ if $z_1 < z_2$. Therefore, we have $\bm{f(x)}\prec_w \bm{f(y)}$ given $\mathbf{x}\prec_w \mathbf{y}$.

Next, we will demonstrate that any complex-valued funtion $f(z)$ that satifies condition provided by Eq.~\eqref{eq2-1: lma: f(x,y) monotone wrt x at diff y} is a monotone decreasing function. There are four cases for the function $f(z)$ to be a monotone decreasing function with respect to the increasing of the variable $z$: 
\begin{eqnarray}\label{eq4: thm: prec w monotone f(x) cond dec f}
\begin{cases}
f_\Re(u_1, v_1) > f_\Re(u_2,v_2), &~\mbox{if $u_1 < u_2$, case I;}\\
\mbox{or} \\
f_\Re(u_1, v_1) = f_\Re(u_2,v_2), f_\Im(u_1, v_1) > f_\Im(u_2,v_2), &~\mbox{if $u_1 < u_2$, case II;}\\
\mbox{or} \\
f_\Re(u_1, v_1) > f_\Re(u_2,v_2), &~\mbox{if $u_1 = u_2$, $v_1 < v_2$, case III;}\\
\mbox{or} \\
f_\Re(u_1, v_1) = f_\Re(u_2,v_2), f_\Im(u_1, v_1) > f_\Im(u_2,v_2), &~\mbox{if $u_1 = u_2$, $v_1 < v_2$, case IV;}\\
\end{cases}
\end{eqnarray}

For case I in Eq.~\eqref{eq4: thm: prec w monotone f(x) cond dec f}, the conditions provided by Eq.~\eqref{eq2-1: lma: f(x,y) monotone wrt x at diff y} will make $f_\Re(u_1, v_1) > f_\Re(u_2,v_2)$ if $u_1 < u_2$. 

For case II in Eq.~\eqref{eq4: thm: prec w monotone f(x) cond dec f}, we have to require $f_\Re(u,v)$ as a constant since $f_\Re$ is independent of the variable $u$ and $v$. Again, the conditions provided by  Eq.~\eqref{eq2-1: lma: f(x,y) monotone wrt x at diff y} will make $f_\Im(u_1, v_1) > f_\Im(u_2,v_2)$ if $u_1 < u_2$. 

For case III in Eq.~\eqref{eq4: thm: prec w monotone f(x) cond dec f}, this requires $f_\Re(u,v)$ to be  monotonic decreasing in $v$ for all $u$. Specifically, for fixed $u$, $f_\Re(u,v)$ must satisfy:
  \[
  \frac{\partial f_\Re}{\partial v} < 0 \quad \text{for all } u.
  \]

Finally, for case IV in Eq.~\eqref{eq4: thm: prec w monotone f(x) cond inc f}, we need $f_\Re(u,v)$ independent of the variable $v$ first. We also need $f_\Im(u,v)$ to be monotonic decreasing in $v$ for all $u$. Specifically, for fixed $u$, $f_\Im(u,v)$ must satisfy:
  \[
  \frac{\partial f_\Im}{\partial v} < 0 \quad \text{for all } u.
  \]
Therefore, any complex-valued function $f(z)$ that satisfies conditions provided by Eq.~\eqref{eq1: thm: prec w monotone f(x) cond dec f} must be a monotone decreasing function. 
$\hfill\Box$

From Theorem~\ref{thm: prec w monotone f(x) cond inc f} and Theorem~\ref{thm: prec w monotone f(x) cond dec f}, we can have the following Corollary~\ref{cor: prec w monotone no diff sum cond} about $\prec_w$ preserving without difference sum conditions. 
\begin{corollary}\label{cor: prec w monotone no diff sum cond}
Given \(\mathbf{x} = [x_1, x_2, \ldots, x_n] \prec_w \mathbf{y} = [y_1, y_2, \ldots, y_n]\) with $x_i \leq y_i$ for $i=1,2,\ldots,n$, we consider a complex-valued function \(f(z = u + v \iota) = f_{\Re}(u, v) + f_{\Im}(u, v)\), where the partial derivatives of \(f_{\Re}(u, v)\) and \(f_{\Im}(u, v)\) exist. If the function \(f(z)\) satisfies conditions given by Eq.~\eqref{eq1: thm: prec w monotone f(x) cond inc f}, we have $f(x)$ to be $\prec_w$ preserving, i.e., $\bm{f(x)}\prec_w\bm{f(y)}$. 

On the other hand, if f the function \(f(z)\) satisfies conditions given by Eq.~\eqref{eq1: thm: prec w monotone f(x) cond dec f}, we have $\bm{f(y)}\prec_w\bm{f(x)}$. 
\end{corollary}
\textbf{Proof:}
Given $\mathbf{x}=[x_1,x_2,\ldots,x_n] \prec_w \mathbf{y}=[y_1,y_2,\ldots,y_n]$, we have
\begin{eqnarray}\label{eq1: cor: prec w monotone no diff sum cond}
\sum_{i=1}^k x_i \leq \sum_{i=1}^k y_i \quad \text{for all } k = 1, 2, \ldots, n.
\end{eqnarray}
Since $f(z)$ is an increasing function from conditions provided by Eq.~\eqref{eq1: thm: prec w monotone f(x) cond inc f}, we have
\begin{eqnarray}\label{eq2: cor: prec w monotone no diff sum cond}
\sum_{i=1}^k f(x_{i}) \leq \sum_{i=1}^k f(y_{i}) \quad \text{for all } k = 1, 2, \ldots, n,
\end{eqnarray}
where $f(x_1) \geq f(x_2) \geq \ldots \geq f(x_n)$ and $f(y_1) \geq f(y_2) \geq \ldots \geq f(y_n)$. Therefore, we have $\bm{f(x)}\prec_w \bm{f(y)}$ given $\mathbf{x}\prec_w \mathbf{y}$. 

If $f(z)$ is a decreasing function from conditions provided by Eq.~\eqref{eq1: thm: prec w monotone f(x) cond dec f}, we have
\begin{eqnarray}\label{eq3: cor: prec w monotone no diff sum cond}
\sum_{i=1}^k f(y_{n-i+1}) \leq \sum_{i=1}^k f(x_{n-i+1}) \quad \text{for all } k = 1, 2, \ldots, n,
\end{eqnarray}
where $f(x_n) \geq f(x_{n-1}) \geq \ldots \geq f(x_1)$ and $f(y_n) \geq f(y_{n-1}) \geq \ldots \geq f(y_1)$. Therefore, we have $\bm{f(y)}\prec_w \bm{f(x)}$ given $\mathbf{x}\prec_w \mathbf{y}$. 
$\hfill\Box$

\section{Jordan Blocks for Function of a Matrix}\label{sec: Joedan Blocks for Function of a Matrix}

The purpose of this section is to investigate the nilpotent part ordering relationship between the matrix $\bm{X}$ and the matrix $f(\bm{X})$. The dominance order relations for Jordan block sizes with same eigenvalue is discussed in Section~\ref{sec: Generalized Dominance Order Relations for Jordan Block Sizes with Same Eigenvalues}. In Section~\ref{sec: Generalized Dominance Order Relations in Nilpotent Part Matrix Representations R X}, we study the generalized dominance order relations in the nilpotent part of matrix representations, denoted by $\mathfrak{R}(f(\mathbf{X}))$. In Section~\ref{sec: under special f}, we present several corollaries based on the general matrix representation $\mathfrak{R}(f(\bm{X}))$ for various special cases of the function $f$.

\subsection{Generalized Dominance Order Relations for Jordan Block Sizes with Identical Eigenvalues}\label{sec: Generalized Dominance Order Relations for Jordan Block Sizes with Same Eigenvalues}

Let us consider a Jordan block $\bm{J}_{n}(\lambda)$ shown below
\begin{eqnarray}\label{eq: Jordan Block}
\bm{J}_{n}(\lambda)= \begin{bmatrix}
   \lambda & 1 & 0 & \cdots & 0 \\
   0 & \lambda & 1 & \cdots & 0 \\
   \vdots & \vdots & \ddots & \ddots & \vdots \\
   0 & 0 & \cdots & \lambda & 1 \\
   0 & 0 & \cdots & 0 & \lambda
   \end{bmatrix}_{n\times n},
\end{eqnarray}
then, given a function $f$ with the derivative existing to the $n-1$-th ordering, we have $f(\bm{J}_{n}(\lambda))$ expressed as
\begin{eqnarray}\label{eq: f Jordan Block}
f(\bm{J}_{n}(\lambda))= \begin{bmatrix}
   f(\lambda) & \frac{f'(\lambda)}{1!} & \frac{f''(\lambda)}{2!} & \cdots & \frac{f^{(n-1)}(\lambda)}{(n-1)!} \\
   0 & f(\lambda) & \frac{f'(\lambda)}{1!} & \cdots & \frac{f^{(n-2)}(\lambda)}{(n-2)!}  \\
   \vdots & \vdots & \ddots & \ddots & \vdots \\
   0 & 0 & \cdots & f(\lambda) & \frac{f'(\lambda)}{1!} \\
   0 & 0 & \cdots & 0 & f(\lambda)
   \end{bmatrix}_{n\times n},
\end{eqnarray}

We have the following Lemma~\ref{lma: GDO Distance single Jordan block} about generalized dominance ordering distance between the matrix $\bm{J}_{n}(\lambda)$ and $f(\bm{J}_{n}(\lambda))$.
\begin{lemma}\label{lma: GDO Distance single Jordan block}
Given a Jordan block provided by Eq.~\eqref{eq: Jordan Block} and a function $f$ with the derivative existing to the $n-1$-th order. If $f^{(\kappa)}(\lambda)\neq 0$ and $f^{(i)}(\lambda) = 0$ for $i=1,2,\ldots,\kappa-1$, then, we have $\bm{m}_1(f(\bm{J}_{n}(\lambda))) \trianglelefteq \bm{m}_1(\bm{J}_{n}(\lambda))$, where $\bm{m}_1(\bm{J}_{n}(\lambda))=[n]$ as there are only one eigenvalue for the Jordan block $\bm{J}_n(\lambda)$ and the generalized dominance ordering distance between $\bm{m}_1(f(\bm{J}_{n}(\lambda)))$ and $\bm{m}_1(\bm{J}_{n}(\lambda))$, denoted by \\
$\mathfrak{D}_{\bm{m}_1(f(\bm{J}_{n}(\lambda))),\bm{m}_1(\bm{J}_{n}(\lambda))}(j)$, is 
\begin{eqnarray}\label{eq1:lma:GDO_Distance_single_Jordan_block}
\mathfrak{D}_{\bm{m}_1(f(\bm{J}_{n}(\lambda))),\bm{m}_1(\bm{J}_{n}(\lambda))}(j) =
\begin{cases} 
n - j\lceil\frac{n}{\kappa} \rceil, & \text{if } j \leq \ell; \\
n+j-\ell-j \lceil\frac{n}{\kappa} \rceil, & \text{if } \ell < j \leq \kappa.
\end{cases}
\end{eqnarray}
where $\ell$ is determined by:
\begin{eqnarray}\label{eq2: lma: GDO Distance single Jordan block}
\ell=n+\kappa-\kappa \left\lceil \frac{n}{\kappa}\right\rceil. 
\end{eqnarray}
\end{lemma}
\textbf{Proof:}
For the matrix $f(\bm{J}_{n}(\lambda))$ given by Eq.~\eqref{eq: f Jordan Block}, we know the algebraic multiplicity for the matrix $f(\bm{J}_{n}(\lambda))$ is $n$. Besides, the rank of the following matrix 
\begin{eqnarray}\label{eq2-0: lma: GDO Distance single Jordan block}
f(\bm{J}_{n}(\lambda))- f(\lambda)\bm{I}= \begin{bmatrix}
   0 & \frac{f'(\lambda)}{1!} & \frac{f''(\lambda)}{2!} & \cdots & \frac{f^{(n-1)}(\lambda)}{(n-1)!} \\
   0 & 0 & \frac{f'(\lambda)}{1!} & \cdots & \frac{f^{(n-2)}(\lambda)}{(n-2)!}  \\
   \vdots & \vdots & \ddots & \ddots & \vdots \\
   0 & 0 & \cdots & 0 & \frac{f'(\lambda)}{1!} \\
   0 & 0 & \cdots & 0 & 0
   \end{bmatrix}_{n\times n},
\end{eqnarray}
is $n-\kappa$ given that $f^{(\kappa)}\neq 0$ and $f^{(i)} = 0$ for $i=1,2,\ldots,\kappa-1$. Therefore, there are $\kappa$ blocks of the Jordan form of the following matrix: 
\begin{eqnarray}\label{eq2-1: lma: GDO Distance single Jordan block}
f(\bm{J}_{n}(\lambda))- f(\lambda)\bm{I}= \begin{bmatrix}
   0 & \ldots & 0 & \frac{f^{(\kappa)}(\lambda)}{\kappa!} & \frac{f^{(\kappa+1)}(\lambda)}{(\kappa+1)!} & \cdots & \frac{f^{(n-1)}(\lambda)}{(n-1)!} \\
   0 & \ldots & 0 & 0 & \frac{f^{(\kappa)}(\lambda)}{\kappa!} & \cdots & \frac{f^{(n-2)}(\lambda)}{(n-2)!}  \\
   \vdots & \vdots & \vdots  & \vdots & \ddots & \ddots & \vdots \\
   0 & \ldots & 0 & 0 & \cdots & 0 & \frac{f^{(\kappa)}(\lambda)}{\kappa!} \\
   0 & \ldots & 0 & 0 & \cdots & 0 & 0 \\
   \vdots & \ldots & 0 & \vdots & \cdots & \vdots &\vdots \\
   0 & \ldots & 0 & 0 & \cdots & 0 & 0
\end{bmatrix}_{n\times n} \define \bm{M}_{f,\kappa},
\end{eqnarray}
where $f^{(\kappa)}(\lambda)$ is the $\kappa$-th derivative at the value $\lambda$. Because the matrix $\bm{M}_{f,\kappa}$ is a nilpotnt matrix, from the non-zero entries structure given by Eq.~\eqref{eq2-1: lma: GDO Distance single Jordan block}, we have 
\begin{eqnarray}\label{eq2-2: lma: GDO Distance single Jordan block}
\begin{cases} 
\bm{M}^{p}_{f,\kappa} = \bm{0}, & \text{if } p \geq \lceil \frac{n}{\kappa} \rceil; \\
\bm{M}^{p}_{f,\kappa} \neq \bm{0}, & \text{if } p < \lceil \frac{n}{\kappa} \rceil.
\end{cases}
\end{eqnarray}
Let $\ell$ be the number sub-Jordan blocks of $\bm{M}_{f,\kappa}$ with size $\lceil \frac{n}{\kappa} \rceil$, the remaining $\kappa-\ell$ sub-Jordan blocks of $\bm{M}_{f,\kappa}$ will have the size $\lceil \frac{n}{\kappa} \rceil-1$. Then, we have
\begin{eqnarray}\label{eq2-3: lma: GDO Distance single Jordan block}
\ell \left\lceil \frac{n}{\kappa} \right\rceil + (\kappa-\ell)\left(\left\lceil \frac{n}{\kappa}\right\rceil-1\right)&=&n.
\end{eqnarray}

By applying the definition of $\mathfrak{D}_{\bm{p},\bm{q}}(j)$ (GDOD) provided by Eq.~\eqref{eq: dom order distance} with the value $\ell$ given by Eq.~\eqref{eq2-3: lma: GDO Distance single Jordan block}, this lemma is proved. 
$\hfill\Box$

The $\mathfrak{D}_{\bm{m}_1(f(\bm{J}_{n}(\lambda))),\bm{m}_1(\bm{J}_{n}(\lambda))}(j)$ given by Lemma~\ref{lma: GDO Distance single Jordan block} is only appliable for one geometry multiplicity Jordan block. Following Lemma~\ref{lma: GDO Distance 2 Jordan block same e-value} is to determine the GDOD for two Jordan blocks sharing same eigenvalues. 

\begin{lemma}\label{lma: GDO Distance 2 Jordan block same e-value}
If we are given two natural numbers with $n_1 \geq n_2$, we have the following evaluations for $\mathfrak{D}_{\bm{m}_1(f(\bm{J}_{n_1}(\lambda)\bigoplus\bm{J}_{n_2}(\lambda))),\bm{m}_1(\bm{J}_{n_1}(\lambda)\bigoplus\bm{J}_{n_2}(\lambda))}(j)$ with respect to different $n_1, n_2$ and $\kappa$, where $\kappa$ is the smallest integer to have $f^{(\kappa)} \neq 0$:
\begin{enumerate}
\item Case I: $\lceil \frac{n_1}{\kappa} \rceil = \lceil \frac{n_2}{\kappa} \rceil$
\begin{eqnarray}\label{eq1: lma: GDO Distance 2 Jordan block same e-value}
\lefteqn{\mathfrak{D}_{\bm{m}_1(f(\bm{J}_{n_1}(\lambda)\bigoplus\bm{J}_{n_2}(\lambda))),\bm{m}_1(\bm{J}_{n_1}(\lambda)\bigoplus\bm{J}_{n_2}(\lambda))}(j)=}\nonumber \\
&&
\begin{cases} 
n_1 - \lceil \frac{n_1}{\kappa} \rceil, & \mbox{if $j=1$}; \\
n_1+n_2 - j \lceil \frac{n_1}{\kappa} \rceil, & \mbox{if $1 < j \leq \ell_1 + \ell_2$}; \\
n_1+n_2 - (\ell_1+\ell_2)\lceil \frac{n_1}{\kappa} \rceil -  (j-\ell_1-\ell_2)(\lceil \frac{n_1}{\kappa} \rceil-1), & \mbox{if $\ell_1 + \ell_2 < j \leq 2\kappa$}; \\
\end{cases}
\end{eqnarray}
where $\ell_1=n_1+\kappa-\kappa \left\lceil \frac{n_1}{\kappa}\right\rceil$ and $\ell_2=n_2+\kappa-\kappa \left\lceil \frac{n_2}{\kappa}\right\rceil$ from Eq.~\eqref{eq2: lma: GDO Distance single Jordan block}.

\item Case II: $\lceil \frac{n_1}{\kappa} \rceil-1 = \lceil \frac{n_2}{\kappa} \rceil$
\begin{eqnarray}\label{eq2: lma: GDO Distance 2 Jordan block same e-value}
\lefteqn{\mathfrak{D}_{\bm{m}_1(f(\bm{J}_{n_1}(\lambda)\bigoplus\bm{J}_{n_2}(\lambda))),\bm{m}_1(\bm{J}_{n_1}(\lambda)\bigoplus\bm{J}_{n_2}(\lambda))}(j)=}\nonumber \\
&&
\begin{cases} 
n_1 - \lceil \frac{n_1}{\kappa} \rceil, & \mbox{if $j=1$}; \\
n_1+n_2 - j \lceil \frac{n_1}{\kappa} \rceil, & \mbox{if $1 < j \leq \ell_1$}; \\
n_1+n_2 - \ell_1 \lceil \frac{n_1}{\kappa} \rceil - (j-\ell_1)\lceil \frac{n_2}{\kappa} \rceil, & \mbox{if $\ell_1 < j \leq \kappa+\ell_2$}; \\
n_1+n_2 - \ell_1\lceil \frac{n_1}{\kappa} \rceil - (\kappa-\ell_1+\ell_2)\lceil \frac{n_2}{\kappa} \rceil \nonumber \\
- (j-\kappa - \ell_2)(\lceil \frac{n_2}{\kappa} \rceil-1), & \mbox{if $\kappa+ \ell_2 < j \leq 2\kappa$}. \\
\end{cases}
\end{eqnarray}

\item Case III: $\lceil \frac{n_1}{\kappa} \rceil-1 > \lceil \frac{n_2}{\kappa} \rceil$
\begin{eqnarray}\label{eq3: lma: GDO Distance 2 Jordan block same e-value}
\lefteqn{\mathfrak{D}_{\bm{m}_1(f(\bm{J}_{n_1}(\lambda)\bigoplus\bm{J}_{n_2}(\lambda))),\bm{m}_1(\bm{J}_{n_1}(\lambda)\bigoplus\bm{J}_{n_2}(\lambda))}(j)=}\nonumber \\
&&
\begin{cases} 
n_1 - \lceil \frac{n_1}{\kappa} \rceil, & \mbox{if $j=1$}; \\
n_1+n_2 - j \lceil \frac{n_1}{\kappa} \rceil, & \mbox{if $1 < j \leq \ell_1$}; \\
n_1+n_2 - \ell_1 \lceil \frac{n_1}{\kappa} \rceil - (j-\ell_1)(\lceil \frac{n_1}{\kappa} \rceil-1), & \mbox{if $\ell_1 < j \leq \kappa$}; \\
n_1+n_2 - \ell_1\lceil \frac{n_1}{\kappa} \rceil - (\kappa-\ell_1)(\lceil \frac{n_1}{\kappa} \rceil-1) \nonumber \\
- (j-\kappa)\lceil \frac{n_2}{\kappa} \rceil, & \mbox{if $\kappa < j \leq \kappa+\ell_2$};\\
n_1+n_2 - \ell_1\lceil \frac{n_1}{\kappa} \rceil - (\kappa-\ell_1)(\lceil \frac{n_1}{\kappa} \rceil-1) \nonumber \\
- \ell_2 \lceil \frac{n_2}{\kappa} \rceil - (j-\kappa-\ell_2)(\lceil \frac{n_2}{\kappa} \rceil-1), & \mbox{if $\kappa+\ell_2 < j \leq 2\kappa$}.\\
\end{cases}
\end{eqnarray}
\end{enumerate}

\end{lemma}
\textbf{Proof:}
We will prove for the Case I only. The remaining cases can be proved similarly. From Lemma~\ref{lma: GDO Distance single Jordan block}, we have 
\begin{eqnarray}\label{eq4: lma: GDO Distance 2 Jordan block same e-value}
\bm{m}_1(f(\bm{J}_{n_1}(\lambda))&=&[\underbrace{\lceil\frac{n_1}{\kappa}\rceil, \lceil\frac{n_1}{\kappa}\rceil, \ldots, \lceil\frac{n_1}{\kappa}\rceil}_{\mbox{$\ell_1$ terms}}, \underbrace{\lceil\frac{n_1}{\kappa}\rceil - 1, \lceil\frac{n_1}{\kappa}\rceil - 1, \ldots, \lceil\frac{n_1}{\kappa}\rceil - 1}_{\mbox{$\kappa-\ell_1$ terms}}],
\end{eqnarray}
where $\ell_1=n_1+\kappa-\kappa \left\lceil \frac{n_1}{\kappa}\right\rceil$. Again, by Lemma~\ref{lma: GDO Distance single Jordan block}, we also have 
\begin{eqnarray}\label{eq5: lma: GDO Distance 2 Jordan block same e-value}
\bm{m}_1(f(\bm{J}_{n_2}(\lambda))&=&[\underbrace{\lceil\frac{n_2}{\kappa}\rceil, \lceil\frac{n_2}{\kappa}\rceil, \ldots, \lceil\frac{n_2}{\kappa}\rceil}_{\mbox{$\ell_2$ terms}}, \underbrace{\lceil\frac{n_2}{\kappa}\rceil - 1, \lceil\frac{n_2}{\kappa}\rceil - 1, \ldots, \lceil\frac{n_2}{\kappa}\rceil - 1}_{\mbox{$\kappa-\ell_2$ terms}}],
\end{eqnarray}
where $\ell_2=n_2+\kappa-\kappa \left\lceil \frac{n_2}{\kappa}\right\rceil$.  

Since $\lceil \frac{n_1}{\kappa} \rceil = \lceil \frac{n_2}{\kappa} \rceil$, we have
\begin{eqnarray}\label{eq5: lma: GDO Distance 2 Jordan block same e-value}
\bm{m}_1(f(\bm{J}_{n_1}(\lambda)\bigoplus\bm{J}_{n_2}(\lambda))&=&\bm{m}_1(f(\bm{J}_{n_1}(\lambda))\bigoplus(\bm{J}_{n_2}(\lambda))\nonumber\\
&=&
[\underbrace{\lceil\frac{n_1}{\kappa}\rceil, \lceil\frac{n_1}{\kappa}\rceil, \ldots, \lceil\frac{n_1}{\kappa}\rceil}_{\mbox{$\ell_1$ terms}}, \underbrace{\lceil\frac{n_2}{\kappa}\rceil, \lceil\frac{n_2}{\kappa}\rceil, \ldots, \lceil\frac{n_2}{\kappa}\rceil}_{\mbox{$\ell_2$ terms}}, \nonumber \\
&& \underbrace{\lceil\frac{n_1}{\kappa}\rceil - 1, \lceil\frac{n_1}{\kappa}\rceil - 1, \ldots, \lceil\frac{n_1}{\kappa}\rceil - 1}_{\mbox{$\kappa-\ell_1$ terms}}, \nonumber \\
&& \underbrace{\lceil\frac{n_2}{\kappa}\rceil - 1, \lceil\frac{n_2}{\kappa}\rceil - 1, \ldots, \lceil\frac{n_2}{\kappa}\rceil - 1}_{\mbox{$\kappa-\ell_2$ terms}}],
\end{eqnarray}
This case is proved by applying the generalized dominance ordering distance defined by Eq.~\eqref{eq: dom order distance} to vectors $\bm{m}_1(f(\bm{J}_{n_1}(\lambda)\bigoplus\bm{J}_{n_2}(\lambda))$ and $\bm{m}_1(\bm{J}_{n_1}(\lambda)\bigoplus\bm{J}_{n_2}(\lambda))$.
$\hfill\Box$

Given the matrix $\bigoplus\limits_{i=1}^{\alpha_{1}^{G}}\bm{J}_{m_{1,i}}(\lambda)$ with the dimension $\alpha_1^{A} \times \alpha_1^{A}$, we can express it as
\begin{eqnarray}
\raisebox{10.5em}{\(\bigoplus\limits_{i=1}^{\alpha_{1}^{G}}\bm{J}_{m_{1,i}}(\lambda)\)=}&&
\begin{tikzpicture}
    \matrix[matrix of math nodes,left delimiter={(},right delimiter={)}] (m) {
      \lambda & 1        & 0   & \ldots   & 0      & 0     & 0      & 0  & 0      & 0  & 0 & 0\\
        0       &  \lambda & 1    & 0 \ldots    & 0     & 0      & 0  & 0      & 0 & 0  & 0 & 0 \\
        0      &  0 & \ddots &   \ddots   & 0      & 0     & 0      & 0  & 0      & 0  & 0 & 0  \\
        0      & \ldots   & 0     & \lambda   &  1   & 0     & 0      & 0  & 0      & 0   & 0  & 0 \\
        0      & 0      & \ldots      &  0 &  \lambda     & 0      & 0  & 0      & 0 & 0  & 0  & 0 \\
        0      & 0      & 0      & 0      &  \ddots   & \ddots     & \ddots    & 0      & 0 & 0  & 0  & 0\\
        0      & 0    & 0   & 0      & 0      &  \ddots   & \ddots       & \ddots   & 0      & 0 & 0  & 0  \\
        0      & 0    & 0   & 0      & 0      & 0     & 0      & \lambda  & 1      & 0  & \ldots  & 0 \\
        0      & 0    & 0  & 0      & 0      & 0     & 0     &0     &  \lambda & 1    & 0 \ldots    & 0  \\
        0      & 0    & 0  & 0      & 0      & 0     & 0      &  0      &  0 & \ddots &   \ddots   & 0 \\
        0      & 0    & 0  & 0      & 0      & 0     & 0      &  0      & \ldots   & 0     & \lambda   &  1 \\
        0      & 0    & 0  & 0      & 0      & 0     & 0      &  0      & 0      & \ldots      &  0 &  \lambda\\
    };

    \draw[dashed] 
        (m-1-1.north west) -- (m-1-5.north east) -- 
        (m-5-5.south east) -- (m-5-1.south west) -- cycle;
\node[anchor=north] at ($(m-5-5.south east) + (0.5,0)$) {\tiny $m_{1,1} \times m_{1,1}$};

    \draw[dashed] 
        (m-8-8.north west) -- (m-8-12.north east) -- 
        (m-12-12.south east) -- (m-12-8.south west) -- cycle;
\node[anchor=north] at ($(m-12-12.south east) + (0.5,0)$){\tiny $m_{1,\alpha_1^{G}} \times m_{1,\alpha_1^{G}}$};

\end{tikzpicture}
\end{eqnarray}
where we have $\alpha_1^{A}=\sum\limits_{i=1}^{\alpha_{1}^{G}}m_{1,i}$.  Following Lemma~\ref{lma: GDO Distance multiple Jordan block same e-value} will derive \\
$\mathfrak{D}_{\bm{m}_1\left(f\left(\bigoplus\limits_{i=1}^{\alpha_1^{(\mathrm{G})}}\bm{J}_{m_{1,i}}(\lambda)\right)\right), \bm{m}_1\left(\bigoplus\limits_{i=1}^{\alpha_1^{(\mathrm{G})}}\bm{J}_{m_{1,i}}(\lambda)\right)}(j)$ for all Jordan blocks sharing the same eigenvalue. 

Before presenting Lemma~\ref{lma: GDO Distance multiple Jordan block same e-value}, we have to present a new operation, denoted by $\myop$, between two vectors with sorted entries to perform merge sorting. Merge sort can efficiently merge two vectors with ordered entries into a single sorted vector. Let us denote the two vectors as:

\[
\bm{l} = [l_1, l_2, \dots, l_m], \quad \bm{r} = [r_1, r_2, \dots, r_n],
\]

where the entries of lists $l$ and $r$ are sorted in non-increasing order:
\[
l_1 \geq l_2 \geq \dots \geq l_m, \quad r_1 \geq r_2 \geq \dots \geq r_n.
\]

The goal is to produce a single sorted vector:
\[
\bm{c} = [c_1, c_2, \dots, c_{m+n}],
\]
where $\bm{c}$ contains all the elements of $l$ and $r$ in sorted order. Then, $\bm{c}$ can be expressed by the operation $\myop$ as
\[
\bm{c} = \bm{l} \myop \bm{r}.
\]

The merging process works by comparing the elements from $L$ and $R$ one at a time, and inserting the smallest element into $C$. This process is as follows:

\begin{enumerate}
    \item Initialize three indices:
    \[
    i = 1 \quad (\text{index for } \bm{l}), \quad j = 1 \quad (\text{index for } \bm{r}), \quad k = 1 \quad (\text{index for } \bm{c}).
    \]

    \item Compare $l_i$ (the current element of $\bm{l}$) with $r_j$ (the current element of $\bm{r}$):
    \[
    \text{If } l_i \leq r_j, \text{ set } c_k = l_i \text{ and increment } i \text{ and } k.
    \]
    \[
    \text{Otherwise, set } c_k = r_j \text{ and increment } j \text{ and } k.
    \]

    \item Repeat the comparison until one of the vectors $\bm{l}$ or $\bm{r}$ is exhausted.

    \item Append any remaining elements from the non-exhausted vector to $\bm{c}$.
\end{enumerate}

\begin{lemma}\label{lma: GDO Distance multiple Jordan block same e-value} 
We consider a matrix with Jordan decomposition as $\bigoplus\limits_{i=1}^{\alpha_1^{(\mathrm{G})}}\bm{J}_{m_{1,i}}(\lambda)$ and its functional matrix  $f\left(\bigoplus\limits_{i=1}^{\alpha_1^{(\mathrm{G})}}\bm{J}_{m_{1,i}}(\lambda)\right)$, where the function $f^{(\kappa)}(\lambda) \neq 0$, and $\kappa$ is the smallest integer such that $f^{(\kappa)}(\lambda) \neq 0$. Then, we have
\begin{eqnarray}\label{eq1: lma: GDO Distance multiple Jordan block same e-value} 
\mathfrak{D}_{\bm{m}_1\left(f\left(\bigoplus\limits_{i=1}^{\alpha_1^{(\mathrm{G})}}\bm{J}_{m_{1,i}}(\lambda)\right)\right), \bm{m}_1\left(\bigoplus\limits_{i=1}^{\alpha_1^{(\mathrm{G})}}\bm{J}_{m_{1,i}}(\lambda)\right)}(j)&=& \begin{cases} 
\sum\limits_{i=1}^j m_{1.i} -\sum\limits_{i=1}^j \eta_{f,\lambda}(i), & \mbox{if $j \leq \alpha_1^{(\mathrm{G})}$}; \\
\alpha_1^{(\mathrm{A})} -\sum\limits_{i=1}^j \eta_{f,\lambda}(i), & \mbox{if $\alpha_1^{(\mathrm{G})} < j \leq \alpha_1^{(\mathrm{G})} \kappa$}. \\
\end{cases}
\end{eqnarray}
where $\eta_{f,\lambda}(i)$ is the $i$-th entry of the vector $\bm{\eta}_{f,\lambda}$ obtained by merge sorting entries from $\alpha_1^{(G)}$ vectors:
\begin{eqnarray}\label{eq2: lma: GDO Distance multiple Jordan block same e-value} 
\bm{\eta}_{f,\lambda}&=&\myop_{i=1}^{\alpha_1^{G}}[ (m_{1,i}-\mathfrak{D}_{\bm{m}_1(f(\bm{J}_{m_{1,i}}(\lambda))),\bm{m}_1(\bm{J}_{m_{1,i}}(\lambda))}(1)),  \nonumber \\
&& (\mathfrak{D}_{\bm{m}_1(f(\bm{J}_{m_{1,i}}(\lambda))),\bm{m}_1(\bm{J}_{m_{1,i}}(\lambda))}(1) - \mathfrak{D}_{\bm{m}_1(f(\bm{J}_{m_{1,i}}(\lambda))),\bm{m}_1(\bm{J}_{m_{1,i}}(\lambda))}(2)),\cdots,\nonumber \\
&& (\mathfrak{D}_{\bm{m}_1(f(\bm{J}_{m_{1,i}}(\lambda))),\bm{m}_1(\bm{J}_{m_{1,i}}(\lambda))}(\kappa-1) - \mathfrak{D}_{\bm{m}_1(f(\bm{J}_{m_{1,i}}(\lambda))),\bm{m}_1(\bm{J}_{m_{1,i}}(\lambda))}(\kappa))
]
\end{eqnarray}
\end{lemma}
\textbf{Proof:}
Without loss of generality, we may assume that numbers having $m_{1,1} \geq m_{1,2} \geq \ldots \geq m_{1,\alpha_1^{(G)}}$, which are entries of the vector $\bm{m}_1\left(\bigoplus\limits_{i=1}^{\alpha_1^{(\mathrm{G})}}\bm{J}_{m_{1,i}}(\lambda)\right)$. 

Because $f\left(\bigoplus\limits_{i=1}^{\alpha_1^{(\mathrm{G})}}\bm{J}_{m_{1,i}}(\lambda)\right)=\bigoplus \limits_{i=1}^{\alpha_1^{(\mathrm{G})}}f(\bm{J}_{m_{1,i}})$ from Jordan block structure, we have 
\begin{eqnarray}\label{eq3: lma: GDO Distance multiple Jordan block same e-value} 
\bm{m}_1\left(f\left(\bigoplus\limits_{i=1}^{\alpha_1^{(\mathrm{G})}}\bm{J}_{m_{1,i}}(\lambda)\right)\right)&=&\bm{m}_1\left(\bigoplus \limits_{i=1}^{\alpha_1^{(\mathrm{G})}}f(\bm{J}_{m_{1,i}})\right).
\end{eqnarray}
Also, from Lemma~\ref{lma: GDO Distance single Jordan block}, Jordan block sizes of $f(\bm{J}_{m_{1,i}}(\lambda))$ (if we order non-increasing by size) is
\begin{eqnarray}\label{eq4: lma: GDO Distance multiple Jordan block same e-value} 
[ (m_{1,i}-\mathfrak{D}_{\bm{m}_1(f(\bm{J}_{m_{1,i}}(\lambda))),\bm{m}_1(\bm{J}_{m_{1,i}}(\lambda))}(1)),  ~~~~~~~~~~~~~~~~~~~~~~~~~~~~~~~~~~~~~~~~~~~~~~~~~~~~~\nonumber \\
 (\mathfrak{D}_{\bm{m}_1(f(\bm{J}_{m_{1,i}}(\lambda))),\bm{m}_1(\bm{J}_{m_{1,i}}(\lambda))}(1) - \mathfrak{D}_{\bm{m}_1(f(\bm{J}_{m_{1,i}}(\lambda))),\bm{m}_1(\bm{J}_{m_{1,i}}(\lambda))}(2)),\cdots,\nonumber \\
(\mathfrak{D}_{\bm{m}_1(f(\bm{J}_{m_{1,i}}(\lambda))),\bm{m}_1(\bm{J}_{m_{1,i}}(\lambda))}(\kappa-1) - \mathfrak{D}_{\bm{m}_1(f(\bm{J}_{m_{1,i}}(\lambda))),\bm{m}_1(\bm{J}_{m_{1,i}}(\lambda))}(\kappa))
].
\end{eqnarray}
Then, we have $\bm{m}_1\left(f\left(\bigoplus\limits_{i=1}^{\alpha_1^{(\mathrm{G})}}\bm{J}_{m_{1,i}}(\lambda)\right)\right) = \bm{\eta}_{f,\lambda}$ from Eq.~\eqref{eq2: lma: GDO Distance multiple Jordan block same e-value}. 

This lemma is proved by applying the generalized dominance ordering distance defined by Eq.~\eqref{eq: dom order distance} to vectors $\bm{m}_1\left(f\left(\bigoplus\limits_{i=1}^{\alpha_1^{(\mathrm{G})}}\bm{J}_{m_{1,i}}(\lambda)\right)\right)$ and $\bm{m}_1\left(\bigoplus\limits_{i=1}^{\alpha_1^{(\mathrm{G})}}\bm{J}_{m_{1,i}}(\lambda)\right)$.
$\hfill\Box$

\begin{remark}
The vector $\bm{\eta}_{f,\lambda}$ is the Jordan sublock sizes that sharing the same eigenvalues ordered by non-decreasing order.
\end{remark}

Below, we will present Example~\ref{exp: lma: GDO Distance multiple Jordan block same e-value} to illustare Lemma~\ref{lma: GDO Distance multiple Jordan block same e-value}. 

\begin{example}\label{exp: lma: GDO Distance multiple Jordan block same e-value}
Let us consider the following matrix
\begin{eqnarray}\label{eq1: exp: lma: GDO Distance multiple Jordan block same e-value}
\bm{J}_{4}(\lambda)\bigoplus\bm{J}_{3}(\lambda)\bigoplus\bm{J}_{2}(\lambda) =
\begin{tikzpicture}[baseline=(m.center)]
    \matrix[matrix of math nodes,left delimiter={(},right delimiter={)},row sep=0.8em, column sep=0.8em] (m) {
      \lambda & 1 & 0 & 0 & 0 & 0 & 0 & 0 & 0 \\
      0  & \lambda & 1 & 0 & 0 & 0 & 0 & 0 & 0 \\
      0   & 0 & \lambda & 1 & 0 & 0 & 0 & 0 & 0  \\
      0   & 0  &   0 & \lambda & 0 & 0 & 0 & 0 & 0  \\
      0  &   0 & 0 & 0  & \lambda & 1 & 0 & 0 & 0  \\
      0 & 0 & 0 & 0 & 0 & \lambda & 1 & 0 & 0 \\
      0  & 0 & 0 & 0 & 0 & 0 & \lambda & 0 & 0 \\
      0 & 0 & 0 & 0 & 0 & 0 & 0 & \lambda & 1 \\
      0 & 0 & 0 & 0 & 0 & 0 & 0 & 0 & \lambda \\
    };

    \draw[dashed] (m-4-1.south west) -- (m-4-4.south east); 
    \draw[dashed] (m-1-4.north east) -- (m-4-4.south east); 

    \draw[dashed] (m-5-5.north west) -- (m-5-7.north east); 
    \draw[dashed] (m-7-5.south west) -- (m-7-7.south east); 
    \draw[dashed] (m-5-5.north west) -- (m-7-5.south west); 
    \draw[dashed] (m-5-7.north east) -- (m-7-7.south east); 

    \draw[dashed] (m-8-8.north west) -- (m-8-9.north east); 
    \draw[dashed] (m-9-8.south west) -- (m-9-9.south east); 
    \draw[dashed] (m-8-8.north west) -- (m-9-8.south west); 
    \draw[dashed] (m-8-9.north east) -- (m-9-9.south east); 
\end{tikzpicture}_{9\times 9},
\end{eqnarray}

Then, we have
\begin{eqnarray}\label{eq2: exp: lma: GDO Distance multiple Jordan block same e-value}
\lefteqn{f(\bm{J}_{4}(\lambda)\bigoplus\bm{J}_{3}(\lambda)\bigoplus\bm{J}_{2}(\lambda))=}\nonumber \\
&&
\begin{tikzpicture}[baseline=(m.center)]
    \matrix[matrix of math nodes,left delimiter={(},right delimiter={)},row sep=0.8em, column sep=0.8em] (m) {
      f(\lambda) & f'(\lambda) & f''(\lambda)/2 & f'''(\lambda)/3! & 0 & 0 & 0 & 0 & 0 \\
      0  & f(\lambda)  & f'(\lambda) & f''(\lambda)/2 & 0 & 0 & 0 & 0 & 0 \\
      0   & 0 & f(\lambda)  & f'(\lambda) & 0 & 0 & 0 & 0 & 0  \\
      0   & 0  &   0 & f(\lambda)  & 0 & 0 & 0 & 0 & 0  \\
      0  &   0 & 0 & 0  & f(\lambda)  & f'(\lambda) & f''(\lambda)/2  & 0 & 0  \\
      0 & 0 & 0 & 0 & 0 & f(\lambda)  & f'(\lambda) & 0 & 0 \\
      0  & 0 & 0 & 0 & 0 & 0 & f(\lambda) & 0 & 0 \\
      0 & 0 & 0 & 0 & 0 & 0 & 0 & f(\lambda)  & f'(\lambda) \\
      0 & 0 & 0 & 0 & 0 & 0 & 0 & 0 & f(\lambda)  \\
    };

    \draw[dashed] (m-4-1.south west) -- (m-4-4.south east); 
    \draw[dashed] (m-1-4.north east) -- (m-4-4.south east); 

    \draw[dashed] (m-5-5.north west) -- (m-5-7.north east); 
    \draw[dashed] (m-7-5.south west) -- (m-7-7.south east); 
    \draw[dashed] (m-5-5.north west) -- (m-7-5.south west); 
    \draw[dashed] (m-5-7.north east) -- (m-7-7.south east); 

    \draw[dashed] (m-8-8.north west) -- (m-8-9.north east); 
    \draw[dashed] (m-9-8.south west) -- (m-9-9.south east); 
    \draw[dashed] (m-8-8.north west) -- (m-9-8.south west); 
    \draw[dashed] (m-8-9.north east) -- (m-9-9.south east); 
\end{tikzpicture}_{9\times 9},
\end{eqnarray}

If we have $f'(\lambda)=0$, we have the following Jordan block sizes:
\begin{eqnarray}\label{eq3: exp: lma: GDO Distance multiple Jordan block same e-value}
\mbox{for $f(\bm{J}_4(\lambda))$}&\rightarrow&[2,2], \nonumber \\
\mbox{for $f(\bm{J}_3(\lambda))$}&\rightarrow&[2,1], \nonumber \\
\mbox{for $f(\bm{J}_2(\lambda))$}&\rightarrow&[1,1], \nonumber \\
\end{eqnarray}
then, we have 
\begin{eqnarray}\label{eq4: exp: lma: GDO Distance multiple Jordan block same e-value}
\bm{\eta}_{f,\lambda}&=&[2,2]\myop[2,1]\myop[1,1]\nonumber \\
&=&[2,2,2,1,1,1],
\end{eqnarray}
Moreover, we also have
\begin{eqnarray}\label{eq4: exp: lma: GDO Distance multiple Jordan block same e-value}
\lefteqn{\mathfrak{D}_{\bm{m}_1\left(f(\bm{J}_{4}(\lambda)\bigoplus\bm{J}_{3}(\lambda)\bigoplus\bm{J}_{2}(\lambda))\right), \bm{m}_1\left(\bm{J}_{4}(\lambda)\bigoplus\bm{J}_{3}(\lambda)\bigoplus\bm{J}_{2}(\lambda)\right)}(j)=}\nonumber \\
&& \begin{cases} 
\sum\limits_{i=1}^j m_{1.i} -\sum\limits_{i=1}^j \eta_{f,\lambda}(i), & \mbox{if $j \leq 3$}; \\
(4+3+2) -\sum\limits_{i=1}^j \eta_{f,\lambda}(i), & \mbox{if $3 < j \leq 3 \times 2$}, \\
\end{cases}
\end{eqnarray}
where $m_{1,1}=4, m_{1,2}=3$ and $m_{1,3}=2$.
\end{example}

\subsection{Matrix Representations $\mathfrak{R}(f(\mathbf{X}))$ and GDOD of Nilpotent Part of $\mathfrak{R}(f(\mathbf{X}))$}\label{sec: Generalized Dominance Order Relations in Nilpotent Part Matrix Representations R X}

In this section, we will study the generalized dominance order relations in nilpotent part matrix representations $\mathfrak{R}(f(\mathbf{X}))$ and discuss some special cases regarding these order relations with respect to properties of $f$ and $\bm{X}$.

Recall Theorem 1 in~\cite{chang2024operatorCH}, we have the following spectral mapping theorem. 
\begin{theorem}\label{thm: Spectral Mapping Theorem for Single Variable}
Given an analytic function $f(z)$ within the domain for $|z| < R$, a matrix $\bm{X}$ with the dimension $m$ and $K$ distinct eigenvalues $\lambda_k$ for $k=1,2,\ldots,K$ such that
\begin{eqnarray}\label{eq1: thm: Spectral Mapping Theorem for Single Variable}
\bm{X}&=&\sum\limits_{k=1}^K\sum\limits_{i=1}^{\alpha_k^{\mathrm{G}}} \lambda_k \bm{P}_{k,i}+
\sum\limits_{k=1}^K\sum\limits_{i=1}^{\alpha_k^{\mathrm{G}}} \bm{N}_{k,i},
\end{eqnarray}
where $\left\vert\lambda_k\right\vert<R$, and matrices $\bm{P}_{k,i}$ and $\bm{N}_{k,i}$ are projector and nilpotent matrices with respect to the $k$-th eigenvalue and its $i$-th geometry component. 

Then, we have
\begin{eqnarray}\label{eq2: thm: Spectral Mapping Theorem for Single Variable}
f(\bm{X})&=&\sum\limits_{k=1}^K \left[\sum\limits_{i=1}^{\alpha_k^{(\mathrm{G})}}f(\lambda_k)\bm{P}_{k,i}+\sum\limits_{i=1}^{\alpha_k^{(\mathrm{G})}}\sum\limits_{q=1}^{m_{k,i}-1}\frac{f^{(q)}(\lambda_k)}{q!}\bm{N}_{k,i}^q\right].
\end{eqnarray}
\end{theorem}

Given a matrix $\bm{X}$ with the following Jordan decompostion form:
\begin{eqnarray}\label{eq1: matrix setup GDO for Nilpotent}
\bm{X}&=& \bm{U}\left(\bigoplus\limits_{k=1}^{K}\bigoplus\limits_{i=1}^{\alpha_{k}^{(\mathrm{G})}}\bm{J}_{m_{k,i}}(\lambda_{k})\right)\bm{U}^{-1},
\end{eqnarray}
we have the following expression of $f(\bm{X})$ from Theorem~\ref{thm: Spectral Mapping Theorem for Single Variable}:
\begin{eqnarray}\label{eq2: matrix setup GDO for Nilpotent}
f(\bm{X})&=& \bm{U}\left(\bigoplus\limits_{k=1}^{K}\bigoplus\limits_{i=1}^{\alpha_{k}^{(\mathrm{G})}}f(\bm{J}_{m_{k,i}}(\lambda_{k}))\right)\bm{U}^{-1},
\end{eqnarray}

We have the following Theorem~\ref{thm: Nilpotent Part of f(X) and GDOD} about the matrix representations $\mathfrak{R}(f(\mathbf{X}))$ and GDOD of nilpotent part of $\mathfrak{R}(f(\mathbf{X}))$. 

\begin{theorem}\label{thm: Nilpotent Part of f(X) and GDOD}
Given a matrix $\bm{X}\in \mathbb{C}^{m \times m}$ with the Jordan decompostion form as Eq.~\eqref{eq1: matrix setup GDO for Nilpotent}, we have 
\begin{eqnarray}\label{eq1: thm: Nilpotent Part of f(X) and GDOD}
\mathfrak{R}(\bm{X})&=&[ [\lambda_1]_{\alpha_1^{\mathrm{A}}},\ldots, [\lambda_K]_{\alpha_K^{\mathrm{A}}}, [m_{1,i}]_{1 \leq i \leq \alpha_1^{\mathrm{G}}},\ldots,[m_{k,i}]_{1 \leq i \leq \alpha_K^{\mathrm{G}}}],
\end{eqnarray}
where $[\lambda_k]_{\alpha_k^{\mathrm{A}}}$ is a vector with entries as $\lambda_k$ and size as $\alpha_k^{\mathrm{A}} = \sum\limits_{i=1}^{\alpha_k^{\mathrm{G}}} m_{k,i}$ for $k=1,2,\ldots,K$. We also have $\sum\limits_{k=1}^{K}\alpha_k^{\mathrm{A}}= m$.
Besides the conditions required by Theorem~\ref{thm: Spectral Mapping Theorem for Single Variable} for the differentiable function $f$, we also assume that the value $\kappa_k$ is the smallest integer such that $f^{(\kappa_k)}(\lambda_k) \neq 0$, where  $k=1,2,\ldots,K$. After applying $f$ to eigenvalues $\lambda_1, \ldots, \lambda_K$, we have the following ordering 
\begin{eqnarray}\label{eq1-1: thm: Nilpotent Part of f(X) and GDOD}
f(\lambda_{\sigma(1)}) \geq f(\lambda_{\sigma(2)}) \geq \ldots \geq f(\lambda_{\sigma(K)}), 
\end{eqnarray}
where $\sigma$ is the permutation operator of indices $1,2,\ldots,K$.

Then, we have 
\begin{eqnarray}\label{eq2: thm: Nilpotent Part of f(X) and GDOD}
\mathfrak{R}(f(\bm{X}))&=&[ [f(\lambda_{\sigma(1)})]_{\alpha_{\sigma(1)}^{\mathrm{A}}},\ldots, [f(\lambda_{\sigma(K)})]_{\alpha_{\sigma(K)}^{\mathrm{A}}}, \bm{\eta}_{f,\lambda_{\sigma(1)}},\ldots,\bm{\eta}_{f,\lambda_{\sigma(K)}}],
\end{eqnarray}
where vectors $\bm{\eta}_{f,\lambda_{\sigma(k)}}$ for $k=1,2,\ldots,K$ are 
\begin{eqnarray}\label{eq3: thm: Nilpotent Part of f(X) and GDOD}
\lefteqn{\bm{\eta}_{f,\lambda_{\sigma(k)}}=}\nonumber \\
&&\myop_{i=1}^{\alpha_{\sigma(k)}^{G}}[ (m_{\sigma(k),i}-\mathfrak{D}_{\bm{m}_{\sigma(k)}(f(\bm{J}_{m_{{\sigma(k)},i}}(\lambda_{\sigma(k)}))),\bm{m}_{\sigma(k)}(\bm{J}_{m_{{\sigma(k)},i}}(\lambda_{\sigma(k)}))}(1)),  \nonumber \\
&& (\mathfrak{D}_{\bm{m}_{\sigma(k)}(f(\bm{J}_{m_{{\sigma(k)},i}}(\lambda_{\sigma(k)}))),\bm{m}_{\sigma(k)}(\bm{J}_{m_{{\sigma(k)},i}}(\lambda_{\sigma(k)}))}(1) - \nonumber \\
&& \mathfrak{D}_{\bm{m}_{\sigma(k)}(f(\bm{J}_{m_{{\sigma(k)},i}}(\lambda_{\sigma(k)}))),\bm{m}_{\sigma(k)}(\bm{J}_{m_{{\sigma(k)},i}}(\lambda_{\sigma(k)}))}(2)),\cdots,\nonumber \\
&& (\mathfrak{D}_{\bm{m}_{\sigma(k)}(f(\bm{J}_{m_{{\sigma(k)},i}}(\lambda_{\sigma(k)}))),\bm{m}_{\sigma(k)}(\bm{J}_{m_{{\sigma(k)},i}}(\lambda_{\sigma(k)}))}(\kappa_{\sigma(k)}-1) - \nonumber \\
&& \mathfrak{D}_{\bm{m}_{\sigma(k)}(f(\bm{J}_{m_{{\sigma(k)},i}}(\lambda_{\sigma(k)}))),\bm{m}_{\sigma(k)}(\bm{J}_{m_{{\sigma(k)},i}}(\lambda_{\sigma(k)}))}(\kappa_{\sigma(k)})).
]
\end{eqnarray}

Moreover, we have an array of generalized dominance distance vectors with respect to each eigenvalue $\lambda_k$, denoted by $\underline{\mathfrak{D}}_{f,\bm{X}}$, and we can express it as:
\begin{eqnarray}\label{eq4: thm: Nilpotent Part of f(X) and GDOD}
\underline{\mathfrak{D}}_{f,\bm{X}}
&=&
\left[\mathfrak{D}_{\bm{m}_{\sigma(k)}\left(f\left(\bigoplus\limits_{i=1}^{\alpha_{\sigma(k)}^{(\mathrm{G})}}\bm{J}_{m_{{\sigma(k)},i}}(\lambda_{\sigma(k)})\right)\right), \bm{m}_{\sigma(k)}\left(\bigoplus\limits_{i=1}^{\alpha_{\sigma(k)}^{(\mathrm{G})}}\bm{J}_{m_{{\sigma(k)},i}}(\lambda_{\sigma(k)})\right)}\right]
\end{eqnarray}
where $k=1,2,\ldots,K$, and the term $\mathfrak{D}_{\bm{m}_{\sigma(k)}\left(f\left(\bigoplus\limits_{i=1}^{\alpha_{\sigma(k)}^{(\mathrm{G})}}\bm{J}_{m_{{\sigma(k)},i}}(\lambda_{\sigma(k)})\right)\right), \bm{m}_{\sigma(k)}\left(\bigoplus\limits_{i=1}^{\alpha_{\sigma(k)}^{(\mathrm{G})}}\bm{J}_{m_{{\sigma(k)},i}}(\lambda_{\sigma(k)})\right)}$ comes from Lemma~\ref{lma: GDO Distance multiple Jordan block same e-value}.
\end{theorem}
\textbf{Proof:}
From Eq.~\eqref{eq: f Jordan Block} and Eq.~\eqref{eq2: matrix setup GDO for Nilpotent}, the spectral part of $\mathfrak{R}(f(\bm{X}))$ will be
\begin{eqnarray}\label{eq5: thm: Nilpotent Part of f(X) and GDOD}
[f(\lambda_{\sigma(1)})]_{\alpha_{\sigma(1)}^{\mathrm{A}}},\ldots, [f(\lambda_{\sigma(K)})]_{\alpha_{\sigma(K)}^{\mathrm{A}}}.
\end{eqnarray}

For the nilpotent part, from Eq.~\eqref{eq2: matrix setup GDO for Nilpotent} again, we know that the nilpotent part is composed by the nilpotent contributed by each eigenvalue. From Lemma~\ref{lma: GDO Distance multiple Jordan block same e-value} and Eq.~\eqref{eq2: lma: GDO Distance multiple Jordan block same e-value}, we have the vector $\bm{\eta}_{f,\lambda_{\sigma(k)}}$ given by Eq.~\eqref{eq3: thm: Nilpotent Part of f(X) and GDOD} under the condition that $\kappa_{\sigma(k)}$ is the smallest integer such that $f^{(\kappa_{\sigma(k)})}(\lambda_{\sigma(k)}) \neq 0$. Accoding to Definition~\ref{def: SNO}, the nilpotent part of $\mathfrak{R}(f(\bm{X}))$ will be
\begin{eqnarray}\label{eq5: thm: Nilpotent Part of f(X) and GDOD}
[\bm{\eta}_{f,\lambda_{\sigma(1)}},\ldots,\bm{\eta}_{f,\lambda_{\sigma(K)}}].
\end{eqnarray}

Finally, the array of generalized dominance distance vectors with respect to each eigenvalue $\lambda_k$, represented by $\underline{\mathfrak{D}}_{f,\bm{X}}$, can be obtained by Eq.~\eqref{eq1: lma: GDO Distance multiple Jordan block same e-value} in Lemma~\ref{lma: GDO Distance multiple Jordan block same e-value}. 
$\hfill\Box$

\subsection{Matrix Representations Under Special $f$}\label{sec: under special f}

In this section, we will present several corollaries derived from Theorem~\ref{thm: Nilpotent Part of f(X) and GDOD} by examining specific cases of the function \( f \).

Following Corollary~\ref{cor: Nilpotent Part of f(X) and GDOD} is used to determine $\mathfrak{R}(f(\bm{X}))$ and the generalized dominance distance between the nilpotent part of $\mathfrak{R}(f(\bm{X}))$ and the nilpotent part of $\mathfrak{R}(\bm{X})$ under different  condtions of the increasing function $f$.  

\begin{corollary}\label{cor: Nilpotent Part of f(X) and GDOD}
Given a matrix $\bm{X}\in \mathbb{C}^{m \times m}$ with the Jordan decompostion form as Eq.~\eqref{eq1: matrix setup GDO for Nilpotent}, we have 
\begin{eqnarray}\label{eq1: cor: Nilpotent Part of f(X) and GDOD}
\mathfrak{R}(\bm{X})&=&[ [\lambda_1]_{\alpha_1^{\mathrm{A}}},\ldots, [\lambda_K]_{\alpha_K^{\mathrm{A}}}, [m_{1,i}]_{1 \leq i \leq \alpha_1^{\mathrm{G}}},\ldots,[m_{k,i}]_{1 \leq i \leq \alpha_K^{\mathrm{G}}}],
\end{eqnarray}
where $[\lambda_k]_{\alpha_k^{\mathrm{A}}}$ is a vector with entries as $\lambda_k$ and size as $\alpha_k^{\mathrm{A}} = \sum\limits_{i=1}^{\alpha_k^{\mathrm{G}}} m_{k,i}$ for $k=1,2,\ldots,K$. We also have $\sum\limits_{k=1}^{K}\alpha_k^{\mathrm{A}}= m$.
Besides the conditions required by Theorem~\ref{thm: Spectral Mapping Theorem for Single Variable} for the increasing function $f$, we also assume that the value $\kappa_k$ is the smallest integer such that $f^{(\kappa_k)}(\lambda_k) \neq 0$. Then, we have 
\begin{eqnarray}\label{eq2: cor: Nilpotent Part of f(X) and GDOD}
\mathfrak{R}(f(\bm{X}))&=&[ [f(\lambda_1)]_{\alpha_1^{\mathrm{A}}},\ldots, [f(\lambda_K)]_{\alpha_K^{\mathrm{A}}}, \bm{\eta}_{f,\lambda_1},\ldots,\bm{\eta}_{f,\lambda_K}],
\end{eqnarray}
where vectors $\bm{\eta}_{f,\lambda_k}$ for $k=1,2,\ldots,K$ are 
\begin{eqnarray}\label{eq3: cor: Nilpotent Part of f(X) and GDOD}
\bm{\eta}_{f,\lambda_k}&=&\myop_{i=1}^{\alpha_k^{G}}[ (m_{k,i}-\mathfrak{D}_{\bm{m}_k(f(\bm{J}_{m_{k,i}}(\lambda_k))),\bm{m}_k(\bm{J}_{m_{k,i}}(\lambda_k))}(1)),  \nonumber \\
&& (\mathfrak{D}_{\bm{m}_k(f(\bm{J}_{m_{k,i}}(\lambda_k))),\bm{m}_k(\bm{J}_{m_{k,i}}(\lambda_k))}(1) - \mathfrak{D}_{\bm{m}_k(f(\bm{J}_{m_{k,i}}(\lambda_k))),\bm{m}_k(\bm{J}_{m_{k,i}}(\lambda_k))}(2)),\cdots,\nonumber \\
&& (\mathfrak{D}_{\bm{m}_k(f(\bm{J}_{m_{k,i}}(\lambda_k))),\bm{m}_k(\bm{J}_{m_{k,i}}(\lambda_k))}(\kappa_k-1) - \mathfrak{D}_{\bm{m}_k(f(\bm{J}_{m_{k,i}}(\lambda_k))),\bm{m}_k(\bm{J}_{m_{k,i}}(\lambda_k))}(\kappa_k)).
]
\end{eqnarray}

Moreover, we have an array of generalized dominance distance vectors with respect to each eigenvalue $\lambda_k$, denoted by $\underline{\mathfrak{D}}_{f,\bm{X}}$, and we can express it as:
\begin{eqnarray}\label{eq4: cor: Nilpotent Part of f(X) and GDOD}
\underline{\mathfrak{D}}_{f,\bm{X}}
&=&
\left[\mathfrak{D}_{\bm{m}_k\left(f\left(\bigoplus\limits_{i=1}^{\alpha_k^{(\mathrm{G})}}\bm{J}_{m_{k,i}}(\lambda_k)\right)\right), \bm{m}_k\left(\bigoplus\limits_{i=1}^{\alpha_k^{(\mathrm{G})}}\bm{J}_{m_{k,i}}(\lambda_k)\right)}\right]
\end{eqnarray}
where $k=1,2,\ldots,K$ and the term $\mathfrak{D}_{\bm{m}_k\left(f\left(\bigoplus\limits_{i=1}^{\alpha_k^{(\mathrm{G})}}\bm{J}_{m_{k,i}}(\lambda_k)\right)\right), \bm{m}_k\left(\bigoplus\limits_{i=1}^{\alpha_k^{(\mathrm{G})}}\bm{J}_{m_{k,i}}(\lambda_k)\right)}$ comes from Lemma~\ref{lma: GDO Distance multiple Jordan block same e-value}.
\end{corollary}
\textbf{Proof:}
From Eq.~\eqref{eq: f Jordan Block} and Eq.~\eqref{eq2: matrix setup GDO for Nilpotent}, the spectral part of $\mathfrak{R}(f(\bm{X}))$ will be
\begin{eqnarray}\label{eq5: cor: Nilpotent Part of f(X) and GDOD}
[f(\lambda_1)]_{\alpha_1^{\mathrm{A}}},\ldots, [f(\lambda_K)]_{\alpha_K^{\mathrm{A}}}.
\end{eqnarray}

For the nilpotent part, from Eq.~\eqref{eq2: matrix setup GDO for Nilpotent} again, we know that the nilpotent part is composed by the nilpotent contributed by each eigenvalue. From Lemma~\ref{lma: GDO Distance multiple Jordan block same e-value} and Eq.~\eqref{eq2: lma: GDO Distance multiple Jordan block same e-value}, we have the vector $\bm{\eta}_{f,\lambda_k}$ given by Eq.~\eqref{eq3: thm: Nilpotent Part of f(X) and GDOD} under the condition that $\kappa_k$ is the smallest integer such that $f^{(\kappa_k)}(\lambda_k) \neq 0$. Accoding to Definition~\ref{def: SNO}, the nilpotent part of $\mathfrak{R}(f(\bm{X}))$ will be
\begin{eqnarray}\label{eq5: cor: Nilpotent Part of f(X) and GDOD}
[\bm{\eta}_{f,\lambda_1},\ldots,\bm{\eta}_{f,\lambda_K}].
\end{eqnarray}

Finally, the array of generalized dominance distance vectors with respect to each eigenvalue $\lambda_k$, represented by $\underline{\mathfrak{D}}_{f,\bm{X}}$, can be obtained by Eq.~\eqref{eq1: lma: GDO Distance multiple Jordan block same e-value} in Lemma~\ref{lma: GDO Distance multiple Jordan block same e-value}. 
$\hfill\Box$

Below, we present Corollary~\ref{cor: kappa 1} and Corollary~\ref{cor: kappa larger or equal m}, highlighting variations in the representation \(\mathfrak{R}(f(\bm{X}))\) based on the properties of the special function \(f\).

\begin{corollary}\label{cor: kappa 1}
Given a matrix $\bm{X}\in \mathbb{C}^{m \times m}$ with the representation $\mathfrak{R}(\bm{X})$ as
\begin{eqnarray}\label{eq1: cor: kappa 1}
\mathfrak{R}(\bm{X})&=&[ [\lambda_1]_{\alpha_1^{\mathrm{A}}},\ldots, [\lambda_K]_{\alpha_K^{\mathrm{A}}}, [m_{1,i}]_{1 \leq i \leq \alpha_1^{\mathrm{G}}},\ldots,[m_{K,i}]_{1 \leq i \leq \alpha_K^{\mathrm{G}}}],
\end{eqnarray}
where $[\lambda_k]_{\alpha_k^{\mathrm{A}}}$ is a vector with entries as $\lambda_k$ and size as $\alpha_k^{\mathrm{A}} = \sum\limits_{i=1}^{\alpha_k^{\mathrm{G}}} m_{k,i}$ for $k=1,2,\ldots,K$. We also have $\sum\limits_{k=1}^{K}\alpha_k^{\mathrm{A}}= m$.
Besides the conditions required by Theorem~\ref{thm: Spectral Mapping Theorem for Single Variable} for the increasing function $f$, we also assume that the value $\kappa_k=1$ with respect to $f(\lambda_k)$, i.e., $f^{(1)}(\lambda_k) \neq 0$ for $k=1,2,\ldots,K$. Then, we have 
\begin{eqnarray}\label{eq2: cor: kappa 1}
\mathfrak{R}(f(\bm{X}))&=&[[f(\lambda_1)]_{\alpha_1^{\mathrm{A}}},\ldots, [f(\lambda_K)]_{\alpha_K^{\mathrm{A}}}, [m_{1,i}]_{1 \leq i \leq \alpha_1^{\mathrm{G}}},\ldots,[m_{K,i}]_{1 \leq i \leq \alpha_K^{\mathrm{G}}}].
\end{eqnarray}

Moreover, we have an array of generalized dominance distance vectors with respect to each eigenvalue $\lambda_k$, denoted by $\underline{\mathfrak{D}}_{f,\bm{X}}$, and we can express it as:
\begin{eqnarray}\label{eq3: cor: kappa 1}
\underline{\mathfrak{D}}_{f,\bm{X}}
&=&
\left[\overbrace{0,0,\ldots,0}^{\mbox{$K$ zeros}}\right].
\end{eqnarray}
\end{corollary}
\textbf{Proof:}
The function of the Jordan block $\bm{J}_{m_{k,i}}(\lambda_k)$ is
\begin{eqnarray}\label{eq4: cor: kappa 1}
f(\bm{J}_{m_{k,i}}(\lambda_k))&=&\begin{pmatrix}
   f(\lambda_k) & \frac{f'(\lambda_k)}{1!} & \frac{f''(\lambda_k)}{2!} & \cdots & \frac{f^{(m_{k,i}-1)}(\lambda_k)}{(m_{k,i}-1)!} \\
   0 & f(\lambda_k) & \frac{f'(\lambda_k)}{1!} & \cdots & \frac{f^{(m_{k,i}-2)}(\lambda_k)}{(m_{k,i}-2)!}  \\
   \vdots & \vdots & \ddots & \ddots & \vdots \\
   0 & 0 & \cdots & f(\lambda) & \frac{f'(\lambda_k)}{1!} \\
   0 & 0 & \cdots & 0 & f(\lambda_k)
   \end{pmatrix}_{m_{k,i} \times m_{k,i}},
\end{eqnarray}
where $i=1,2,\ldots,\alpha_k^{\mathrm{G}}$. Because $f^{(1)}(\lambda_k) \neq 0$ for $k=1,2,\ldots,K$, the rank of $f(\bm{J}_{m_{k,i}}(\lambda_k)) - f(\lambda_k)\bm{I}$ becomes $m_{k,i}-1$. Then, the Jordan block for the eigenvalue $f(\lambda_k)$ will be
\begin{eqnarray}\label{eq5: cor: kappa 2}
\begin{pmatrix}
   f(\lambda_k) & 1 & 0 & \cdots & 0 \\
   0 & f(\lambda_k) & 1 & \cdots & 0  \\
   \vdots & \vdots & \ddots & \ddots & \vdots \\
   0 & 0 & \cdots & f(\lambda) & 1 \\
   0 & 0 & \cdots & 0 & f(\lambda_k)
   \end{pmatrix}_{m_{k,i} \times m_{k,i}},
\end{eqnarray}
Therefore, the nilpotent part of $\mathfrak{R}(f(\bm{X}))$ will still be $[[m_{1,i}]_{1 \leq i \leq \alpha_1^{\mathrm{G}}},\ldots,[m_{K,i}]_{1 \leq i \leq \alpha_K^{\mathrm{G}}}]$. 

Finally, by applying Lemma~\ref{lma: GDO Distance single Jordan block} and Lemma~\ref{lma: GDO Distance multiple Jordan block same e-value}, we have desired result provided by Eq.~\eqref{eq3: cor: kappa 1}.
$\hfill\Box$

\begin{corollary}\label{cor: kappa larger or equal m}
%
Given a matrix $\bm{X}\in \mathbb{C}^{m \times m}$  with the representation $\mathfrak{R}(\bm{X})$ as
\begin{eqnarray}\label{eq1: cor: kappa larger or equal m}
\mathfrak{R}(\bm{X})&=&[ [\lambda_1]_{\alpha_1^{\mathrm{A}}},\ldots, [\lambda_K]_{\alpha_K^{\mathrm{A}}}, [m_{1,i}]_{1 \leq i \leq \alpha_1^{\mathrm{G}}},\ldots,[m_{k,i}]_{1 \leq i \leq \alpha_K^{\mathrm{G}}}],
\end{eqnarray}
where $[\lambda_k]_{\alpha_k^{\mathrm{A}}}$ is a vector with entries as $\lambda_k$ and size as $\alpha_k^{\mathrm{A}} = \sum\limits_{i=1}^{\alpha_k^{\mathrm{G}}} m_{k,i}$ for $k=1,2,\ldots,K$. We also have $\sum\limits_{k=1}^{K}\alpha_k^{\mathrm{A}}= m$.
Besides the conditions required by Theorem~\ref{thm: Spectral Mapping Theorem for Single Variable} for the increasing function $f$, we also assume that the value $\kappa_k \geq  \max\limits_{\substack{k=1,2,\ldots,K \\ i=1,2,\ldots,\alpha_k^{\mathrm{G}}}}
m_{k,i}$ with respect to $f(\lambda_k)$, i.e., $f^{(\kappa_k)}(\lambda_k) \neq 0$ for $k=1,2,\ldots,K$. 

Then, we have 
\begin{eqnarray}\label{eq2: cor: kappa larger or equal m}
\mathfrak{R}(f(\bm{X}))&=&[ [f(\lambda_1)]_{\alpha_1^{\mathrm{A}}},\ldots, [f(\lambda_K)]_{\alpha_K^{\mathrm{A}}}, [\overbrace{1,\ldots,1}^{\mbox{$\alpha_1^{\mathrm{A}}$ ones}}],\ldots,
[\overbrace{1,\ldots,1}^{\mbox{$\alpha_K^{\mathrm{A}}$ ones}}]],
\end{eqnarray}

Moreover, we have an array of generalized dominance distance vectors with respect to each eigenvalue $\lambda_k$, denoted by $\underline{\mathfrak{D}}_{f,\bm{X}}$, and we can express it as:
\begin{eqnarray}\label{eq3: cor: kappa larger or equal m}
\underline{\mathfrak{D}}_{f,\bm{X}}
&=&
\left[\mathfrak{D}_{\bm{m}_k\left(f\left(\bigoplus\limits_{i=1}^{\alpha_k^{(\mathrm{G})}}\bm{J}_{m_{k,i}}(\lambda_k)\right)\right), \bm{m}_k\left(\bigoplus\limits_{i=1}^{\alpha_k^{(\mathrm{G})}}\bm{J}_{m_{k,i}}(\lambda_k)\right)}(j)=\begin{cases} 
\sum\limits_{i=1}^j m_{k,i} - j, & \mbox{if $1 \leq j \leq \alpha_k^{\mathrm{G}}$}; \\
\alpha_k^{\mathrm{A}} - j, & \mbox{if $\alpha_k^{\mathrm{G}} < j \leq \alpha_k^{\mathrm{A}}$}; \\
\end{cases}\right],\nonumber \\
\end{eqnarray}
where $k=1,2,\ldots,K$.
\end{corollary}
\textbf{Proof:}
The function of the Jordan block $\bm{J}_{m_{k,i}}(\lambda_k)$ is
\begin{eqnarray}\label{eq4: cor: kappa larger or equal m}
f(\bm{J}_{m_{k,i}}(\lambda_k))&=&\begin{pmatrix}
   f(\lambda_k) & \frac{f'(\lambda_k)}{1!} & \frac{f''(\lambda_k)}{2!} & \cdots & \frac{f^{(m_{k,i}-1)}(\lambda_k)}{(m_{k,i}-1)!} \\
   0 & f(\lambda_k) & \frac{f'(\lambda_k)}{1!} & \cdots & \frac{f^{(m_{k,i}-2)}(\lambda_k)}{(m_{k,i}-2)!}  \\
   \vdots & \vdots & \ddots & \ddots & \vdots \\
   0 & 0 & \cdots & f(\lambda) & \frac{f'(\lambda_k)}{1!} \\
   0 & 0 & \cdots & 0 & f(\lambda_k)
   \end{pmatrix}_{m_{k,i} \times m_{k,i}},
\end{eqnarray}
where $i=1,2,\ldots,\alpha_k^{\mathrm{G}}$. Because $\kappa_k \geq  \max\limits_{\substack{k=1,2,\ldots,K \\ i=1,2,\ldots,\alpha_k^{\mathrm{G}}}}m_{k,i}$, this implies that $f'(\lambda_k) = f''(\lambda_k) = \ldots = f^{(m_{k,i}-1)}(\lambda_k)=0$.   The matrix $f(\bm{J}_{m_{k,i}}(\lambda_k))$ becomes 
\begin{eqnarray}\label{eq5: cor: kappa larger or equal m}
f(\bm{J}_{m_{k,i}}(\lambda_k))&=&\begin{pmatrix}
   f(\lambda_k) &0 & 0 & \cdots & 0 \\
   0 & f(\lambda_k) & 0& \cdots & 0  \\
   \vdots & \vdots & \ddots & \ddots & \vdots \\
   0 & 0 & \cdots & f(\lambda) & 0 \\
   0 & 0 & \cdots & 0 & f(\lambda_k)
   \end{pmatrix}_{m_{k,i} \times m_{k,i}},
\end{eqnarray}
then, $\bm{m}_k(f(\bm{J}_{m_{k,i}}(\lambda_k)))=[\overbrace{1,1,\ldots,1}^{\mbox{$m_{k,i}$ ones}}]$. Moreover, we also have 
\begin{eqnarray}\label{eq6: cor: kappa larger or equal m}
\raisebox{13em}{\(f\left(\bigoplus\limits_{i=1}^{\alpha_{k}^{G}}\bm{J}_{m_{k,i}}(\lambda_k)\right)\)=}&&
\begin{tikzpicture}
    \matrix[matrix of math nodes,left delimiter={(},right delimiter={)}] (m) {
      f(\lambda_k) & 0       & 0   & \ldots   & 0      & 0     & 0      & 0  & 0      & 0  & 0 & 0\\
        0       &  f(\lambda_k)& 0    & 0 \ldots    & 0     & 0      & 0  & 0      & 0 & 0  & 0 & 0 \\
        0      &  0 & \ddots &   \ddots   & 0      & 0     & 0      & 0  & 0      & 0  & 0 & 0  \\
        0      & \ldots   & 0     & f(\lambda_k)   &  0   & 0     & 0      & 0  & 0      & 0   & 0  & 0 \\
        0      & 0      & \ldots      &  0 &  f(\lambda_k)  & 0      & 0  & 0      & 0 & 0  & 0  & 0 \\
        0      & 0      & 0      & 0      &  \ddots   & \ddots     & \ddots    & 0      & 0 & 0  & 0  & 0\\
        0      & 0    & 0   & 0      & 0      &  \ddots   & \ddots       & \ddots   & 0      & 0 & 0  & 0  \\
        0      & 0    & 0   & 0      & 0      & 0     & 0      & f(\lambda_k)  & 0      & 0  & \ldots  & 0 \\
        0      & 0    & 0  & 0      & 0      & 0     & 0     &0     &  f(\lambda_k) & 0    & 0 \ldots    & 0  \\
        0      & 0    & 0  & 0      & 0      & 0     & 0      &  0      &  0 & \ddots &   \ddots   & 0 \\
        0      & 0    & 0  & 0      & 0      & 0     & 0      &  0      & \ldots   & 0     & f(\lambda_k)   &  0 \\
        0      & 0    & 0  & 0      & 0      & 0     & 0      &  0      & 0      & \ldots      &  0 &  f(\lambda_k) \\
    };

    \draw[dashed] 
        (m-1-1.north west) -- (m-1-5.north east) -- 
        (m-5-5.south east) -- (m-5-1.south west) -- cycle;
\node[anchor=north] at ($(m-5-5.south east) + (0.0,0)$) {\tiny $m_{k,1} \times m_{k,1}$};

    \draw[dashed] 
        (m-8-8.north west) -- (m-8-12.north east) -- 
        (m-12-12.south east) -- (m-12-8.south west) -- cycle;
\node[anchor=north] at ($(m-12-12.south east) + (0.0,0)$){\tiny $m_{k,\alpha_k^{G}} \times m_{k,\alpha_k^{G}}$};

\end{tikzpicture} \nonumber 
\end{eqnarray}
which indicates that 
\begin{eqnarray}\label{eq7: cor: kappa larger or equal m}
\bm{m}_k\left(f\left(\bigoplus\limits_{i=1}^{\alpha_k^{(\mathrm{G})}}\bm{J}_{m_{k,i}}(\lambda_k)\right)\right)=
[\overbrace{1,\ldots,1}^{\mbox{$\alpha_k^{\mathrm{A}}$ ones}}], 
\end{eqnarray}
due to the fact that $\alpha_k^{\mathrm{A}} = \sum\limits_{i=1}^{\mathrm{G}}m_{k,i}$.

Finally, by applying the Generalized Dominance Ordering Distance definition provided by Eq.~\eqref{eq: dom order distance}, we have desired result provided by Eq.~\eqref{eq3: cor: kappa larger or equal m}.
$\hfill\Box$

The following Corollary~\ref{cor: GDOD between f X and g Y} will discuss the generalized dominance distance relationship between $f(\bm{X})$ and $g(\bm{Y})$.

\begin{corollary}\label{cor: GDOD between f X and g Y}
Given a matrix $\bm{X}\in \mathbb{C}^{m \times m}$ with the representation $\mathfrak{R}(\bm{X})$ as
\begin{eqnarray}\label{eq1: cor: GDOD between f X and g Y}
\mathfrak{R}(\bm{X})&=&[ [\lambda_{x,1}]_{\alpha_{x,1}^{\mathrm{A}}},\ldots, [\lambda_{x,K_{x}}]_{\alpha_{x,K_{x}}^{\mathrm{A}}}, [m_{x,1,i}]_{1 \leq i \leq \alpha_{x,1}^{\mathrm{G}}},\ldots,[m_{x,K_{x},i}]_{1 \leq i \leq \alpha_{x,K_{x}}^{\mathrm{G}}}],
\end{eqnarray}
where $[\lambda_{x,k}]_{\alpha_{x,k}^{\mathrm{A}}}$ is a vector with entries as eigenvalues $\lambda_{x,k}$ of $\bm{X}$ and size as $\alpha_{x,k}^{\mathrm{A}} = \sum\limits_{i=1}^{\alpha_{x,k}^{\mathrm{G}}} m_{x,k,i}$ for $k=1,2,\ldots,K_{x}$. Note that $\bm{m}_{x,k}=[m_{x,k,i}]_{1 \leq i \leq \alpha_{x,k}^{\mathrm{G}}}$, where $k=1,2,\ldots,K_x$. We also have $\sum\limits_{k=1}^{K_{x}}\alpha_{x,k}^{\mathrm{A}}= m$. The function $f$ satisfies the conditions required by Theorem~\ref{thm: Spectral Mapping Theorem for Single Variable}, and we also assume that the value $\kappa_{x,k}$ is the smallest integer such that $f^{(\kappa_{x,k})}(\lambda_{x,k}) \neq 0$. 

On the other hand, we also have another matrix $\bm{Y}\in \mathbb{C}^{m \times m}$ with the representation $\mathfrak{R}(\bm{Y})$ as
\begin{eqnarray}\label{eq2: cor: GDOD between f X and g Y}
\mathfrak{R}(\bm{Y})&=&[ [\lambda_{y,1}]_{\alpha_{y,1}^{\mathrm{A}}},\ldots, [\lambda_{y,K_{y}}]_{\alpha_{y,K_{y}}^{\mathrm{A}}}, [m_{y,1,i}]_{1 \leq i \leq \alpha_{y,1}^{\mathrm{G}}},\ldots,[m_{y,K_{y},i}]_{1 \leq i \leq \alpha_{y,K_{y}}^{\mathrm{G}}}],
\end{eqnarray}
where $[\lambda_{y,k}]_{\alpha_{y,k}^{\mathrm{A}}}$ is a vector with entries as eigenvalues $\lambda_{y,k}$ of $\bm{X}$ and size as $\alpha_{y,k}^{\mathrm{A}} = \sum\limits_{i=1}^{\alpha_{y,k}^{\mathrm{G}}} m_{y,k,i}$ for $k=1,2,\ldots,K_{y}$. Note that $\bm{m}_{y,k}=[m_{y,k,i}]_{1 \leq i \leq \alpha_{y,k}^{\mathrm{G}}}$, where $k=1,2,\ldots,K_y$. We also have $\sum\limits_{k=1}^{K_{y}}\alpha_{y,k}^{\mathrm{A}}= m$. The function $g$ satisfies the conditions required by Theorem~\ref{thm: Spectral Mapping Theorem for Single Variable}, and we also assume that the value $\kappa_{y,k}$ is the smallest integer such that $g^{(\kappa_{y,k})}(\lambda_{y,k}) \neq 0$. 

For the eigenvalues \(\lambda_{x,k}\) and \(\lambda_{y,k}\), after applying the functions \(f\) and \(g\), respectively, we denote by \(\sigma_x\) and \(\sigma_y\) the permutations that reorder these eigenvalues according to their complex-valued total ordering. Thus, we have
\begin{eqnarray}\label{eq2-1: cor: GDOD between f X and g Y}
\lambda_{x,\sigma_x(1)} \geq \lambda_{x,\sigma_x(2)} \geq \ldots \geq \lambda_{x,\sigma_x(K_x)},\nonumber \\
\lambda_{y,\sigma_y(1)} \geq \lambda_{y,\sigma_y(2)} \geq \ldots \geq \lambda_{y,\sigma_y(K_y)}.
\end{eqnarray}
Then, we have 
\begin{eqnarray}\label{eq3: cor: GDOD between f X and g Y}
\mathfrak{R}(f(\bm{X}))&=&[ [f(\lambda_{x,\sigma_x(1)})]_{\alpha_{x,\sigma_x(1)}^{\mathrm{A}}},\ldots, [f(\lambda_{x,\sigma_x(K_x)})]_{\alpha_{x,\sigma_x(K_x)}^{\mathrm{A}}}, \bm{\eta}_{f,\lambda_{x,\sigma_x(1)}},\ldots,\bm{\eta}_{f,\lambda_{x,\sigma_x(K_x)}}],
\end{eqnarray}
where vectors $\bm{\eta}_{f,\lambda_{x,\sigma_x(k)}}$ for $k=1,2,\ldots,K_x$ are 
\begin{eqnarray}\label{eq4: cor: GDOD between f X and g Y}
\lefteqn{\bm{\eta}_{f,\lambda_{x,\sigma_x(k)}}=}\nonumber \\
&& \myop_{i=1}^{\alpha_{x,\sigma_x(k)}^{G}}[ (m_{x,\sigma_x(k),i}-\mathfrak{D}_{\bm{m}_{x,\sigma_x(k)}(f(\bm{J}_{m_{x,\sigma_x(k),i}}(\lambda_{x,\sigma_x(k)}))),\bm{m}_{x,\sigma_x(k)}(\bm{J}_{m_{x,\sigma_x(k),i}}(\lambda_{x,\sigma_x(k)}))}(1)),  \nonumber \\
&& (\mathfrak{D}_{\bm{m}_{x,\sigma_x(k)}(f(\bm{J}_{m_{x,\sigma_x(k),i}}(\lambda_{x,\sigma_x(k)}))),\bm{m}_{x,\sigma_x(k)}(\bm{J}_{m_{x,\sigma_x(k),i}}(\lambda_{x,\sigma_x(k)}))}(1) \nonumber \\
&& - \mathfrak{D}_{\bm{m}_{x,\sigma_x(k)}(f(\bm{J}_{m_{x,\sigma_x(k),i}}(\lambda_{x,\sigma(k)}))),\bm{m}_{x,\sigma_x(k)}(\bm{J}_{m_{x,\sigma_x(k),i}}(\lambda_{x,\sigma_x(k)}))}(2)),\cdots,\nonumber \\
&& (\mathfrak{D}_{\bm{m}_{x,\sigma_x(k)}(f(\bm{J}_{m_{x,\sigma_x(k),i}}(\lambda_{x,\sigma_x(k)}))),\bm{m}_{x,\sigma_x(k)}(\bm{J}_{m_{x,\sigma_x(k),i}}(\lambda_{x,\sigma_x(k)}))}(\kappa_{x,k}-1) \nonumber \\
&& - \mathfrak{D}_{\bm{m}_{x,\sigma_x(k)}(f(\bm{J}_{m_{x,\sigma_x(k),i}}(\lambda_{x,\sigma_x(k)}))),\bm{m}_{x,\sigma_x(k)}(\bm{J}_{m_{x,\sigma_x(k),i}}(\lambda_{x,\sigma_x(k)}))}(\kappa_{x,k})).
]
\end{eqnarray}
Similarly, we also have 
\begin{eqnarray}\label{eq5: cor: GDOD between f X and g Y}
\mathfrak{R}(g(\bm{Y}))&=&[ [g(\lambda_{y,\sigma_y(1)})]_{\alpha_{y,\sigma_y(1)}^{\mathrm{A}}},\ldots, [g(\lambda_{y,\sigma_y(K_y)})]_{\alpha_{y,\sigma_y(K_y)}^{\mathrm{A}}}, \bm{\eta}_{g,\lambda_{y,\sigma_y(1)}},\ldots,\bm{\eta}_{g,\lambda_{y,\sigma_y(K_y)}}],
\end{eqnarray}
where vectors $\bm{\eta}_{g,\lambda_{y,\sigma_y(k)}}$ for $k=1,2,\ldots,K_y$ are 
\begin{eqnarray}\label{eq6: cor: GDOD between f X and g Y}
\lefteqn{\bm{\eta}_{g,\lambda_{y,\sigma_y(k)}}=}\nonumber \\
&& \myop_{i=1}^{\alpha_{y,\sigma_y(k)}^{G}}[ (m_{y,\sigma_y(k),i}-\mathfrak{D}_{\bm{m}_{y,\sigma_y(k)}(g(\bm{J}_{m_{y,\sigma_y(k),i}}(\lambda_{y,\sigma_y(k)}))),\bm{m}_{y,\sigma_y(k)}(\bm{J}_{m_{y,\sigma_y(k),i}}(\lambda_{y,\sigma_y(k)}))}(1)),  \nonumber \\
&& (\mathfrak{D}_{\bm{m}_{y,\sigma_y(k)}(g(\bm{J}_{m_{y,\sigma_y(k),i}}(\lambda_{y,\sigma_y(k)}))),\bm{m}_{y,\sigma_y(k)}(\bm{J}_{m_{y,\sigma_y(k),i}}(\lambda_{y,\sigma_y(k)}))}(1) \nonumber \\
&& - \mathfrak{D}_{\bm{m}_{y,\sigma_y(k)}(g(\bm{J}_{m_{y,\sigma_y(k),i}}(\lambda_{y,\sigma(k)}))),\bm{m}_{y,\sigma_y(k)}(\bm{J}_{m_{y,\sigma_y(k),i}}(\lambda_{y,\sigma_y(k)}))}(2)),\cdots,\nonumber \\
&& (\mathfrak{D}_{\bm{m}_{y,\sigma_y(k)}(g(\bm{J}_{m_{y,\sigma_y(k),i}}(\lambda_{y,\sigma_y(k)}))),\bm{m}_{y,\sigma_y(k)}(\bm{J}_{m_{y,\sigma_y(k),i}}(\lambda_{y,\sigma_y(k)}))}(\kappa_{y,k}-1) \nonumber \\
&& - \mathfrak{D}_{\bm{m}_{y,\sigma_y(k)}(g(\bm{J}_{m_{y,\sigma_y(k),i}}(\lambda_{y,\sigma_y(k)}))),\bm{m}_{y,\sigma_y(k)}(\bm{J}_{m_{y,\sigma_y(k),i}}(\lambda_{y,\sigma_y(k)}))}(\kappa_{y,k})).
]
\end{eqnarray}

Moreover, we have an array of generalized dominance distance vectors with respect to nilpotent part of $f(\bm{X})$ and $g(\bm{Y})$, denoted by $\underline{\mathfrak{D}}_{f,\bm{X};g,\bm{Y}}$, and we can express it as:
\begin{eqnarray}\label{eq7: cor: GDOD between f X and g Y}
\lefteqn{\underline{\mathfrak{D}}_{f,\bm{X};g,\bm{Y}}=}\nonumber \\
&&
\Bigg[
\begin{cases} 
\mathfrak{D}_{\bm{\eta}_{g,\lambda_{y,\sigma_y(k)}},\bm{\eta}_{f,\lambda_{x,\sigma_x(k)}}}(j)=\sum\limits_{i=1}^j (\bm{\eta}_{f,\lambda_{x,\sigma_x(k)}}(i)-\bm{\eta}_{g,\lambda_{y,\sigma_y(k)}}(i)), & \mbox{if $\bm{\eta}_{g,\lambda_{y,\sigma_y(k)}}\trianglelefteq \bm{\eta}_{f,\lambda_{x,\sigma_x(k)}}$} \\
\mathfrak{D}_{\bm{\eta}_{f,\lambda_{x,\sigma_x(k)}},\bm{\eta}_{g,\lambda_{y,\sigma_y(k)}}}(j)=\sum\limits_{i=1}^j (\bm{\eta}_{g,\lambda_{y,\sigma_y(k)}}(i)-\bm{\eta}_{f,\lambda_{x,\sigma_x(k)}}(i)), & \mbox{if $\bm{\eta}_{f,\lambda_{x,\sigma_x(k)}}\trianglelefteq \bm{\eta}_{g,\lambda_{y,\sigma_y(k)}}$} \\
\end{cases}
\bigg],\nonumber \\
\end{eqnarray}
where $k=1,2,\ldots,\max(K_x,K_y)$. As usual, the shorter nilpotent part will be filled with zero vectors.
\end{corollary}
\textbf{Proof:}
Eq.~\eqref{eq3: cor: GDOD between f X and g Y} and Eq.~\eqref{eq5: cor: GDOD between f X and g Y} are obtained by applying Eq.~\eqref{eq2: thm: Nilpotent Part of f(X) and GDOD} and Eq.~\eqref{eq3: thm: Nilpotent Part of f(X) and GDOD} from Theorem~\ref{thm: Nilpotent Part of f(X) and GDOD} after reordering eigenvalues $\lambda_{x,k}$ and $\lambda_{y,k}$ according to Eq.~\eqref{eq2-1: cor: GDOD between f X and g Y}.

Finally, the array of generalized dominance distance vectors with respect to nilpotent part of $f(\bm{X})$ and $g(\bm{Y})$ is obtained by applying the definition of the generalized dominance ordering distance  provided by Eq.~\eqref{eq: dom order distance} to vectors given by Eq.~\eqref{eq4: cor: GDOD between f X and g Y} and Eq.~\eqref{eq6: cor: GDOD between f X and g Y}.
$\hfill\Box$

\section{Monotone and Convex Conditions}\label{sec: Monotone and Convex Conditions}

In this section, the conditions for a function with monotonicity under SNO is discussed in Section~\ref{sec: Monotonicity}, and the conditions for a function with convexity uner SNO is presented in Section~\ref{sec: Convexity}.

\subsection{Monotonicity}\label{sec: Monotonicity}

Given two matrices $\bm{X}$ and $\bm{Y}$ with $\bm{X} \preceq_{\mbox{\tiny SN}} \bm{Y}$, we will discuss the monotone function $f$ conditions to have $f(\bm{X}) \preceq_{\mbox{\tiny SN}} f(\bm{Y})$ based on Definition~\ref{def: SNO}.

\begin{theorem}\label{thm: monotonicity conditions}
Given a matrix $\bm{X}\in \mathbb{C}^{m \times m}$ with the representation $\mathfrak{R}(\bm{X})$ as
\begin{eqnarray}\label{eq1: thm: monotonicity conditions}
\mathfrak{R}(\bm{X})&=&[ [\lambda_{x,1}]_{\alpha_{x,1}^{\mathrm{A}}},\ldots, [\lambda_{x,K_{x}}]_{\alpha_{x,K_{x}}^{\mathrm{A}}}, [m_{x,1,i}]_{1 \leq i \leq \alpha_{x,1}^{\mathrm{G}}},\ldots,[m_{x,K_{x},i}]_{1 \leq i \leq \alpha_{x,K_{x}}^{\mathrm{G}}}],
\end{eqnarray}
where $[\lambda_{x,k}]_{\alpha_{x,k}^{\mathrm{A}}}$ is a vector with entries as eigenvalues $\lambda_{x,k}$ of $\bm{X}$ and size as $\alpha_{x,k}^{\mathrm{A}} = \sum\limits_{i=1}^{\alpha_{x,k}^{\mathrm{G}}} m_{x,k,i}$ for $k=1,2,\ldots,K_{x}$. On the other hand, we also have another matrix $\bm{Y}\in \mathbb{C}^{m \times m}$ with the representation $\mathfrak{R}(\bm{Y})$ as
\begin{eqnarray}\label{eq2: thm: monotonicity conditions}
\mathfrak{R}(\bm{Y})&=&[ [\lambda_{y,1}]_{\alpha_{y,1}^{\mathrm{A}}},\ldots, [\lambda_{y,K_{y}}]_{\alpha_{y,K_{y}}^{\mathrm{A}}}, [m_{y,1,i}]_{1 \leq i \leq \alpha_{y,1}^{\mathrm{G}}},\ldots,[m_{y,K_{y},i}]_{1 \leq i \leq \alpha_{y,K_{y}}^{\mathrm{G}}}],
\end{eqnarray}
where $[\lambda_{y,k}]_{\alpha_{y,k}^{\mathrm{A}}}$ is a vector with entries as eigenvalues $\lambda_{y,k}$ of $\bm{Y}$ and size as $\alpha_{y,k}^{\mathrm{A}} = \sum\limits_{i=1}^{\alpha_{y,k}^{\mathrm{G}}} m_{y,k,i}$ for $k=1,2,\ldots,K_{y}$. We also have $\sum\limits_{k=1}^{K_{y}}\alpha_{y,k}^{\mathrm{A}}= m$. 

For eigenvalues $\lambda_{x,k}$ and $\lambda_{y,k}$, after acting by the function $f$, we use $\sigma_x$ and $\sigma_y$ to indicate these eigenvalues reordering permutation of $\lambda_{x,k}$ and $\lambda_{y,k}$ by their complex-valued total ordering. We have
\begin{eqnarray}\label{eq2-1: thm: monotonicity conditions}
f(\lambda_{x,\sigma_x(1)})\geq f(\lambda_{x,\sigma_x(2)})\geq \ldots \geq f(\lambda_{x,\sigma_x(K_x)}),\nonumber \\
f(\lambda_{y,\sigma_y(1)})\geq f(\lambda_{y,\sigma_y(2)})\geq \ldots \geq f(\lambda_{y,\sigma_y(K_y)}).
\end{eqnarray}
Then, we have 
\begin{eqnarray}\label{eq3: thm: monotonicity conditions}
\mathfrak{R}(f(\bm{X}))&=&[ [f(\lambda_{x,\sigma_x(1)})]_{\alpha_{x,\sigma_x(1)}^{\mathrm{A}}},\ldots, [f(\lambda_{x,\sigma_x(K_x)})]_{\alpha_{x,\sigma_x(K_x)}^{\mathrm{A}}}, \bm{\eta}_{f,\lambda_{x,\sigma_x(1)}},\ldots,\bm{\eta}_{f,\lambda_{x,\sigma_x(K_x)}}],
\end{eqnarray}
and
\begin{eqnarray}\label{eq4: thm: monotonicity conditions}
\mathfrak{R}(f(\bm{Y}))&=&[ [f(\lambda_{y,\sigma_y(1)})]_{\alpha_{y,\sigma_y(1)}^{\mathrm{A}}},\ldots, [f(\lambda_{y,\sigma_y(K_y)})]_{\alpha_{y,\sigma_y(K_y)}^{\mathrm{A}}}, \bm{\eta}_{f,\lambda_{y,\sigma_y(1)}},\ldots,\bm{\eta}_{f,\lambda_{y,\sigma_y(K_y)}}],
\end{eqnarray}

There are two cases corresponding to $\bm{X} \prec_{\mbox{\tiny SN}} \bm{Y}$.

\textbf{Case I:  $[ [\lambda_{x,1}]_{\alpha_{x,1}^{\mathrm{A}}},\ldots, [\lambda_{x,K_{x}}]_{\alpha_{x,K_{x}}^{\mathrm{A}}}] \prec_w [ [\lambda_{y,1}]_{\alpha_{y,1}^{\mathrm{A}}},\ldots, [\lambda_{y,K_{y}}]_{\alpha_{y,K_{y}}^{\mathrm{A}}}]$}. For Case I, we have the following three possibilities to have $f(\bm{X}) \prec_{\mbox{\tiny SN}} f(\bm{Y})$.  

(A) If the complex-valued increasing function $f$ satisfies conditions given by Theorem~\ref{thm: prec w monotone f(x) cond inc f}, we have $f(\bm{X}) \prec_{\mbox{\tiny SN}} f(\bm{Y})$.  

(B) If the complex-valued decreasing function $f$ satisfies conditions given by Theorem~\ref{thm: prec w monotone f(x) cond dec f}, we have $f(\bm{X}) \prec_{\mbox{\tiny SN}} f(\bm{Y})$.  

(C) If the complex-valued function $f$ satsifies the following:
\begin{eqnarray}\label{eq5-1: thm: monotonicity conditions}
\lefteqn{[ [f(\lambda_{x,\sigma_x(1)})]_{\alpha_{x,\sigma_x(1)}^{\mathrm{A}}},\ldots, [f(\lambda_{x,\sigma_x(K_x)})]_{\alpha_{x,\sigma_x(K_x)}^{\mathrm{A}}}]=}\nonumber \\
&& [ [f(\lambda_{y,\sigma_y(1)})]_{\alpha_{y,\sigma_y(1)}^{\mathrm{A}}},\ldots, [f(\lambda_{y,\sigma_y(K_y)})]_{\alpha_{y,\sigma_y(K_y)}^{\mathrm{A}}}],
\end{eqnarray}
with $K_x = K_y$, $\alpha_{x,\sigma_x(k)}^{\mathrm{A}}=\alpha_{y,\sigma_y(k)}^{\mathrm{A}}$, and $[1]_{1 \leq i \leq \alpha_{y,k}^{\mathrm{G}}} \triangleleft [m_{y,k,i}]_{1 \leq i \leq \alpha_{y,k}^{\mathrm{G}}} $ for some $k=1,2,\ldots,K_y$. Moreover, $\kappa_{x,k} \geq \max\limits_{\substack{k=1,2,\ldots,K_x \\ i=1,2,\ldots,\alpha_{x,k}^{\mathrm{G}}}}m_{x,k,i}$  and  $\kappa_{y,k}=1$ for all $k=1,2,\ldots,K_y$, then, we have $f(\bm{X}) \prec_{\mbox{\tiny SN}} f(\bm{Y})$.

\textbf{Case II:  $[ [\lambda_{x,1}]_{\alpha_{x,1}^{\mathrm{A}}},\ldots, [\lambda_{x,K_{x}}]_{\alpha_{x,K_{x}}^{\mathrm{A}}}] = [ [\lambda_{y,1}]_{\alpha_{y,1}^{\mathrm{A}}},\ldots, [\lambda_{y,K_{y}}]_{\alpha_{y,K_{y}}^{\mathrm{A}}}]$ but \\
$[[m_{x,1,i}]_{1 \leq i \leq \alpha_{x,1}^{\mathrm{G}}},\ldots,[m_{x,K_{x},i}]_{1 \leq i \leq \alpha_{x,K_{x}}^{\mathrm{G}}}]\prec_{\mbox{\tiny N}}[[m_{y,1,i}]_{1 \leq i \leq \alpha_{y,1}^{\mathrm{G}}},\ldots,[m_{y,K_{y},i}]_{1 \leq i \leq \alpha_{y,K_{y}}^{\mathrm{G}}}]$}. For Case II, we have the following two  possibilities to have $f(\bm{X}) \prec_{\mbox{\tiny SN}} f(\bm{Y})$.  

(D) If the complex-valued increasing function $f$ satisfying $\kappa_{x, k} > 1$ for some $k \in \{1,2,\ldots,K_x\}$ and $\kappa_{y,k}=1$ for all $k=1,2,\ldots,K_y$, we have $f(\bm{X}) \prec_{\mbox{\tiny SN}} f(\bm{Y})$.  

(E) If we have $[m_{x,k,i}]_{1 \leq i \leq \alpha_{x,k}^{\mathrm{G}}} \triangleleft [m_{y,k,i}]_{1 \leq i \leq \alpha_{y,k}^{\mathrm{G}}}$ for all $k=1,2,\ldots,K_x$, and the complex-valued decreasing function $f$ satisfying $\kappa_{x,k} > 1$ for some $k \in \{1,2,\ldots,K_x\}$ and $\kappa_{y,k}=1$ for all $k=1,2,\ldots,K_y$~\footnote{Note that we have $\alpha_{x,k}^{\mathrm{A}}=\alpha_{y,k}^{\mathrm{A}}$ and $K_x=K_y$ here.}, we have $f(\bm{X}) \prec_{\mbox{\tiny SN}} f(\bm{Y})$.
\end{theorem}
\textbf{Proof:}
For (A), since the function $f$ satisifies those conidtions given by Theorem~\ref{thm: prec w monotone f(x) cond inc f}, we have
\begin{eqnarray}\label{eq6: thm: monotonicity conditions}
[ [f(\lambda_{x,1})]_{\alpha_{x,1}^{\mathrm{A}}},\ldots, [f(\lambda_{x,K_{x}})]_{\alpha_{x,K_{x}}^{\mathrm{A}}}] \prec_w [ [f(\lambda_{y,1})]_{\alpha_{y,1}^{\mathrm{A}}},\ldots, [f(\lambda_{y,K_{y}})]_{\alpha_{y,K_{y}}^{\mathrm{A}}}],
\end{eqnarray}
The statement $f(\bm{X}) \prec_{\mbox{\tiny SN}} f(\bm{Y})$ is true by applying Definitoin~\ref{def: SNO} with Eq.~\eqref{eq6: thm: monotonicity conditions}.

For (B), since the function $f$ satisifies those conidtions given by Theorem~\ref{thm: prec w monotone f(x) cond dec f}, we have
\begin{eqnarray}\label{eq7: thm: monotonicity conditions}
[ [f(\lambda_{x,K_{x}})]_{\alpha_{x,K_{x}}^{\mathrm{A}}},\ldots,[f(\lambda_{x,1})]_{\alpha_{x,1}^{\mathrm{A}}}] \prec_w [ [f(\lambda_{y,K_{y}})]_{\alpha_{y,K_{y}}^{\mathrm{A}}},\ldots,[f(\lambda_{y,1})]_{\alpha_{y,1}^{\mathrm{A}}}],
\end{eqnarray}
The statement $f(\bm{X}) \prec_{\mbox{\tiny SN}} f(\bm{Y})$ is true by applying Definitoin~\ref{def: SNO} with Eq.~\eqref{eq6: thm: monotonicity conditions}.

For (C), since the function $f$ satisfies
\begin{eqnarray}\label{eq8: thm: monotonicity conditions}
\lefteqn{[ [f(\lambda_{x,\sigma_x(1)})]_{\alpha_{x,\sigma_x(1)}^{\mathrm{A}}},\ldots, [f(\lambda_{x,\sigma_x(K_x)})]_{\alpha_{x,\sigma_x(K_x)}^{\mathrm{A}}}]=}\nonumber \\
&& [ [f(\lambda_{y,\sigma_y(1)})]_{\alpha_{y,\sigma_y(1)}^{\mathrm{A}}},\ldots, [f(\lambda_{y,\sigma_y(K_y)})]_{\alpha_{y,\sigma_y(K_y)}^{\mathrm{A}}}],
\end{eqnarray}
with $K_x = K_y$, $\alpha_{x,\sigma_x(k)}^{\mathrm{A}}=\alpha_{y,\sigma_y(k)}^{\mathrm{A}}$, we need to verify that other conditions of $f$ will make the nilpotent parts satisfying $[\bm{\eta}_{f,\lambda_{x,\sigma_x(1)}},\ldots,\bm{\eta}_{f,\lambda_{x,\sigma_x(K_x)}}] \prec_{\mbox{\tiny N}} [\bm{\eta}_{f,\lambda_{y,\sigma_y(1)}},\ldots,\bm{\eta}_{f,\lambda_{y,\sigma_y(K_y)}}]$. 

We have
\begin{eqnarray}\label{eq9: thm: monotonicity conditions}
[\bm{\eta}_{f,\lambda_{x,\sigma_x(1)}},\ldots,\bm{\eta}_{f,\lambda_{x,\sigma_x(K_x)}}]&=_1&[[\overbrace{1,\ldots,1}^{\mbox{$\alpha_{x,\sigma_x(1)}^{\mathrm{A}}$ ones}}],\ldots,
[\overbrace{1,\ldots,1}^{\mbox{$\alpha_{x,\sigma_x(K_x)}^{\mathrm{A}}$ones}}]]\nonumber \\
&\prec_{\mbox{\tiny N}, 2}&[[m_{y,\sigma_y(1),i}]_{1 \leq i \leq \alpha_{y,\sigma_y(1)}^{\mathrm{G}}},\ldots,[m_{y,\sigma_y(K_{y}),i}]_{1 \leq i \leq \alpha_{y,\sigma_y(K_{y})}^{\mathrm{G}}}]\nonumber \\
&=_3&[\bm{\eta}_{f,\lambda_{y,\sigma_y(1)}},\ldots,\bm{\eta}_{f,\lambda_{y,\sigma_y(K_y)}}],
\end{eqnarray}
where $=_1$ comes from $\kappa_{x,k} \geq \max\limits_{\substack{k=1,2,\ldots,K_x \\ i=1,2,\ldots,\alpha_{x,k}^{\mathrm{G}}}}m_{x,k,i}$ and Corollary~\ref{cor: kappa larger or equal m}, $\preceq_{\mbox{\tiny N}, 2}$ comes from $[1]_{1 \leq i \leq \alpha_{y,k}^{\mathrm{G}}} \triangleleft [m_{y,k,i}]_{1 \leq i \leq \alpha_{y,k}^{\mathrm{G}}}$ for some $k=1,2,\ldots,K_y$, and $=_3$ comes from $\kappa_{y,k}=1$ for all $k=1,2,\ldots,K_y$, Therefore, we have $f(\bm{X}) \prec_{\mbox{\tiny SN}} f(\bm{Y})$ from Definitoin~\ref{def: SNO}.

For (D), since $[ [\lambda_{x,1}]_{\alpha_{x,1}^{\mathrm{A}}},\ldots, [\lambda_{x,K_{x}}]_{\alpha_{x,K_{x}}^{\mathrm{A}}}] = [ [\lambda_{y,1}]_{\alpha_{y,1}^{\mathrm{A}}},\ldots, [\lambda_{y,K_{y}}]_{\alpha_{y,K_{y}}^{\mathrm{A}}}]$, for the complex-valued increasing function $f$, we must have $[ [f(\lambda_{x,1})]_{\alpha_{x,1}^{\mathrm{A}}},\ldots, [f(\lambda_{x,K_{x}})]_{\alpha_{x,K_{x}}^{\mathrm{A}}}] = [[f(\lambda_{y,1})]_{\alpha_{y,1}^{\mathrm{A}}},\ldots, [f(\lambda_{y,K_{y}})]_{\alpha_{y,K_{y}}^{\mathrm{A}}}]$. 

We have
\begin{eqnarray}\label{eq10: thm: monotonicity conditions}
[\bm{\eta}_{f,\lambda_{x,1}},\ldots,\bm{\eta}_{f,\lambda_{x,K_x}}]
&\preceq_{\mbox{\tiny N}, 1}&[[m_{x,1,i}]_{1 \leq i \leq \alpha_{x,1}^{\mathrm{G}}},\ldots,[m_{x,K_{x},i}]_{1 \leq i \leq \alpha_{x,K_{x}}^{\mathrm{G}}}]\nonumber \\
&\prec_{\mbox{\tiny N}, 2}&[[m_{y,1,i}]_{1 \leq i \leq \alpha_{y,1}^{\mathrm{G}}},\ldots,[m_{y,K_{y},i}]_{1 \leq i \leq \alpha_{y,K_{y}}^{\mathrm{G}}}]\nonumber \\
&=_3&[\bm{\eta}_{f,\lambda_{y,1}},\ldots,\bm{\eta}_{f,\lambda_{y,K_y}}],
\end{eqnarray}
where $\preceq_{\mbox{\tiny N}, 1}$ comes from $\kappa_{x_k} > 1$ for some $k \in \{1,2,\ldots,K_x\}$, $\prec_{\mbox{\tiny N}, 2}$ comes from Casse II condition, and $=_3$ comes from $\kappa_{y,k}=1$ for all $k=1,2,\ldots,K_y$. Therefore, we have $f(\bm{X}) \prec_{\mbox{\tiny SN}} f(\bm{Y})$ from Definitoin~\ref{def: SNO}.

For (E), similar to the situation of (D), we also have \\$[[f(\lambda_{x,K_{x}})]_{\alpha_{x,K_{x}}^{\mathrm{A}}},\ldots,[f(\lambda_{x,1})]_{\alpha_{x,1}^{\mathrm{A}}}] = [[f(\lambda_{y,K_{y}})]_{\alpha_{y,K_{y}}^{\mathrm{A}}},\ldots,[f(\lambda_{y,1})]_{\alpha_{y,1}^{\mathrm{A}}}]$. 

We have
\begin{eqnarray}\label{eq11: thm: monotonicity conditions}
[\bm{\eta}_{f,\lambda_{x,K_x}},\ldots,\bm{\eta}_{f,\lambda_{x,1}}]
&\preceq_{\mbox{\tiny N}, 1}&[[m_{x,K_{x},i}]_{1 \leq i \leq \alpha_{x,K_{x}}^{\mathrm{G}}},\ldots,[m_{x,1,i}]_{1 \leq i \leq \alpha_{x,1}^{\mathrm{G}}}]\nonumber \\
&\prec_{\mbox{\tiny N}, 2}&[[m_{y,K_{y},i}]_{1 \leq i \leq \alpha_{y,K_{y}}^{\mathrm{G}}},\ldots,[m_{y,1,i}]_{1 \leq i \leq \alpha_{y,1}^{\mathrm{G}}}]\nonumber \\
&=_3&[\bm{\eta}_{f,\lambda_{y,K_y}},\ldots,\bm{\eta}_{f,\lambda_{y,1}}],
\end{eqnarray}
where $\preceq_{\mbox{\tiny N}, 1}$ comes from $\kappa_{x,k} > 1$ for some $k \in \{1,2,\ldots,K_x\}$, $\prec_{\mbox{\tiny N}, 2}$ comes from $[m_{x,k,i}]_{1 \leq i \leq \alpha_{x,k}^{\mathrm{G}}} \trianglelefteq [m_{y,k,i}]_{1 \leq i \leq \alpha_{y,k}^{\mathrm{G}}}$ for all $k=1,2,\ldots,K_x$, and $=_3$ comes from $\kappa_{y,k}=1$ for all $k=1,2,\ldots,K_y$. Therefore, we have $f(\bm{X}) \prec_{\mbox{\tiny SN}} f(\bm{Y})$ from Definitoin~\ref{def: SNO}.
$\hfill\Box$

\begin{remark}\label{rmk2}
\begin{enumerate}
\item For Case II, it is impossible to have $[[f(\lambda_{x,\sigma_x(1)})]_{\alpha_{x,\sigma_x(1)}^{\mathrm{A}}},\ldots, [f(\lambda_{x,\sigma_x(K_x)})]_{\alpha_{x,\sigma_x(K_x)}^{\mathrm{A}}}] \preceq_w 
[ [f(\lambda_{y,\sigma_y(1)})]_{\alpha_{y,\sigma_y(1)}^{\mathrm{A}}},\ldots, [f(\lambda_{y,\sigma_y(K_y)})]_{\alpha_{y,\sigma_y(K_y)}^{\mathrm{A}}}]$ since the complex-valued function $f$ cannot be mapped to distinct values for the same inputs. 
\item From Theorem~\ref{thm: monotonicity conditions} statement (A) and (D), if the complex-valued function satisfies conditions given by Theorem~\ref{thm: prec w monotone f(x) cond inc f} with $\kappa_{x, k} > 1$ for some $k \in \{1,2,\ldots,K_x\}$ and $\kappa_{y,k}=1$ for all $k=1,2,\ldots,K_y$, we always have $f(\bm{X}) \prec_{\mbox{\tiny SN}} f(\bm{Y})$.  
\item From Theorem~\ref{thm: monotonicity conditions} statement (B) and (E), if we have $[m_{x,k,i}]_{1 \leq i \leq \alpha_{x,k}^{\mathrm{G}}} \triangleleft [m_{y,k,i}]_{1 \leq i \leq \alpha_{y,k}^{\mathrm{G}}}$ for all $k=1,2,\ldots,K_x$, and the complex-valued function $f$ satisfying Theorem~\ref{thm: prec w monotone f(x) cond dec f} conditions, $\kappa_{x,k} > 1$ for some $k \in \{1,2,\ldots,K_x\}$ and $\kappa_{y,k}=1$ for all $k=1,2,\ldots,K_y$, we always have $f(\bm{X}) \prec_{\mbox{\tiny SN}} f(\bm{Y})$.
\end{enumerate}
\end{remark}

\subsection{Convexity}\label{sec: Convexity}

In this section, we define and examine matrix function convexity under the spectral and nilpotent structures ordering (SNO) introduced in Definition~\ref{def: SNO}.

Let \( f : \mathbb{C} \to \mathbb{C}\) be an analytic function. The function \( f \) is called matrix convex (or operator convex) under SNO if for any matrices \(\bm{A}, \bm{B}\) of the same size, and any \( t \in [0, 1] \):
\begin{eqnarray}\label{eq: convexity by SNO}
f(t\bm{A} + (1-t)\bm{B}) \preceq_{\mbox{\tiny SN}} t f(\bm{A}) + (1-t) f(\bm{B}).
\end{eqnarray}
This means the image of the convex combination of \( \bm{A} \) and \( \bm{B} \) under \( f \) is less than or equal to the convex combination of their images, in the spectral and nilpotent structures ordering. 

This section aims to characterize operator convexity using the approach introduced by Hansen and Pedersen in~\cite{hansen2003jensen,hiai2010matrix}, which employs \( 2 \times 2 \) block matrices. However, we extend their method to general matrices, rather than restricting to Hermitian matrices as required in~\cite{hansen2003jensen,hiai2010matrix}. The proof presented below relies on the following two lemmas.

\begin{lemma}\label{lma: 2.5.1-1}
Given an analytic function $f(z)$ within the domain for $|z| < R$, a matrix $\bm{X}$ with the dimension $m$ and $K$ distinct eigenvalues $\lambda_k$ for $k=1,2,\ldots,K$ such that the eigenvalues of $\bm{X}$ satisfies $\left\vert\lambda_k\right\vert<R$, we have
\begin{eqnarray}\label{eq1: lma: 2.5.1-1}
f(\bm{U}^{\mathrm{H}}\bm{X}\bm{U})&=&\bm{U}^{\mathrm{H}}f(\bm{X})\bm{U},
\end{eqnarray}
where $\mathrm{H}$ is the conjugate transpose operation, and $\bm{U}\in\mathbb{C}^{m \times m}$ is an unitary matrix.
\end{lemma}
\textbf{Proof:}
From Theorem~\ref{thm: Spectral Mapping Theorem for Single Variable}, we have
\begin{eqnarray}\label{eq2: lma: 2.5.1-1}
\bm{X}&=&\sum\limits_{k=1}^K\sum\limits_{i=1}^{\alpha_k^{\mathrm{G}}} \lambda_k \bm{P}_{k,i}+
\sum\limits_{k=1}^K\sum\limits_{i=1}^{\alpha_k^{\mathrm{G}}} \bm{N}_{k,i},
\end{eqnarray}
where $\left\vert\lambda_k\right\vert<R$, then, we have
\begin{eqnarray}\label{eq2: lma: 2.5.1-2}
\bm{U}^{\mathrm{H}}\bm{X}\bm{U}&=&\sum\limits_{k=1}^K\sum\limits_{i=1}^{\alpha_k^{\mathrm{G}}} \lambda_k \bm{U}^{\mathrm{H}}\bm{P}_{k,i}\bm{U}+
\sum\limits_{k=1}^K\sum\limits_{i=1}^{\alpha_k^{\mathrm{G}}} \bm{U}^{\mathrm{H}} \bm{N}_{k,i}\bm{U}.
\end{eqnarray}

Therefore, we obtain
\begin{eqnarray}\label{eq3: lma: 2.5.1-2}
f(\bm{U}^{\mathrm{H}}\bm{X}\bm{U})&=&\sum\limits_{k=1}^K \left[\sum\limits_{i=1}^{\alpha_k^{(\mathrm{G})}}f(\lambda_k)\bm{U}^{\mathrm{H}}\bm{P}_{k,i}\bm{U}+\sum\limits_{i=1}^{\alpha_k^{(\mathrm{G})}}\sum\limits_{q=1}^{m_{k,i}-1}\frac{f^{(q)}(\lambda_k)}{q!}\bm{U}^{\mathrm{H}}\bm{N}_{k,i}^q\bm{U}\right]\nonumber \\
&=&\bm{U}^{\mathrm{H}}f(\bm{X})\bm{U}
\end{eqnarray}
$\hfill\Box$

\begin{lemma}\label{lma: 2.5.1-2}
Given an analytic function $f(z)$ within the domain for $|z| < R$, a matrix $\bm{X}$ with the dimension $m$ and $K$ distinct eigenvalues $\lambda_k$ for $k=1,2,\ldots,K$ such that the eigenvalues of $\bm{X}^{\mathrm{H}}\bm{X}$ and $\bm{X}\bm{X}^{\mathrm{H}}$ within in $\mathbb{D} = \{ z \in \mathbb{C} : |z| < R \}$, we have
\begin{eqnarray}\label{eq1: lma: 2.5.1-2}
\bm{X}f(\bm{X}^{\mathrm{H}}\bm{X})&=&f(\bm{X}\bm{X}^{\mathrm{H}})\bm{X}.
\end{eqnarray}
\end{lemma}
\textbf{Proof:}
For any positive integer $k$, we have 
\begin{eqnarray}\label{eq2: lma: 2.5.1-2}
\bm{X}(\bm{X}^{\mathrm{H}}\bm{X})^k &=& (\bm{X}\bm{X}^{\mathrm{H}})^k\bm{X}.
\end{eqnarray}
Hence, we have Eq.~\eqref{eq1: lma: 2.5.1-2} if $f$ is a polynomial. Let $\lambda'_k$ be eigenvalues of $\bm{X}^{\mathrm{H}}\bm{X}$, for any analytic function $f(z)$, we can construct polynomial $p(z)$ via Hermite interpolation method to achieve the following conditions: 
\begin{eqnarray}\label{eq3: lma: 2.5.1-2}
f(\lambda'_k)=p(\lambda'_k),\quad,  f^{(1)}(\lambda'_k)=p^{(1)}(\lambda'_k), \quad, \ldots, \quad f^{(m_{k,i})}(\lambda'_k)=p^{(m_{k,i})}(\lambda'_k),
\end{eqnarray}
where $k$ is the index for the distinct eigenvalues of $\bm{X}^{\mathrm{H}}\bm{X}$, $i$ is the index for Jordan blocks with eigenvalue as $\lambda'_k$, and $m_{k,i}$ are positive integers corresponding to the nilpotent degree of the nilpotent $\bm{N}_{k,i}$, i.e., $m_{k,i}$ is the minimum integer of $\ell$ to have $\bm{N}^\ell_{k,i} = \bm{O}$, of the matrix $\bm{X}^{\mathrm{H}}\bm{X}$.

Therefore, from the spectral mapping theorem~\ref{thm: Spectral Mapping Theorem for Single Variable}, we have
\begin{eqnarray}\label{eq4: lma: 2.5.1-2}
\bm{X}f(\bm{X}^{\mathrm{H}}\bm{X})=\bm{X}p(\bm{X}^{\mathrm{H}}\bm{X})
=p(\bm{X}\bm{X}^{\mathrm{H}})\bm{X}=f(\bm{X}\bm{X}^{\mathrm{H}})\bm{X},
\end{eqnarray}
which is the desired result. 
$\hfill\Box$

\begin{remark}
Lemma~\ref{lma: 2.5.1-1} and Lemma~\ref{lma: 2.5.1-2} are extension from self-adjoint (Hermitian) matrices (Lemma 2.5.1 in~\cite{hiai2010matrix}) to arbitrary matrix via the spectral mapping theorem~\ref{thm: Spectral Mapping Theorem for Single Variable}.
\end{remark}

Below, we will provide the following Theorem~\ref{thm: HP convexity characterization} to characterize convexity based on Hansen and Pedersen's method.

\begin{theorem}\label{thm: HP convexity characterization}
Given an analytic function $f(z)$ within the domain for $|z| < R$, all eigenvalues of the matrices $\bm{X}, \bm{X}_i, \bm{C}, \bm{C}_i, \bm{Y}, \bm{C}^\mathrm{H}\bm{X}\bm{C},\sum\limits_{i=1}^{\ell}\bm{C}_i^\mathrm{H}\bm{X}_i\bm{C}_i, (\bm{I}-\bm{C}\bm{C}^\mathrm{H})^{1/2}\bm{X}(\bm{I}-\bm{C}\bm{C}^\mathrm{H})^{1/2},\bm{P}\bm{X}\bm{P}+(\bm{I}-\bm{P})\bm{Y}(\bm{I}-\bm{P})$, where $i \in \{1,2,\ldots,\ell\}$ and $\bm{P}$ is any orthogonal matrix, are within in $\mathbb{D} = \{ z \in \mathbb{C} : |z| < R \}$, we have the following equivalence:
\begin{enumerate}
\item $f$ is operator convex under SNO;
\item for every matrix $\bm{C}$ with spectral norm less or equal than one, i.e., $\left\Vert\bm{C}^\mathrm{H}\bm{C}\right\Vert_{s} \leq 1$, we have 
\begin{eqnarray}\label{eq1: thm: HP convexity characterization}
f(\bm{C}^{\mathrm{H}}\bm{X}\bm{C})\preceq_{\mbox{\tiny SN}}\bm{C}^{\mathrm{H}}f(\bm{X})\bm{C};
\end{eqnarray}
\item for every $\ell \in \mathbb{N}$ and every matrix $\bm{C}_i$ satisfying $\left\Vert\sum\limits_{i=1}^\ell \bm{C}_i^{\mathrm{H}}\bm{C}_i\right\Vert_s \leq 1$, we have
\begin{eqnarray}\label{eq2: thm: HP convexity characterization}
f\left(\sum\limits_{i=1}^{\ell}\bm{C}_i^{\mathrm{H}}\bm{X}_i\bm{C}_i\right)\preceq_{\mbox{\tiny SN}}\sum\limits_{i=1}^{\ell}\bm{C}_i^{\mathrm{H}}f(\bm{X}_i)\bm{C}_i;
\end{eqnarray}
\item for every orthogonal projection matrix $\bm{P}$, we have
\begin{eqnarray}\label{eq3: thm: HP convexity characterization}
f(\bm{P}\bm{X}\bm{P} + (\bm{I}-\bm{P})\bm{Y} (\bm{I}-\bm{P}))\preceq_{\mbox{\tiny SN}}
\bm{P}f(\bm{X})\bm{P} + (\bm{I}-\bm{P})f(\bm{Y})(\bm{I}-\bm{P}).
\end{eqnarray}
\end{enumerate}
\end{theorem}
\textbf{Proof:}
Without loss of generality, we will assume $f(0)=0$. If not, we can add some constant to the function $f$ to make $f(0)+c=0$, and such vertical shifting will not affect its convexity.

\textbf{Item 1 $\Rightarrow$ Item 2}

Let us construct the following matrices for later proof:
\begin{eqnarray}\label{eq4: thm: HP convexity characterization}
\overline{\bm{X}}&\define& \begin{pmatrix}
   \bm{X} & \bm{O} \\
   \bm{O} & \bm{O}
   \end{pmatrix}, \nonumber \\
\bm{W}_1&\define& \begin{pmatrix}
   \bm{C} & (\bm{I}-\bm{C}\bm{C}^{\mathrm{H}})^{1/2} \\
   (\bm{I}-\bm{C}^{\mathrm{H}}\bm{C})^{1/2} & -\bm{C}^{\mathrm{H}}
   \end{pmatrix}, \nonumber \\
\bm{W}_2&\define& \begin{pmatrix}
   \bm{C} & -(\bm{I}-\bm{C}\bm{C}^{\mathrm{H}})^{1/2} \\
   (\bm{I}-\bm{C}^{\mathrm{H}}\bm{C})^{1/2} & \bm{C}^{\mathrm{H}}
   \end{pmatrix}.
\end{eqnarray}
From these new matrices $\overline{\bm{X}}, \bm{W}_1$ and $\bm{W}_2$, we have the following facts. \\
(1) By applying Lemma~\ref{lma: 2.5.1-2} with $\left\Vert\bm{C}\right\Vert_{s} \leq 1$  to the function $(1-z)^{1/2}$, we have $\bm{C}(\bm{I}-\bm{C}^{\mathrm{H}}\bm{C})^{1/2}=(\bm{I}-\bm{C}\bm{C}^{\mathrm{H}})^{1/2}\bm{C}$. \\
(2) From above fact (1), we have $\bm{W}_1^{\mathrm{H}}\bm{W}_1 = \bm{W}_2^{\mathrm{H}}\bm{W}_2 = \bm{I}$.\\ 
(3) The matrix $\bm{W}_1^{\mathrm{H}}\overline{\bm{X}}\bm{W}_1$ can be expressed as
\begin{eqnarray}\label{eq5: thm: HP convexity characterization}
\bm{W}_1^{\mathrm{H}}\overline{\bm{X}}\bm{W}_1&=& \begin{bmatrix}
   \bm{C}^{\mathrm{H}}\bm{X}\bm{C} & \bm{C}^{\mathrm{H}}\bm{X}(\bm{I}-\bm{C}\bm{C}^{\mathrm{H}})^{1/2} \\
   (\bm{I}-\bm{C}\bm{C}^{\mathrm{H}})^{1/2}\bm{X}\bm{C} & (\bm{I}-\bm{C}\bm{C}^{\mathrm{H}})^{1/2}\bm{X}(\bm{I}-\bm{C}\bm{C}^{\mathrm{H}})^{1/2}
   \end{bmatrix}.
\end{eqnarray}
\\
(4) The matrix $\bm{W}_2^{\mathrm{H}}\overline{\bm{X}}\bm{W}_2$ can be expressed as
\begin{eqnarray}\label{eq6: thm: HP convexity characterization}
\bm{W}_2^{\mathrm{H}}\overline{\bm{X}}\bm{W}_2&=& \begin{bmatrix}
   \bm{C}^{\mathrm{H}}\bm{X}\bm{C} & -\bm{C}^{\mathrm{H}}\bm{X}(\bm{I}-\bm{C}\bm{C}^{\mathrm{H}})^{1/2} \\
   -(\bm{I}-\bm{C}\bm{C}^{\mathrm{H}})^{1/2}\bm{X}\bm{C} & (\bm{I}-\bm{C}\bm{C}^{\mathrm{H}})^{1/2}\bm{X}(\bm{I}-\bm{C}\bm{C}^{\mathrm{H}})^{1/2}
   \end{bmatrix}.
\end{eqnarray}

From these facts, we have
\begin{eqnarray}\label{eq7: thm: HP convexity characterization}
\lefteqn{\begin{bmatrix}
   f(\bm{C}^{\mathrm{H}}\bm{X}\bm{C}) & \bm{O} \\
  \bm{O} & f((\bm{I}-\bm{C}\bm{C}^{\mathrm{H}})^{1/2}\bm{X}(\bm{I}-\bm{C}\bm{C}^{\mathrm{H}})^{1/2})
\end{bmatrix}}\nonumber \\
&=&f\left(\begin{bmatrix}
   \bm{C}^{\mathrm{H}}\bm{X}\bm{C} & \bm{O} \\
  \bm{O} & (\bm{I}-\bm{C}\bm{C}^{\mathrm{H}})^{1/2}\bm{X}(\bm{I}-\bm{C}\bm{C}^{\mathrm{H}})^{1/2}
\end{bmatrix}\right)\nonumber \\
&=_1&f\left(\frac{\bm{W}_1^{\mathrm{H}}\overline{\bm{X}}\bm{W}_1+\bm{W}_2^{\mathrm{H}}\overline{\bm{X}}\bm{W}_2}{2}\right)\nonumber \\
&\preceq_{\mbox{\tiny SN}, 2}&\frac{f(\bm{W}_1^{\mathrm{H}}\overline{\bm{X}}\bm{W}_1)}{2}+\frac{f(\bm{W}_2^{\mathrm{H}}\overline{\bm{X}}\bm{W}_2)}{2}\nonumber \\
&=_3&\frac{1}{2}\bm{W}_1^{\mathrm{H}}\begin{bmatrix}
   f(\bm{X}) & \bm{O} \\
  \bm{O} & \bm{O}
\end{bmatrix}\bm{W}_1 + \frac{1}{2}\bm{W}_2^{\mathrm{H}}\begin{bmatrix}
   f(\bm{X}) & \bm{O} \\
  \bm{O} & \bm{O}
\end{bmatrix}\bm{W}_2 \nonumber \\
&=& \begin{bmatrix}
   \bm{C}^{\mathrm{H}}f(\bm{X})\bm{C} & \bm{O} \\
  \bm{O} & (\bm{I}-\bm{C}\bm{C}^{\mathrm{H}})^{1/2}f(\bm{X})(\bm{I}-\bm{C}\bm{C}^{\mathrm{H}})^{1/2}
\end{bmatrix},
\end{eqnarray}
where we apply above facts (3) and (4) in $=_1$, $\preceq_{\mbox{\tiny SN}, 2}$ comes from the convexity assumption of the function $f$, and $=_3$ comes from Lemma~\ref{lma: 2.5.1-1} with fact (2) about unitary matrices $\bm{W}_1$ and $\bm{W}_2$. Therefore, we must require $f(\bm{C}^{\mathrm{H}}\bm{X}\bm{C})\preceq_{\mbox{\tiny SN}}\bm{C}^{\mathrm{H}}f(\bm{X})\bm{C}$ to make the inequality $\preceq_{\mbox{\tiny SN}, 2}$ valid by comparing the upper-left block of the matrix $\begin{bmatrix}
   f(\bm{C}^{\mathrm{H}}\bm{X}\bm{C}) & \bm{O} \\
  \bm{O} & f((\bm{I}-\bm{C}\bm{C}^{\mathrm{H}})^{1/2}\bm{X}(\bm{I}-\bm{C}\bm{C}^{\mathrm{H}})^{1/2})
\end{bmatrix}$ and the matrix \\
 $\begin{bmatrix}
   \bm{C}^{\mathrm{H}}f(\bm{X})\bm{C} & \bm{O} \\
  \bm{O} & (\bm{I}-\bm{C}\bm{C}^{\mathrm{H}})^{1/2}f(\bm{X})(\bm{I}-\bm{C}\bm{C}^{\mathrm{H}})^{1/2}
\end{bmatrix}$.

\textbf{Item 2 $\Rightarrow$ Item 3}

Given $\ell$ matrices $\bm{X}_i$ and $\ell$ matrices $\bm{C}_i$, we construct matrices $\overline{\overline{\bm{X}}}$ and $\overline{\bm{C}}$ as follows:
\begin{eqnarray}\label{eq7: thm: HP convexity characterization}
\overline{\overline{\bm{X}}}&\define& \begin{bmatrix}
\bm{X}_1 & \bm{O}    &  \ldots & \bm{O} \\
\bm{O}    & \bm{X}_2 &  \ldots & \bm{O} \\
\vdots    & \vdots &  \ddots & \vdots  \\
\bm{O}   & \bm{O} & \ldots & \bm{X}_\ell \\
\end{bmatrix},
\end{eqnarray}
and
\begin{eqnarray}\label{eq8: thm: HP convexity characterization}
\overline{\bm{C}}&\define& \begin{bmatrix}
\bm{C}_1 \\
\bm{C}_2 \\
\vdots  \\
\bm{C}_\ell \\
\end{bmatrix}.
\end{eqnarray}
From $\left\Vert\sum\limits_{i=1}^\ell \bm{C}_i^{\mathrm{H}}\bm{C}_i\right\Vert_s \leq 1$, we have
\begin{eqnarray}\label{eq9: thm: HP convexity characterization}
f\left(\sum\limits_{i=1}^{\ell}\bm{C}_i^{\mathrm{H}}\bm{X}_i\bm{C}_i\right)&=&f(\overline{\bm{C}}^{\mathrm{H}}\overline{\overline{\bm{X}}}\overline{\bm{C}})\nonumber \\
&\preceq_{\mbox{\tiny SN}}&
\overline{\bm{C}}^{\mathrm{H}}f(\overline{\overline{\bm{X}}})\overline{\bm{C}}\nonumber \\
&=&\sum\limits_{i=1}^{\ell}\bm{C}_i^{\mathrm{H}}f(\bm{X}_i)\bm{C}_i
\end{eqnarray}

\textbf{Item 3 $\Rightarrow$ Item 4}

If we set $\bm{C}_1 = \bm{P}$, and $\bm{C}_2 = \bm{I} - \bm{P}$, Item 4 is obtained by applying Item 3 with $\ell=2$. \\

Finally, we will show \textbf{Item 4 $\Rightarrow$ Item 1}.\\

Let us construct the following matrices for later proof:
\begin{eqnarray}\label{eq10: thm: HP convexity characterization}
\acute{\bm{X}}&\define& \begin{pmatrix}
   \bm{X} & \bm{O} \\
   \bm{O} & \bm{Y}
   \end{pmatrix}, \nonumber \\
\bm{W}_3&\define& \begin{pmatrix}
   \sqrt{t} \bm{I} & - \sqrt{1 - t} \bm{I} \\
   \sqrt{1 - t} \bm{I} & \sqrt{t} \bm{I}
   \end{pmatrix}, \nonumber \\
\bm{P}&\define& \begin{pmatrix}
   \bm{I}   & \bm{O} \\
   \bm{O} & \bm{O} 
   \end{pmatrix}.
\end{eqnarray}
where we assume that $0 < t < 1$. Note that the matrix $\bm{W}_3$ is a unitary matrix. Then, we have
\begin{eqnarray}\label{eq11: thm: HP convexity characterization}
\lefteqn{\begin{bmatrix}
   f(t\bm{X}+(1-t)\bm{Y}) & \bm{O} \\
   \bm{O} & f((1-t)\bm{X}+t\bm{Y}) 
   \end{bmatrix}}\nonumber \\
&=&f(\bm{P}\bm{W}_3^{\mathrm{H}}\acute{\bm{X}}\bm{W}_3\bm{P} + (\bm{I}-\bm{P})\bm{W}_3^{\mathrm{H}}\acute{\bm{X}}\bm{W}_3(\bm{I}-\bm{P}))\nonumber \\
&\preceq_{\mbox{\tiny SN},1}&\bm{P}f(\bm{W}_3^{\mathrm{H}}\acute{\bm{X}}\bm{W}_3)\bm{P} + (\bm{I}-\bm{P})f(\bm{W}_3^{\mathrm{H}}\acute{\bm{X}}\bm{W}_3)(\bm{I}-\bm{P}) \nonumber \\
&=_2&\bm{P}\bm{W}_3^{\mathrm{H}}\begin{bmatrix}
   f(\bm{X}) & \bm{O} \\
   \bm{O} & f(\bm{Y}) 
   \end{bmatrix} \bm{W}_3\bm{P}\nonumber \\
&&+ (\bm{I}-\bm{P})\bm{W}_3^{\mathrm{H}}\begin{bmatrix}
   f(\bm{X}) & \bm{O} \\
   \bm{O} & f(\bm{Y}) 
   \end{bmatrix} \bm{W}_3(\bm{I}-\bm{P})\nonumber \\
&=&\begin{bmatrix}
   tf(\bm{X})+(1-t)f(\bm{Y}) & \bm{O} \\
   \bm{O} & (1-t)f(\bm{X})+tf(\bm{Y}) 
   \end{bmatrix},
\end{eqnarray}
where we apply Item 4 in $\preceq_{\mbox{\tiny SN}, 1}$ and Lemma~\ref{lma: 2.5.1-1} with the unitary matrix $\bm{W}_3$ to $=_2$. Therefore, we must require $f(t\bm{X}+(1-t)\bm{Y})\preceq_{\mbox{\tiny SN}} tf(\bm{X}+(1-t)f(\bm{Y})$ to make the inequality $\preceq_{\mbox{\tiny SN}, 1}$ valid by comparing the upper-left block of the matrix $\begin{bmatrix}
   f(t\bm{X}+(1-t)\bm{Y}) & \bm{O} \\
   \bm{O} & f((1-t)\bm{X}+t\bm{Y}) 
   \end{bmatrix}$ and the matrix \\
 $\begin{bmatrix}
   tf(\bm{X})+(1-t)f(\bm{Y}) & \bm{O} \\
   \bm{O} & (1-t)f(\bm{X})+tf(\bm{Y}) 
\end{bmatrix}$.
$\hfill\Box$

\begin{remark}
The characterization method builds on the work in~\cite{hansen2003jensen,hiai2010matrix}, particularly Theorems 2.5.2 and 2.5.7 from~\cite{hiai2010matrix}. However, Theorem~\ref{thm: HP convexity characterization} is developed within a significantly more general ordering framework, SNO.
\end{remark}

\bibliographystyle{IEEETran}
\bibliography{ComplexNumberCompare_MatrixArg_Bib}

\end{document}